\documentclass[11pt,lettersize,reqno,fullpage]{article}
\usepackage{amsmath,amssymb,amsfonts,mathrsfs,verbatim,enumitem,pstricks,amsthm, graphicx}
\usepackage[all]{xy}
\usepackage{centernot, xr, accents}
\usepackage{hyperref}

\title{{\LARGE Counting curves in a linear system with upto eight singular points}}
\author{Somnath Basu and Ritwik Mukherjee }
\date{}
\usepackage{hyperref}
\hypersetup{
	colorlinks,
	citecolor=red,
	filecolor=black,
	linkcolor=blue,
	urlcolor=black
}

\theoremstyle{plain}
\newtheorem{thm}{Theorem}[section]

\newtheorem{lmm}[thm]{Lemma}

\newtheorem{que}[thm]{Question}

\newtheorem{mresult}[thm]{Main Result}

\theoremstyle{definition}
\newtheorem{rem}[thm]{Remark} 
\newtheorem{defn}[thm]{Definition}

\def \hf{\hspace*{0.5cm}}                      % Somnath
\def\bge{\begin{equation}}                % Somnath
\def\ede{\end{equation}}                % Somnath
\def\bgd{\begin{displaymath}}         % Somnath
\def\edd{\end{displaymath}}            % Somnath
\def\bgee{\begin{equation*}}           % Somnath
\def\edee{\end{equation*}}           % Somnath

\def \ni{\noindent}

\def\lra{\longrightarrow}

\def\lan{\langle}
\def\ran{\rangle}

\def\BA{\begin{eqnarray}}
\def\EA{\end{eqnarray}}
\def\BAA{\begin{eqnarray*}}
\def\EAA{\end{eqnarray*}}
\def\Bal{\begin{align*}}
\def\Eal{\end{align*}}

\def \C{\mathbb{C}}

\def \P{\mathbb{P}}
\def \A{\mathcal{A}}
\def \E{\mathcal{E}}
\def \Num{N}
\def \PP{\mathcal{P}}
\def \f{\frac}

\def \B{\mathcal{B}}

\def \25node{A_5} 
\def \62node{A_6}

\def \D{\mathcal{D}}
\def \DD{\mathcal{D}}

\def \X{\mathfrak{X}}

\def \ov{\overline}

\def \DD{\hat{\D}}

\def \q{f}

\def \U{\mathcal{U}}

\def \N{\nabla}

\def \sq{s}
\def \rq{r}

\def \hxt2{\hat{x}_{t_2}} 
\def \hyt2{\hat{y}_{t_2}}

\def \hf{\hspace*{0.5cm}}                      % Somnath
\def\bge{\begin{equation}}                % Somnath
\def\ede{\end{equation}}                % Somnath
\def\bgd{\begin{displaymath}}         % Somnath
\def\edd{\end{displaymath}}            % Somnath
\def\bgee{\begin{equation*}}           % Somnath
\def\edee{\end{equation*}}           % Somnath

\def \P{\mathbb{P}}
\def \PP{\mathcal{P}}
\def \X{X} 
\def \A{A}
\def \D{D}
\def \E{E}
\def \N{N}
\def \lra{\longrightarrow}
\def \U{\mathcal{U}}
\def \C{\mathbb{C}}
\def \ni{\noindent}
\def \DD{\mathcal{D}}
\def \lan{\langle} 
\def \ran{\rangle} 
\def \B{\mathcal{B}}

\begin{document}

\maketitle

\begin{abstract} 
In this paper, we develop a systematic approach to enumerate 
curves with a certain number of nodes and one further singularity which maybe 
more degenerate. 
As a result, we obtain an explicit formula for the number of curves
in a sufficiently ample linear system, passing through the right number of generic points,
that have $\delta$ nodes and one singularity of codimension $k$, for all  
$\delta+k \leq 8$. 
In particular, we recover the formulas for curves with upto six nodal points obtained by 
Vainsencher. Moreover, all the codimension seven numbers we have obtained agree with the formulas 
obtained by Kazarian. Finally, in codimension eight, we recover the formula 
of A.Weber, M.Mikosz and P.Pragacz 
for 
curves with one singular point  and we also recover the formula of 
%eight nodal curves obtained  by 
Kleiman and Piene for eight nodal curves. 
All the other codimension eight numbers we have obtained are new. 
%(to the best of our knowledge). 
%(with the exception of the formula 
%for eight nodal curves 
%and the formula 
%for curves with one singular point which have also been computed by 
%A.Weber, M.Mikosz and P.Pragacz).
%, namely curves with  
%a singularity of type $A_8, D_8, E_8$ and $X_9$). 
%. singularity). 
%Our starting point is a classical fact from 
%differential topology: 
%the number of zeros of a generic smooth section of a vector 
%bundle $V$ over $M$, counted with signs, is 
%the Euler class of $V$ evaluated on the fundamental class of $M$.
%We then go on to use a topological method to compute the degenerate contribution to 
%the Euler class. 
\end{abstract}

\tableofcontents

\section{Introduction}
\hf\hf Enumerative geometry is a branch of mathematics concerned with the following question: 
\begin{center}
{\it How many geometric objects are there which satisfy prescribed constraints?  }
\end{center}
It is a classical subject that dates back to 
over $150$ years ago. 
%The general goal of this  subject is to count how many geometric objects 
%are there that satisfy certain conditions. 
%As early $1873$, Zeuthen calculated th number of quartics with upto three nodes.   
%It has been an active field 
%of research since then. 
%nineteenth century. In fact, 
Hilbert's fifteenth problem was to lay a rigorous 
foundation for enumerative Schubert calculus.   
While the problems in this
field are typically easy to state, solutions to almost
all of them require various deep concepts from various branches of mathematics.\\
%singular curves in $\P^2$ (complex projective space) is a question 
%that has been studied for more than hundred years. In \cite{Zeuthen} 
%problem in algebraic geometry. 
\hf \hf In order to motivate the results of this paper, let us first state a 
very general problem in enumerative geometry: 
%is as follows: 
\begin{que}
\label{main_question}
Let $L \lra \X$ be a holomorphic 
line bundle over a compact complex surface and 
$\mathcal{D} :=\P H^0(\X,L) \approx \P^{\delta_L}$ be the space of 
non zero holomorphic sections upto scaling. What is $\N(\A_1^{\delta} \mathfrak{X})$,
the number of curves in $\X$, that belong to the linear system 
$H^0(\X, L)$, passing through $\delta_L-(k+\delta)$ generic points 
and having $\delta$ distinct nodes and one singularity of type $\mathfrak{X}$ whose codimension is $k$.
%where $k$ is the codimension\footnote{By codimension we mean the number of 
%conditions that is imposed on the space of curves by having that particular singularity. 
%For example, a node is a codimension one singularity, a cusp is a codimesnion two singularity and a taconde is a codimension three 
%singularity.} of the singularity?
\end{que}
\hf \hf The above question has been studied by several mathematicians over the 
last thirty years; it becomes increasingly hard when the codimension $\delta+k$ 
increases and/or when the singularity $\mathfrak{X}$ becomes more degenerate. 
In this paper, we give develop a systematic method to approach this question 
and 
%In particular, we 
obtain 
a complete answer to the above question till codimension eight. 
Before stating the main result of this paper, let us make  a couple of definitions: 
%giving the formulas for $\N(\A_1^{\delta} \mathfrak{X})$, let us 
%make a few definitions. 

\begin{defn}
\label{singularity_defn}
Let $L \lra \X$ be a holomorphic 
line bundle over a complex surface and 
and $f \in H^0(\X, L)$ a holomorphic 
section. 
%Let $[f] \in \D$ and $q \in \X$. 
A point $q \in f^{-1}(0)$ \textsf{is of singularity type} $\A_k$,
$\D_k$, $\E_6$, $\E_7$, $\E_8$ or $\X_9$ if there exists a coordinate system
$(x,y) :(\U,q) \lra (\C^2,0)$ such that $f^{-1}(0) \cap \U$ is given by 
\begin{align*}
\A_k: y^2 + x^{k+1}   &=0  \qquad k \geq 0, \qquad \D_k: y^2 x + x^{k-1} =0  \qquad k \geq 4, \\
\E_6: y^3+x^4 &=0,  \qquad \E_7: y^3+ y x^3=0, 
\qquad \E_8: y^3 + x^5=0, \\
\X_9: x^4 + y^4 &=0.   
\end{align*}
\end{defn}
\hf \hf In more common terminology, $q$ is a {\it smooth} point of $f^{-1}(0)$ if 
it is a singularity of type $\A_0$; a {\it simple node} (or just node) if its singularity type is $\A_1$; 
a {\it cusp} if its type is $\A_2$; a {\it tacnode} if its type is $\A_3$; a {\it an ordinary triple point} 
if its type is $\D_4$; and an {\it ordinary quadruple point} if its type is $X_9$. \\

\begin{defn}
A holomorphic line bundle $L \lra \X$ over a compact complex manifold $\X$ 
is \textsf{sufficiently $k$-ample} if $L \approx L_1^{\otimes n} \otimes \xi \lra \X$ for some $n\geq k$, 
%is of the form 
%\[ L= L_1^{\otimes n} \otimes L_2 \] 
where $L_1 \lra \X$ is a very ample line bundle and $\xi \lra \X$ is a line bundle 
such that the linear system $H^0(\X, \xi)$ is base point free.
%and 
%the linear system $H^0(X, L_2)$ is base point free, 
%i.e. given any point $q\in X$, there exists a holomorphic 
%section $f \in H^0(\X, L_2)$, such that $f(q) \neq 0$.
\end{defn}

\begin{rem}
%Before stating the main result of this paper, let us set up a somewhat loose terminology. 
We will frequently use the phrase ``a singularity of codimension $k$''. Roughly speaking, 
this refers to the number of conditions having that singularity imposes on the space of curves. 
More precisely, it is the expected codimension of the equisingular strata.  
Hence, a singularity of type $A_k$, $D_k$ or $E_{k \leq 8}$ is a singularity of codimension $k$. However, a 
singularity of type $X_9$ is of codimension $8$, not $9$. The reason that it is standard 
to denote an ordinary quadruple point by $X_9$ and not $X_8$ is because singularities 
are traditionally indexed by their \textsf{Milnor number} and not by the codimension of the 
equisingular strata. For $A_k$, $D_k$ and $E_{k \leq 8}$ singularities, 
these two numbers happen to be the same. In 
general, these two numbers need not be the same. 
%For our purposes, the codimension of the equisingular 
%strata is going to play a more imp
%notation 
\end{rem}

The main result of this paper 
%(cf. Theorem \ref{one_pt_chern_class_codim_7} and 
%Theorem \ref{two_pt_chern_class_codim_7}) 
is as follows: 

\begin{mresult}
\label{main_codim_8}
%Let $\mathfrak{X}$ be a singularity of codimension $k$. 
%Let $\N(\A_1^{\delta} \mathfrak{X})$ denote the 
%If $\delta \leq 8$ and  
%$\delta+k \leq 7$, then 
Let $L \lra X$ be a line bundle 
%that is sufficiently ample 
and 
$\mathfrak{X}$ be one of the singularities defined in \ref{singularity_defn} whose codimension is $k$.  
%singularity. 
%\footnote{By codimension we mean the number of 
%conditions having that particular singularity imposes on the space of curves. 
%For example, a node is a codimension one 
%singularity, a cusp is a codimesnion two singularity, a taconde is a codimension three 
%singularity and so on.}
We obtain an explicit formula for 
$\N(\A_1^{\delta} \mathfrak{X})$ if $\delta +k \leq 8$, provided the line bundle is 
sufficiently $(2 \delta + C_{\mathfrak{X}})$-ample, where 
\[ C_{A_k} = k, ~~C_{D_4} = 3, ~~C_{D_k} = k-2 ~~\textnormal{if} ~~k>4, \qquad ~~C_{E_k} = k-3 \qquad 
\textnormal{and} \qquad C_{X_9} = 4.\]
%provided $L \lra \X$ is sufficiently ample.
\end{mresult}

\begin{rem}
The ampleness requirement is there to prove that the sections we encounter are 
transverse to the zero set. 
However, the bound 
we impose is not the optimal bound. In particular, our bound is coarser than the 
bound given by Kleiman and Piene in \cite{KP}. 
%to prove enumerativity of their formula for 
%eight nodal curves. 
One of the reasons our bound is coarse, is because while proving transversality, we restrict ourselves 
exclusively to the tangent space of the curves; we do not utilize the fact that we can also move the marked 
point. This makes the proof of transversality much easier, but it requires a higher bound on the 
ampleness of the line bundle. 
We believe that with some further effort we should be able to improve our bound (namely by moving 
both the curve and the point); 
we intend to investigate that question in the future. 
\end{rem}

\begin{rem}
In order to keep this paper of a reasonable length, 
%we 
%To keep the paper of a reasonable length, 
we 
have decided to split the proof of our main result into two papers. 
In this paper, we focus exclusively 
on the Euler class computation. We precisely state the relevant transversality, closure 
and multiplicity claims and then go on to compute the enumerative numbers.  
Our second paper \cite{BM_closure_of_seven_points} 
is devoted exclusively to proving the technical results on
closure and multiplicities that  we use in this paper. We also prove in \cite{BM_closure_of_seven_points} 
the necessary transversality 
results that we use in this paper. The transversality results actually follow without too much effort. 
In \cite{BM13_one_singular_point_published} and 
\cite{BM_Detail} we prove the relevant transversality results for $\mathbb{P}^2$ when there are 
one or two singular points. It is easy to see from the proof presented in \cite{BM_Detail}, how the 
proof of transversality generalizes to more than one singular point. In \cite{RM_Hypersurfaces_published} 
we show how to generalize the proof of transversality to the case of general linear system that 
is sufficiently ample. In particular, the proof of the relevant transversality results that we use in this paper 
is very similar to the proof presented in our earlier three papers (\cite{BM13_one_singular_point_published}, 
\cite{BM_Detail} and \cite{RM_Hypersurfaces_published}). 
%essentially
%follows from what we have proved in our earlier three papers.
\end{rem}

\section{A survey of related results} 
\label{survey}
\subsection{Counting singular curves in a linear system}
%\hf\hf ``Enumerative Geometry 
%of Singular Curves''  is a classical subject and has been studied extensively using tools of 
%algebraic geometry. As we explained in the introduction,
%section\ref{method_description}, 
%the crucial aspect of our method is the way we compute the degenerate contribution to the 
%Euler class. We perturb the section \textit{smoothly} (as opposed to 
%holomorphically) 
%and count (with a sign) how many zeros are there in a neighborhood 
%of the degenerate locus. Hence, the method is \textit{topological} as opposed to  
%\textit{algebro-geometric}.  \\ 
%namely. 
%\subsection{The results of Kazarian, Vainsencher and Kleiman and Piene}
\hf \hf We now give a brief survey of  
related results in this area of mathematics.
%, namely ``Enumerative Geometry of Singular Curves''. 
As early as $1873$, Zeuthen computed the number of nodal cubics and the number of tri nodal 
quartics in $\mathbb{P}^2$, passing through the appropriate number of generic points (cf \cite{Zeuthen}). 
A lot of progress has since been made in enumerating plane nodal curves; the question 
is now very well understood. 
%The question of enumerating plane nodal curves is now very well understood.
%Next, 
Ziv Ran and then 
Caporaso and Harris completely solved the problem of enumerating $\delta$-nodal curves in $\mathbb{P}^2$ 
(cf \cite{Ran1}, \cite{Ran2} and \cite{CH}). Caprasso and Harris gave a recursive formula for 
$\delta$-nodal curves, in terms of 
the number of curves with lesser number of nodes, but which have contacts of arbitrary 
order with a given generic line. Later on their result for $\mathbb{P}^2$ was extended by Ravi Vakil 
to certain del-Pezzo surfaces (cf. \cite{Vakil_CH}). An alternative proof of the Caporasso-Harris formula was also   
given by Ionel-Parker (cf. \cite{IP_SSF}) and Tehrani-Zinger (cf. \cite{MT_AZ}), 
using the Symplectic Sum formula.\footnote{It would be interesting to 
see if the result of Ravi Vakil can also be obtained using the Symplectic Sum Formula.}\\ 
\hf \hf Let us now look at the more general question of enumerating curves in an algebraic surface 
(such as $\mathbb{P}^1 \times \mathbb{P}^1$). In \cite{Van}, Vainsencher
obtained an explicit formula for curves in a sufficiently ample linear system that have upto six nodes.\\ 
%\hf \hf Next, let us look at the results 
%considers a general linear system $L \lra \X$ and enumerates curves that have upto six nodes. 
%We start by looking at the results 
%of M. Kazarian. 
\hf \hf Next, 
Kazarian computed \textit{all possible} codimension seven 
numbers in \cite{Kaz}.  
More precisely, he obtained a formula for the number of curves 
%(passing through $\frac{d(d+3)}{2}-\kappa$ generic points)  
in a sufficiently ample linear system (passing through the right 
number of generic points) and 
%and 
having a 
singularity of type $\chi_{k_1}, \chi_{k_2}, \ldots, \chi_{k_n}$, provided 
$\kappa:= k_1 + k_2+ \ldots+ k_n \leq 7$. 
Here $\chi_{k_i}$ is a singularity of codimension $k_i$.\\
\hf \hf We now mention the only two results that we are aware of in codimension eight.  
In \cite{WP1}, \cite{WP2} and \cite{WMP3}, 
A.Weber, M.Mikosz and P.Pragacz extended the method of Kazarian to enumerate curves 
with one singular point. In his slides \cite[page 52]{A_Weber}, Weber has explicitly written down the  formula for 
%For the problem of curves with one singular point, a formula has been obtained 
%\cite[page 263, (1.2)]{RT}
$N(A_8)$ and  $N(E_8)$; %(page $52$); 
our formulas agree with theirs.\\ 
%and 
%(and we believe it is very likely that they must have also obtained a formula for 
%$N(D_8)$ and $N(X_9)$).  .\\
%$by A.Weber (cf *). \\ 
%obtains a formula for the number of curves in a linear system with on singularity of type $A_8$.\\ 
\hf \hf Finally, we look at the results of S.Kleiman and R.Piene. 
%Let us now describe the results  S. Kleiman and R. Piene. 
In \cite[Theorem $1.2$]{KP}, Kleiman and Piene obtained 
an explicit formula for \textit{eight} nodal curves in a suitably ample linear 
system! 
%passing through the appropriate number of points that have 
%\textit{eight} nodes! 
They also obtained 
a formula for $\Num(\A_1^{\delta} D_4)$, $\Num(\A_1^{\delta} D_6)$ 
and $\Num(\A_1^{\delta} E_7)$, with $\delta$ at most $3$, $1$ 
and $0$ respectively (basically till codimension seven).\footnote{The paper by Kleiman and Piene is earlier than 
the paper by Kazarian; so their formulas for the codimension seven numbers were all new.} 
%many of the 
%formulas we have obtained in Theorems \ref{main_result}; 
%namely $\Num(A_1), \Num(D_4), \Num(D_6)$ and $\Num(E_7)$.
%, \N(A_1A_1)$, 
%$\N(A_1\D_4)$ and $\N(A_1 \D_6)$ 
%(it is part of Theorem $1.2$ in \cite{KP}). 
%They indeed consider a general linear system 
%$L \lra \X$, but 
%They study a \textit{completely different class of enumerative problems}
%(with the exception of $\N(A_1)$, $\N(A_2)$ and $\N(A_1A_1)$; our results 
%are consistent with theirs in this special case.) 
%They in fact obtain a formula for  
%Similarly, S.Klienman and R.Piene in their paper \cite{KP} consider a 
%general linear system $L \lra \X$ and 
%enumerating curves that have upto eight simple nodes, or one triple point 
%and up to three simple nodes, or one singularity of type $D_6$ and upto one simple node. 
%Their results are applicable for a linear system that is suitably ample. 
They also give a very precise condition for the ampleness requirement of their line bundle.  
They consider a line bundle 
$L \lra X$, 
%over a smooth projective surface $X$ 
that is 
of the form $L := M^{\otimes m} \otimes N$, 
where $m$ is at least three times the expected codimension of the 
equisingular strata and where $N$ is spanned by global holomorphic sections. 
%They then make the curves pass through the right number of generic points 
%so that the expected number is finite. 
Using this ampleness condition, they prove that the numbers they 
compute are indeed enumerative, i.e. each curve appears with a multiplicity 
of one in this linear system. \\
%Neither of them consider curves with a singularity other than a simple node 
%a cusp or a triple point. 
%(or in one case a cusp or a triple point).
%Our results are different from theirs; in Theorem \ref{main_result},  
%\ref{one_pt_chern_class_codim_7}, \ref{two_pt_chern_class_codim_7} 
%we enumerate curves 
%with upto one arbitrarily degenerate singularity (till a total codimension of seven). \\
%\hf\hf To summarize, the results of Kazarian, Vainsencher and Kleiman and Piene are applicable 
%for enumerating singular curves in an arbitrary smooth projective surface 
%(and not just for curves in $\P^2$). \\
%\hf \hf Next, we note that in 
%\cite{Ran1}, \cite{Ran2} and \cite{CH}, Z. Ran, L. Caporaso and J. Harris have 
%obtained a formula for the number of curves of degree $d$ in $\P^2$ (through the right 
%number of generic points) having $r$ simple nodes, for any $r$ (provided 
%$d$ is sufficiently larger than $r$). \\
%However, 
%their 
%results are only for $\P^2$. Moreover, 
%the allowed singularities 
%in their cases are simple nodes only and nothing more degenerate. \\
%\hf \hf The numbers $ \Num(\X)$ have also been computed by Dmitry Kerner for 
%many singularities in his paper \cite{Ker1} (for curves in $\P^2$). \\ 
%In some of the 
%formulas obtained by Kerner there is a numerical error made by the author. 
%In a personal communication \cite{Kerner_communication} the author informed us that he agrees that 
%he has made some error while computing few of these numbers. \\
\hf \hf Finally we note that a great deal of progress has been made in proving that 
universal formulas exist in terms of Chern classes 
for the number of curves in a sufficiently ample 
linear system passing through the right number of generic points 
and having singularities of type $\chi_{k_1}, \chi_{k_2}, \ldots, \chi_{k_n}$. 
The fact that a universal formula exists was conjectured by G{\"o}ttsche 
when the singularities $\chi_{k_1}, \ldots, \chi_{k_{n}}$ 
are all simple nodes. 
%There are 
%now quite a few proofs of this conjecture;  
Two independent proofs of this conjecture have now been given using methods of   
BPS calculus and using degeneration methods by Kool, Shende and Thomas (\cite{KST}) 
and Tzeng (\cite{Tz}) respectively. There is also an earlier approach by Liu in 
\cite{L2} and \cite{L1}. 
%In \cite{KST}, Kool, Shende and Thomas prove the conjecture by using BPS calculus. A different 
%proof was given by Tzeng in \cite{Tz} using degeneration methods (see also the approach of Liu, \cite{L1} 
%and \cite{L2}).  
%\cite{L1}, \cite{L2}, \cite{KST} and \cite{Tz}. 
Recently, Li and Tzeng gave a proof for the existence of universal formulas  
for any collection of singularities $\chi_{k_1}, \chi_{k_2}, \ldots, \chi_{k_n}$ in \cite{Tzeng_Li}, 
thereby generalizing the G{\"o}ttsche conjecture.  
A further generalization of the G{\"o}ttsche conjecture was proved by Rennemo,  
where he shows that there is a universal polynomial to count hypersurfaces in a 
sufficiently ample linear system, with any collection of isolated singularities (\cite[Proposition 7.8]{Rennemo3}).\\ 
\hf \hf It should be noted that even for $\delta$-nodal curves in $\P^2$, it is nontrivial question as to 
when the 
polynomials obtained 
by the G{\"o}ttsche conjecture are actually enumerative. They are enumerative 
if $d$ is large. However, G{\"o}ttsche conjectured that these polynomials are actually enumerative for 
all $d$ roughly greater than $\frac{\delta}{2}$. This has now been proven by Kleiman and Shende in \cite{KlSh}.\\  
\hf \hf For a more detailed overview of this subject, 
%(enumerative geometry of singular curves), 
we direct the reader to the comprehensive survey article \cite{Kl_survey2} by Kleiman.  

\subsection{Counting curves of a fixed genus: Gromov-Witten theory} 
\hf \hf Let us now consider another classical enumerative geometry problem; namely 
counting curves of a fixed genus. 
%which has also been investigated for over a 
%hundred years. This problem comprises of fixing the genus of a curve and 
%imposing enough constraints on the space of curves to get a number.
%\hf \hf Let us now discuss a different type of problem; namely counting curves of 
%a fixed genus. 
To begin with, let us look at genus zero curves (i.e. rational curves). 
Let $X$ be a compact complex surface and $\beta \in H_2(X, \mathbb{Z})$ a given homology class. 
A classical question is to ask what is $n_{\beta}$, the number of genus zero degree 
$\beta$ curves in $X$ passing through %$\delta_{\beta} :=% 
$\lan c_1(TX), ~\beta\ran -1$ 
generic points?
%, where 
%\[ \delta_{\beta} := \lan c_1(TX), ~\beta\ran -1. \]
%This question has also been studied since the late nineteenth century. 
For the case 
of $\mathbb{P}^2$, 
%plane rational degree $d$ curves 
%(i.e. when $X:= \mathbb{P}^2$), 
the numbers 
$n_1$, $n_2$, $n_3$ and $n_4$ were computed by Zeuthen in the late nineteenth century (\cite{Zeuthen}). 
There was essentially no progress in this question, until in the early $1990^{'\textnormal{s}}$ 
a formula for $n_{d}$ was announced by Kontsevich-Manin and Ruan-Tian. 
More precisely, in \cite{KoMa} and \cite{RT}, the authors obtain a formula for 
$n_{\beta}$ for any complex del-Pezzo surface.\footnote{Strictly speaking, Kontsevich-Manin compute the 
genus zero Gromov-Witten invariants of del-Pezzo surfaces. 
The fact that the Gromov-Witten invariants  are actually enumerative 
for all del-Pezzo surfaces was later proved by Pandharipande and  G\"ottsche in \cite{Rahul_Gottsche}.} 
%(which is $\mathbb{P}^1 \times \mathbb{P}^1$ or 
%$\mathbb{P}^2$ blown up at upto eight points in general position). 
Both the 
approaches involve Symplectic Geometry; it involves using a gluing theorem which can be algebraically 
reformulated as the associativity of Quantum Cohomology. A lot of progress has since been 
made in this question over the last thirty years using Symplectic and Algebraic geometry 
and more recently using Tropical geometry (cf \cite{GregMikh}). \\ 
\hf \hf The question of enumerating genus zero curves with higher 
singularities is a much harder question. 
%The question of enumerating 
%Rational cuspidal curves has been studied by 
Pandahripande and later Kock computed the number of rational cuspidal 
curves in $\mathbb{P}^2$ and $\mathbb{P}^1 \times \mathbb{P}^1$   
using algebro geometric methods (cf \cite{Rahul1} and \cite{JKock}). 
Later on, Zinger 
used methods from Symplectic Geometry to approach the question of enumerating rational curves 
with cusps and higher singularities (cf \cite{g2p2and3}, \cite{g3}, \cite{g1} and \cite{g0pr}). 
It should also be noted that although Zinger obtains his 
formulas for curves in $\mathbb{P}^2$, his methods are applicable in a much more general setup. 
In fact, his methods apply to surfaces where the Gromov-Witten Theory and its relationship with 
Enumerative geometry 
%of rational curves (without any singularities) 
is well understood. In particular, his methods 
apply very nicely to del-Pezzo surfaces. A few of Zinger's results have been 
extended by the second author to del-Pezzo surfaces (cf. \cite{BVSRM}, \cite{BMT1} and  \cite{BMT2}).  \\
\hf \hf Let us now explain how this question is relevant to the problem we are studying in this paper. 
The problem  of counting curves of a fixed genus enables us to subject our formulas to several 
low degree checks. For example, let us look at $n_d$, the number of degree $d$ rational curves in 
$\mathbb{P}^2$ passing through $3d-1$ generic points. A little bit of thought shows that  
$n_3$ is same as the number of nodal cubics through $8$ points; hence, this number helps us verify 
the formula for $N(A_1)$. Next, we consider $n_4$. A little bit of thought shows that this is 
same as the number of irreducible quartics with three (unordered) nodes through $11$ points. 
Hence, the number $n_4$ helps us verify the the formula for $N(A_1^3)$. And so on. \\ 
\hf \hf Let us now consider higher singularities. For example,in \cite{g2p2and3}, Zinger has computed the number 
of rational degree $d$ curves with one $E_6$-singularity. For $d=3$, this ought to be zero 
(which it is) and for $d=4$, his answer answer ought to be the same as the number of quartics 
through $8$ points having an $E_6$-singularity (again, that is indeed the case). 
Hence, these results for higher singularities also 
help us verify many of our formulas. \\
\hf \hf In section \ref{low_degree_checks}, we subject our formulas to several low degree checks. 
In section \ref{ldc_van}, we correct some of the low degree checks made by Vainsencher in his paper 
\cite{Van} and show how the numbers are consistent with the values predicted by the 
formula of Kontsevich and Manin in \cite{KoMa}. 
%We welcome additional ideas of low degree checks to verify our formulas. 

%\section{Description of the method we use} 

\section{A brief overview of the method we use} 
\hf \hf We now give a brief description of the method we use to compute the enumerative numbers. 
Our starting point is the following fact from differential topology 
(cf. \cite[Proposition 12.8]{BoTu}):
%, Proposition 12.8)
\begin{thm} 
\label{Main_Theorem_topology} 
Let $V\lra X$ be an oriented vector bundle over a compact oriented manifold $X$ and suppose
$s:X \lra V$ is a smooth section that is transverse to the zero set. Then the 
Poincar\'{e} dual of $[s^{-1}(0)]$ is the Euler class of $V$. In particular, if the 
rank of $V$ is same as the dimension of $X$ then 
%Let $V\lra X$ be a smooth ($C^\infty$) vector bundle over a smooth manifold $X$. Then the following are true: \\
%\hspace*{0.5cm}(1) A generic smooth section $s: X\lra V$ is transverse to the zero section. \\
%\hspace*{0.5cm}(2) Furthermore, if $V$ and $X$ are oriented with $X$ compact then the zero set of such 
%a section defines an integer homology class in $X$, whose Poincar\'{e} dual is the Euler class of $V$. 
%In particular, if the rank of $V$ is same as the dimension of $X$, 
the signed cardinality of $s^{-1}(0)$ is the Euler class of $V$, evaluated on the fundamental class of 
$X$, i.e., 
\bgd
|\pm s^{-1}(0)| = \lan e(V), [X] \ran. 
\edd
\end{thm}
Our goal is to enumerate curves with certain singularities. The fact that a curve has a specific  
singularity means that certain derivatives vanish (an example of this fact is the Implicit Function Theorem 
and the Morse Lemma). We interpret these derivatives as the section of some bundle.   
Hence, our enumerative numbers are the zeros  
of a section of some bundle restricted to the \textit{open} part of a manifold (or variety).
One of the reasons we typically have to restrict ourselves to the open part of a manifold is because 
we are usually enumerating curves with more than one singular point; hence we have to consider the 
space of curves with a collection of marked points, where the marked points are all distinct.  \\
\hf \hf Next, we observe that if the line bundle is sufficiently ample, then this section (induced by taking 
derivatives) will be transverse 
to the zero set restricted to the open part.\footnote{The open part of 
the variety we consider is always going to be smooth in our case.}  
We then evaluate the Euler class of this bundle on the fundamental class   
of the manifold/variety\footnote{Any algebraic variety defines a homology class since the singularities have at least 
real codimension two. This follows from the main result of \cite{pseudo_cycle}, namely 
that any singular space whose singularities are of real codimesnion two  or more 
(i.e. a pseudocycle)
defines 
a homology class.} 
and hope that this number is what we want. 
As one might expect, that hope 
is almost always too naive. 
This is because, the section will typically vanish on the boundary 
and hence give an excess (degenerate) contribution to the Euler class. 
The real challenge is therefore 
to compute this degenerate contribution to the Euler class. \\ 
\hf \hf The most famous and well known method of computing degenerate contributions to the 
Euler class is global excess intersection theory, which is developed in Fulton's book \cite{F}. 
This approach is global and rigid in nature and 
involves blowing up the degenerate loci to smooth them out. 
This method has been successfully applied to solve a large class of enumerative problems; 
in particular it has been applied by  
%formulas for $\delta$-nodal curves obtained by 
Vainsencher and Kleiman and Piene to obtain their formula for $\delta$-nodal curves. \\ 
%involve applying Fulton's intersection theory. \\ 
\hf \hf In his PhD thesis \cite{Zinger_Thesis}, Aleksey Zinger developed a completely different approach to 
address the problem of computing degenerate contributions to the Euler class; he  
developed the method of \textit{local intersection theory}. 
This approach is soft and local in nature and does not involve 
taking any blowups. 
%involves blowing up the degenerate loci to smooth them out. 
The method involves perturbing the relevant section 
\textit{smoothly} and counting (with a sign) the number of zeros of the section near the degenerate 
locus. He applied this method successfully to solve a very large class of enumerative problems; 
these include counting rational curves in $\mathbb{P}^2$ (and even $\mathbb{P}^n$) with various types 
of singularities (which include cusp, tacnode, triple point and $E_6$-singularities). He also 
used local intersection theory to compute the degenerate contribution to the Symplectic Invariant 
of $\mathbb{P}^n$ 
%(and even $\mathbb{P}^n$) 
and computed the 
genus $g$ \textit{Enumerative} Invariant of 
$\mathbb{P}^n$
%number of genus $g$ curves in 
%$\mathbb{P}^n$ with a fixed complex structure passing through the appropriate number of points, 
for $g=2$ and $3$, 
thereby extending the result of Pandharipande and Ionel for genus one (cf. \cite{Rahul_genus_one} 
and \cite{Ionel_genus_one}).\\ 
\hf \hf Zinger further applied the machinery of local intersection theory in his paper \cite{Z1}
to enumerate plane degree $d$ curves with singularities till codimension three. 
In particular, he computed 
all the possible codimension three numbers $N(A_1^{\delta} \mathfrak{X})$ in his paper by using 
local intersection theory.\\ %(his paper is earlier than the paper by Kazarian).  \\ 
\hf \hf The second author's PhD thesis (cf \cite{RM_thesis}) comprises of extending the approach
developed in \cite{Z1} and 
applying \textit{local} intersection 
theory to compute $N(A_1^{\delta} \mathfrak{X})$ for higher codimension. In his thesis, 
the second 
author obtained a formula for $N(A_1^{\delta} \mathfrak{X})$ for $\mathbb{P}^2$ when $\delta+k \leq 7$; 
part of that thesis has been published in  \cite{BM13_one_singular_point_published} and 
\cite{BM13_2pt_published} (namely the cases $\delta=0$ and $1$ have been published). 
In this paper we finally obtain 
a formula for all the numbers $N(A_1^{\delta} \mathfrak{X})$ for any sufficiently ample 
linear system 
till codimension $8$, thereby going beyond the results of the second author's thesis.  
Furthermore, we give \textit{alternative} proofs for the results of Vainsencher, Kleiman-Piene 
and Kazarian\footnote{However, we should mention that we only compute numbers of the form 
$N(A_1^{\delta} \mathfrak{X})$; Kazarian computes all numbers till codimension seven; such as $N(A_3 A_4)$ for 
example.} 
and we also 
%recovering many of the earlier results by others and also 
obtain new results in the end (basically all the codimension eight numbers are new 
except for the formula for eight nodal curves and the formulas for one singular point). \\
%the ones obtained by Kleiman and Piene and ). 
%It should be mentioned here that although the numbers we have obtained in 
%codimension eight are all new (except for the ones obtained by Kleiman and Piene and *), we give a 
%completely different proof for all the other numbers we have obtained. In particular our method is 
%completely different from 
%from Kazarian's method. Furthermore, our proof for $r$-nodal curves is completely different 
%from those given by Vainsincher, and Klieman and Piene.\\ 
%Hence, we also obtain alternative proofs of 
%their results.  
\hf \hf Finally, let us give a brief description of the method used by Kazarian. 
His method works on the principle that there exists a universal formula (in terms of Chern classes) 
for 
the Thom polynomial associated to a given singularity.
%these 
%enumerative numbers in terms of the Chern classes. 
He then goes on to consider 
enough special cases to find out what that exact combination is. As an example, 
suppose there is a polynomial of degree $m$. To find out what the polynomial is, 
we simply have to find the value of the polynomial at enough points (\cite[page 667]{Kaz}).  
One of the challenges of this method is to prove the existence of such a universal 
formula. This fact is shown by R.Rimanyi and A. Szucs in \cite{Rimanyi1}, \cite{Rimanyi2} and \cite{Rimanyi_Szucs3}. 
This method has been successfully applied to compute all the codimension seven numbers. 
It has also been extended by A.Weber, M.Mikosz and P.Pragacz 
to enumerate curves with one singular point in codimension eight 
(cf. \cite{WP1}, \cite{WP2} and \cite{WMP3}). 
%We also believe it is usually difficult  
%to think of enough special cases in a 
%given situation. 

%\section{Organization of our paper(s)} 
%Let us now explain how our paper is organized. To keep the paper of a reasonable length, we 
%have decided to split the proof of our main result into two papers. In this paper, we focus exclusively 
%on the Euler class computation. We state the relevant transversality, closure 
%and multiplicity claims precisely and then go on to compute the enumerative numbers. \\ 
%\hf \hf Our second paper \cite{BM_closure_of_seven_points} 
%is devoted exclusively to proving the technical results on, 
%closure and multiplicities that  we use in this paper. We also prove the necessary transversality 
%results that we use in this paper. The transversality results actually follow without too much effort. 
%In \cite{BM13_one_singular_point_published} and 
%\cite{BM_Detail} we prove the relevant transversality results for $\mathbb{P}^2$ when there are 
%one or two singular points. It is easy to see from the proof presented in \cite{BM_Detail}, how the 
%proof of transversality generalizes to more than one singular point. In \cite{RM_Hypersurfaces_published} 
%we show how to generalize the proof of transversality to the case of general linear system that 
%is sufficiently ample. Hence, the relevant transversality results that we use in this paper essentially
%follows from what we have proved in our earlier three papers.
%However, for the convenience of the reader, 
%we intend to rewrite 
%\cite{BM13_one_singular_point_published}, \cite{BM_Detail} 
%and \cite{RM}

\section{Necessary and sufficient criteria for a singularity}
\label{condition_for_sing}
\hf\hf For the convenience of the reader, let us recapitulate from the arXiv version of our 
paper \cite[section 5]{BM13}, 
a necessary and sufficient criterion for a curve $f^{-1}(0)$ to have a singularity 
of type $\mathfrak{X}$ at the point $q$.\footnote{The arXiv version contains a few details 
that were omitted from our published paper \cite{BM13_one_singular_point_published}; 
hence we will often refer to the 
arXiv version of our paper.} 
Let $\q=\q(x,y)$ be a holomorphic function defined on a neighborhood
of the origin in $\C^2$ and $i,j$ be non-negative integers. We define
\bgd
\q_{ij}:=\frac{\partial^{i+j} \q}{\partial^i x\partial^j y}\bigg|_{(x,y)=(0,0)}\,.
\edd
Let us now define the following directional derivatives, which are functions of $\q_{ij}$:  
\begin{align}
\A^{\q}_3&:= \q_{30}\,,\qquad
\A^{\q}_4 := \q_{40}-\frac{3 \q_{21}^2}{\q_{02}}\,, \qquad
\A^{\q}_5:= \q_{50} -\frac{10 \q_{21} \q_{31}}{\q_{02}} + 
\frac{15 \q_{12} \q_{21}^2}{\q_{02}^2}, \nonumber \\
\A^{\q}_6 &:= \q_{60}- \f{ 15 \q_{21} \q_{41}}{\q_{02}}-\f{10 \q_{31}^2}{\q_{02}} + \f{60 \q_{12} \q_{21} \q_{31}}{\q_{02}^2}
   +
   \f{45 \q_{21}^2 \q_{22}}{\q_{02}^2} - \f{15 \q_{03} \q_{21}^3}{\q_{02}^3}
   -\f{90 \q_{12}^2 \q_{21}^2}{\q_{02}^3}, \nonumber \\  
\A^{\q}_7 &:= \q_{70} - \frac{21 \q_{21} \q_{51}}{\q_{02}} 
- \frac{35 \q_{31} \q_{41}}{\q_{02}} + \frac{105 \q_{12} \q_{21} \q_{41}}{\q_{02}^2} + \f{105 \q_{21}^2 \q_{32}}{\q_{02}^2} + 
\f{70 \q_{12} \q_{31}^2}{\q_{02}^2}+ \f{210 \q_{21}\q_{22}\q_{31}}{\q_{02}^2} \nonumber \\
&
-\f{105 \q_{03} \q_{21}^2 \q_{31}}{ \q_{02}^3}
-\f{420 \q_{12}^2 \q_{21} \q_{31}}{\q_{02}^3}
-\f{630 \q_{12} \q_{21}^2 \q_{22}}{\q_{02}^3}
-\f{105 \q_{13} \q_{21}^3}{\q_{02}^3}
+ \f{315 \q_{03} \q_{12} \q_{21}^3}{\q_{02}^4}
+ \f{630 \q_{12}^3 \q_{21}^2}{\q_{02}^4}, \nonumber
\end{align}
\begin{align}
\label{Ak_sections}
%\A^{\q}_8 &= -\frac{f(2,1)^2 f(1,2)^4}{8 f(0,2)^5}+\frac{f(2,1) f(3,1) f(1,2)^3}{12 f(0,2)^4}- *\\ 
%         & \frac{f(0,3) f(2,1)^3 f(1,2)^2}{8
%   f(0,2)^5}-\frac{f(3,1)^2 f(1,2)^2}{72 f(0,2)^3}+\frac{3 f(2,1)^2 f(2,2) f(1,2)^2}{16 f(0,2)^4} *\\ 
%   & -\frac{f(2,1) f(4,1)
%   f(1,2)^2}{48 f(0,2)^3}+\frac{f(1,3) f(2,1)^3 f(1,2)}{16 f(0,2)^4}+\frac{f(0,3) f(2,1)^2 f(3,1) f(1,2)}{16
%   f(0,2)^4}-\frac{f(2,1) f(2,2) f(3,1) f(1,2)}{12 f(0,2)^3} * \\ 
%   & -\frac{f(2,1)^2 f(3,2) f(1,2)}{24 f(0,2)^3}+\frac{f(3,1) f(4,1)
%   f(1,2)}{144 f(0,2)^2}+\frac{f(2,1) f(5,1) f(1,2)}{240 f(0,2)^2}-\frac{f(0,3)^2 f(2,1)^4}{128 f(0,2)^5}+ *\\ 
%   & \frac{f(0,4)
%   f(2,1)^4}{384 f(0,2)^4}-\frac{f(2,1)^2 f(2,2)^2}{32 f(0,2)^3}-\frac{f(0,3) f(2,1) f(3,1)^2}{144 f(0,2)^3}+ *\\ 
%    & \frac{f(2,2)
%   f(3,1)^2}{144 f(0,2)^2}-\frac{f(4,1)^2}{1152 f(0,2)}+\frac{f(0,3) f(2,1)^3 f(2,2)}{32 f(0,2)^4}- *\\ 
%   & \frac{f(2,1)^3 f(2,3)}{96
%   f(0,2)^3}-\frac{f(1,3) f(2,1)^2 f(3,1)}{48 f(0,2)^3}+\frac{f(2,1) f(3,1) f(3,2)}{72 f(0,2)^2}- *\\ 
%   & \frac{f(0,3) f(2,1)^2
%   f(4,1)}{192 f(0,2)^3}+\frac{f(2,1) f(2,2) f(4,1)}{96 f(0,2)^2}+\frac{f(2,1)^2 f(4,2)}{192 f(0,2)^2}- *\\ 
%   & \frac{f(3,1) f(5,1)}{720
%   f(0,2)}-\frac{f(2,1) f(6,1)}{1440 f(0,2)}+\frac{f(8,0)}{40320}*
\A^{\q}_8 &:= \q_{80} -\frac{28 \q_{21} \q_{61}}{\q_{02}}
-\frac{56 \q_{31} \q_{51}}{\q_{02}} + \frac{210 \q_{21}^2 \q_{42}}{\q_{02}^2} 
+\frac{420 \q_{21} \q_{22} \q_{41}}{\q_{02}^2}-\frac{210 \q_{03} \q_{21}^2
   \q_{41}}{\q_{02}^3} +\frac{560 \q_{21} \q_{31} \q_{32}}{\q_{02}^2}- \nonumber \\ 
   & \frac{840 \q_{13} \q_{21}^2 \q_{31}}{\q_{02}^3}- \frac{ 420 \q_{21}^3 \q_{23}}{\q_{02}^3} 
   +\frac{1260 \q_{03} \q_{21}^3 \q_{22}}{\q_{04}^4}-\frac{35 \q_{41}^2}{\q_{02}}+ 
   \frac{280 \q_{22}
   \q_{31}^2}{\q_{02}^2}-\frac{280 \q_{03} \q_{21} \q_{31}^2}{\q_{02}^3}
   -\frac{1260 \q_{21}^2 \q_{22}^2}{\q_{02}^3}+  \nonumber \\ 
   & \frac{105 \q_{04}
   \q_{21}^4}{\q_{02}^4}-\frac{315 \q_{03}^2 \q_{21}^4}{\q_{02}^5}
   +\frac{168 \q_{21} \q_{51} \q_{12}}{\q_{02}^2}+\frac{280 \q_{31} \q_{41}
   \q_{12}}{\q_{02}^2}-\frac{1680 \q_{21}^2 \q_{32} \q_{12}}{\q_{02}^3}
   -\frac{3360 \q_{21} \q_{22} \q_{31} \q_{12}}{\q_{02}^3} + \nonumber \\
   & \frac{2520 \q_{03} \q_{21}^2 \q_{31} \q_{12}}{\q_{02}^4}+  
   \frac{2520 \q_{13} \q_{21}^3 \q_{12}}{\q_{02}^4} 
    -\frac{840 \q_{21} \q_{41}
   \q_{12}^2}{\q_{02}^3} +\frac{7560 \q_{21}^2 \q_{22} \q_{12}^2}{\q_{02}^4}
   -\frac{560 \q_{31}^2 \q_{12}^2}{\q_{02}^3}-\frac{5040 \q_{03} \q_{21}^3 \q_{12}^2}{
   \q_{02}^5} \nonumber \\ 
   & +\frac{3360 \q_{21} \q_{31} \q_{12}^3}{\q_{02}^4} 
   -\frac{5040 \q_{21}^2 \q_{12}^4}{\q_{02}^5}
\end{align} 
and 
\begin{align}
\label{Formula_Dk_again}
\D^{\q}_6 &:=  \q_{40},\qquad 
\D^{\q}_7 :=   \q_{50} -\f{5 \q_{31}^2}{3 \q_{12}},\qquad
\D^{\q}_8 := \q_{60} + 
\frac{5 \q_{03} \q_{31} \q_{50}}{3 \q_{12}^2} 
-\frac{5 \q_{31} \q_{41}}{\q_{12}} - \frac{10 \q_{03} \q_{31}^3}{3 \q_{12}^3} 
+ \frac{5 \q_{22} \q_{31}^2}{\q_{12}^2}.
\end{align}

%\begin{lmm}
%\label{ift}
%Let $\q =\q(x,y)$ be a holomorphic function defined on a neighborhood
%of the origin in $\C^2$ such that $\q_{00} =0$ and $\nabla \q |_{(0,0)} \neq 0$. 
%Then the curve $\q^{-1}(0)$ has an $\A_0$-node at the origin. 
%Then there exists a coordinate 
%chart $(u,v)$ centered at the origin so that $\q(u,v)= v^2 + u.$ 
%\end{lmm}
%\pf Follows immediately by considering the Taylor expansion of $\q$. \qed   
%\begin{cor}
%\label{A0_node_condition_cor}
%A curve $\q^{-1}(0)$ has an $\A_0$-node at the origin if and only if it satisfies the hypothesis of Lemma \ref{ift}. 
%\end{cor}
%\begin{lmm}
%\label{ml}
%Let $\q = \q(x,y)$ be a holomorphic function defined on a neighbourhood
%of the origin in $\C$ such that $\q_{00}, \nabla \q |_{(0,0)} =0 $ and
%$\nabla^2 \q|_{(0,0)}$ is non-degenerate.  Then there exists a coordinate 
%chart $(u,v)$ centered at the origin so that $\q(u,v) = v^2 + u^2.$
%\end{lmm}
%\pf This is the Morse Lemma, which again follows by considering the Taylor expansion of $\q$.   \qed 
%\begin{cor}
%\label{A1_node_condition_cor}
%A curve $\q^{-1}(0)$ has an $\A_1$-node if and only if it satisfies the hypothesis of Lemma \ref{ml}. 
%\end{cor}
We will now state a necessary and sufficient criteria for a curve to have a specific singularity. 

\begin{lmm}
\label{ift}
Let $\q =\q(x,y)$ be a holomorphic function defined on a neighborhood
of the origin in $\C^2$ such that $\q_{00} =0$ and $\nabla \q |_{(0,0)} \neq 0$. 
Then the curve has an $\A_0$ singularity at the origin (i.e. a smooth point). 
%Then 
%there exists a coordinate 
%chart $(u,v)$ centered at the origin so that $\q(u,v)= v^2 + u.$ 
\end{lmm}
%\pf Follows immediately by considering the Taylor expansion of $\q$. \qed   
%\begin{cor}
%\label{A0_node_condition_cor}
%A curve $\q^{-1}(0)$ has an $\A_0$-node at the origin if and only if it satisfies the hypothesis of Lemma \ref{ift}. 
%\end{cor}
\begin{lmm}
\label{ml}
Let $\q = \q(x,y)$ be a holomorphic function defined on a neighbourhood
of the origin in $\C$ such that $\q_{00}, \nabla \q |_{(0,0)} =0 $ and
$\nabla^2 \q|_{(0,0)}$ is non-degenerate.  
Then the curve has an $\A_1$ singularity at the origin. 
%Then there exists a coordinate 
%chart $(u,v)$ centered at the origin so that $\q(u,v) = v^2 + u^2.$
\end{lmm}

\begin{rem}
Lemma \ref{ift} is also known as the \textbf{Implicit Function Theorem} and Lemma \ref{ml} 
is also known as the \textbf{Morse Lemma}. 
%This is the Morse Lemma, which again follows by considering the Taylor expansion of $\q$.   \qed 
\end{rem}

\ni We now state the remaining Lemmas, which can be thought of as a continuation of Lemma \ref{ml}. 

\begin{lmm}\label{fstr_prp}
Let $\q =\q(\rq, \sq)$ be a holomorphic function defined on a neighbourhood of the origin in $\C$ such that 
$\q(0,0), ~\nabla \q|_{(0,0)}=0$ and there exists a non-zero vector $w=(w_1,w_2)$ such that at 
the origin $ \nabla^2 f (w, \cdot)=0$, i.e., the Hessian is degenerate.
Let $x = w_1 \rq + w_2 \sq, y = -\ov{w}_2 \rq + \ov{w}_1 \sq $ and $\q_{ij}$ be the partial derivatives with 
respect to the new variables $x$ and $y$. 
Then, the curve $\q^{-1}(0)$ has an $\A_k$ singularity at the origin  
(for $2 \leq k \leq 7$) if $\q_{02} \neq 0$ and 
the directional derivatives 
$\A^{\q}_i$ defined in \eqref{Ak_sections} are zero for all $i \leq k$ and $\A^{\q}_{k+1} \neq 0$. 
\end{lmm}
%\begin{rem}
%In terms of the new coordinates we have $\q_{00}= \q_{10}= \q_{01}= \q_{20}= \q_{11} =0$ and $\q_{02} \neq 0.$ 
%Here $\partial_x + 0 \partial_y = (1,0)$ is the distinguished direction along which the Hessian is degenerate. 
%\end{rem}  

\begin{lmm}\label{fstr_prp_D4}
Let $\q = \q(x,y)$ be a holomorphic function defined on a neighbourhood
of the origin in $\C$ such that $\q_{00}, \nabla \q|_{(0,0)}, \nabla^2 \q|_{(0,0)}=0$ and 
there does not exist a non-zero vector $w=(w_1,w_2)$ such that 
at the origin $\nabla^3 \q (w,w,\cdot) =0$.
Then the curve $\q^{-1}(0)$ has a $\D_4$ singularity at the origin. 
\end{lmm}

\begin{lmm}\label{fstr_prp_Dk}
Let $\q=\q(\rq,\sq)$ be a holomorphic function defined on a neighbourhood 
of the origin in $\C$ such that $\q_{00}, \nabla \q|_{(0,0)}, \nabla^2 \q|_{(0,0)}=0$ and 
there exists a non-zero vector $w=(w_1,w_2)$ such that at the origin $\nabla^3 \q(w,w,\cdot) =0$.
Let $x = w_1 \rq + w_2 \sq, ~y = -\ov{w}_2 \rq + \ov{w}_1 \sq$ and $\q_{ij}$ be the partial derivatives with respect to 
the new variables $x$ and $y$. 
Then, the curve $\q^{-1}(0)$ has a $\D_k$-node at the origin  
(for $5 \leq k \leq 7$) if $\q_{12} \neq 0$ and 
the directional derivatives 
$\D^{\q}_i$ defined in \eqref{Formula_Dk_again} are zero for all $i \leq k$ and $\D^{\q}_{k+1} \neq 0$. 
\end{lmm}
%\begin{rem}
%Note that in terms of the new coordinates we have  
%\bgd
%\q_{00} = \q_{10}= \q_{01}= \q_{20}= \q_{11}= \q_{02}= \q_{30}= \q_{21} =0, \q_{12} \neq 0.
%\edd 
%\end{rem} 

\begin{lmm}\label{fstr_prp_E6}
Let $\q = \q(\rq,\sq)$ be a holomorphic function defined on a neighbourhood
of the origin in $\C$ such that $\q_{00}= \nabla \q|_{(0,0)}= \nabla^2 \q|_{(0,0)}=0$ and 
there exists a non-zero vector $w=(w_1,w_2)$ such that at the origin $\nabla^3 \q(w,w,\cdot) =0$. 
Let $x = w_1 \rq + w_2 \sq,  ~y = -\ov{w}_2 \rq + \ov{w}_1 \sq$ and $\q_{ij}$ be partial derivatives with respect 
to the new coordinates, $x$ and $y$. 
Then, the curve $\q^{-1}(0)$ has an $\E_6$ singularity at the origin 
if $\q_{12} =0$ and  $\q_{03}\neq 0, \q_{40} \neq 0$. 
\end{lmm}
%\begin{rem}
%Note that in terms of the new coordinates $x$ and $y$, we get that 
%\bgd
%\q_{00} =\q_{10}= \q_{01}=\q_{20}= \q_{11}= \q_{02}= \q_{30}= \q_{21}= \q_{12} =0, \q_{03}\neq 0, \q_{40} \neq 0.
%\edd 
%\end{rem}

\begin{lmm}\label{fstr_prp_E7}
Let $\q =\q(\rq,\sq)$ be a holomorphic function defined on a neighbourhood of the origin in $\C$ such that 
$\q_{00}, \nabla \q|_{(0,0)}, \nabla^2 \q|_{(0,0)}=0$ and there exists a non-zero vector $w=(w_1,w_2)$ such that 
at the origin $\nabla^3 \q(w,w,\cdot) =0$. Let $x = w_1 \rq + w_2 \sq, ~y = -\ov{w}_2 \rq + \ov{w}_1 \sq.$ Let $\q_{ij}$ be the partial derivatives with respect 
to the new variables $x$ and $y$. 
Then, the curve $\q^{-1}(0)$ has an $\E_7$ singularity at the origin 
if $\q_{12} =0$ and  $\q_{03}\neq 0, \q_{31} \neq 0$. 
\end{lmm}
%\begin{rem}
%In terms of the new coordinates $x$ and $y$, we get that 
%\bgd
%\q_{00}=\q_{10}= \q_{01}=\q_{20}= \q_{11}= \q_{02} =\q_{30}= \q_{21}= \q_{12}= \q_{40} =0, \q_{03}\neq 0, \q_{31} \neq 0.
%\edd
%\end{rem} 

\begin{lmm}\label{fstr_prp_E8}
Let $\q =\q(\rq,\sq)$ be a holomorphic function defined on a neighbourhood of the origin in $\C$ such that 
$\q_{00}, \nabla \q|_{(0,0)}, \nabla^2 \q|_{(0,0)}=0$ and there exists a non-zero vector $w=(w_1,w_2)$ such that 
at the origin $\nabla^3 \q(w,w,\cdot) =0$. Let $x = w_1 \rq + w_2 \sq, ~y = -\ov{w}_2 \rq + \ov{w}_1 \sq.$
Let $\q_{ij}$ be the partial derivatives with respect 
to the new variables $x$ and $y$. 
Then, the curve $\q^{-1}(0)$ has an $\E_8$ singularity at the origin 
if $\q_{12} =0, \q_{31} =0$ and  $\q_{03}\neq 0, \q_{50} \neq 0$. 
\end{lmm}

\textbf{Proofs of Lemmas \ref{fstr_prp} to \ref{fstr_prp_E7}}: This is proved in section $3$ of our paper 
\cite{BM13_one_singular_point_published}. The proof of Lemma \ref{fstr_prp_E8} is very similar to 
the proof of the other Lemmas; it involves looking at the Taylor expansion and making an appropriate 
change of coordinates. \qed  

\section{Notation}\label{notation}
\hf\hf We will now 
introduce some basic notation that will be used for stating our recursive formulas.  
Our setup is as follows; we consider a sufficiently ample line bundle $L \lra X$ 
over a compact complex surface $X$. 
Let
\[ \DD \approx \P H^0(X, L) \approx \P^{\delta_{L}}, \]
denote the space of holomorphic sections of the line bundle $L\lra X$, upto scaling. 
Let $\mathfrak{X}$ be a singularity of a given type. 
We will also denote $\mathfrak{X}$ to be the space of curves and a marked point 
such that the curve has a singularity of type $\mathfrak{X}$ at the marked point. 
More precisely, 
\begin{align*}
\mathfrak{X} &:= \{ ([f], q) \in \DD \times \X: f ~~\textnormal{has a signularity of type $\mathfrak{X}$ at $q$} \}. 
\end{align*}
For example, 
\begin{align*}
A_2 &:= \{ ([f], q) \in \DD \times \X: f ~~\textnormal{has a signularity of type $A_2$ at $q$} \}. 
\end{align*}
Next, given $n$ subsets $S_1, S_2, \ldots, S_n$ of $\DD \times \mathbb{P}^2$, we define 
\begin{align*}
S_1 \circ S_2 \circ \ldots \circ S_n &:= \{ ([f], q_1, \ldots, q_n) \in \DD \times (\mathbb{P}^2)^n:  
([f], q_1) \in S_1, \ldots, ([f], q_n) \in S_n,  \\ %\qquad \forall i = 1 ~~\textnormal{to} ~~n, \\ 
& \qquad \qquad \qquad \qquad \qquad \qquad \textnormal{and} \qquad q_1, \ldots, q_n 
\qquad \textnormal{are all distinct}\}.
\end{align*}
%More generally, we define $A_1^{\delta} \circ \mathfrak{X}$ to be the subspace of curves with $\delta+1$ 
%ordered points, such that the curve has a node at the first 
%$\delta$ points and a singularity of type $\mathfrak{X}$ at the last point and all the $\delta+1$ 
%points are distinct. 
For example, $A_1^{2} \circ A_2$ is the space of curves with 
three ordered points, where the curve has a simple node at the first two points and a cusp at the 
last point and all the three points are distinct. 
Similarly, $A_1^{2} \circ \overline{A}_2$
is the space of curves with 
three distinct ordered points, where the curve has a simple node at the first two points and a 
singularity at least as degenerate as a cusp at the 
last point; the curve could have a tacnode at the last marked point. \\ 
%, i.e.,
%\begin{align*}
%A_1^{2} \circ A_3 &:= \{ ([f], q_1, q_2, q_3) \in \DD \times (\X)^3: f ~~\textnormal{has a signularity of type $A_1$ at 
%$q_1$ and $q_2$ and} \\ 
%& \qquad \qquad \qquad \qquad \qquad \qquad \qquad \textnormal{a singularity 
%of type $A_3$ at $q_3$ and $q_1, q_2, q_3$ are distinct}\}. 
%\end{align*}
\hf \hf It will also be useful to consider the space of curves with a singularity and a specific direction 
along which certain directional derivatives vanish. More precisely, let $\P T\X \lra X$ denote the projectivization 
of the tangent bundle $\X$ and let $\mathfrak{X}$ be a singularity of type $A_{k \geq 2}$ or $D_{k \geq 5}$ 
or $E_6, E_7$ or $E_8$. Then we define 
\begin{align*}
%\mathfrak{X} &:= \{ ([f], q) \in \DD \times \X: f ~~\textnormal{has a signularity of type $\mathfrak{X}$ at $q$} \}, 
%\\
%\hat{\mathfrak{X}} &:= \{ ([f], l_q) \in \DD \times \mathbb{P}T\X: ([f], q) \in \mathfrak{X} \}, 
%\\
\mathcal{P}\A_k &:= \{ ([f], l_q) \in \DD \times \mathbb{P} T\X: ([f], q) \in \A_k, ~~\nabla^2 f(v, \cdot) =0 
\qquad \forall v \in l_q\} \qquad \textnormal{if ~~$k \geq 2$}, \\
%\mathcal{P}\D_4 &:= \{ ([f], l_q) \in \DD \times \mathbb{P} T\X: ([f], q) \in \D_4, ~~\nabla^3 f(v,v, v) =0 
%\quad \forall v \in l_q\}, \\
\mathcal{P}\D_k &:= \{ ([f], l_q) \in \DD \times \mathbb{P} T\X: ([f], q) \in \D_k, ~~\nabla^3 f(v,v, \cdot) =0 
\quad \forall v \in l_q\} \qquad \textnormal{if ~~$k \geq 5$}, \\
\mathcal{P}\E_k &:= \{ ([f], l_q) \in \DD \times \mathbb{P} T\X: ([f], q) \in \E_k, ~~\nabla^3 f(v,v, \cdot) =0 
\quad \forall v \in l_q\} \qquad \textnormal{if ~~$k=6,7,8$}.
%\PP\D_k^{\vee} &:= \{ ([f], q) \in \D \times \P T\P^2: ([f], q) \in \D_k, 
%~~\nabla^3 f|_q(v,v,v) =0, ~~\nabla^3 f|_q(v,v,w) \neq 0,  \\
%& \qquad \qquad \forall ~~v  \in  l_{q} -0 ~~\textup{and} ~~w \in (T_{q}\X)/l_{q} -0 \}, \qquad \textnormal{if $~k>4$}. 
\end{align*}
For example, $\PP A_2$ is the space of curves with a marked point and a marked direction, such that 
the curve has a cusp at the marked point and the marked direction belongs to the kernel of the Hessian. 
Note that in all the above cases, 
%we have defined above, 
the projection map $\pi : \PP \mathfrak{X} \lra \mathfrak{X}$ 
is one to one. 
%2(it is however not one to one when we pass to the closure).
%$\footnote{The}$
%$\footnote{It is however, not one to one when we pass to the closure.}$ 
Next, let us define 
\begin{align*}
%\mathfrak{X} &:= \{ ([f], q) \in \DD \times \X: f ~~\textnormal{has a signularity of type $\mathfrak{X}$ at $q$} \}, 
%\\
%\hat{\mathfrak{X}} &:= \{ ([f], l_q) \in \DD \times \mathbb{P}T\X: ([f], q) \in \mathfrak{X} \}, 
%\\
\mathcal{P}\A_1 &:= \{ ([f], l_q) \in \DD \times \mathbb{P} T\X: ([f], q) \in \A_1, ~~\nabla^2 f(v, v) =0 
\qquad \forall v \in l_q\}, \\
\mathcal{P}\D_4 &:= \{ ([f], l_q) \in \DD \times \mathbb{P} T\X: ([f], q) \in \D_4, ~~\nabla^3 f(v,v, v) =0 
\quad \forall v \in l_q\}, \\
\mathcal{P}\X_9 &:= \{ ([f], l_q) \in \DD \times \mathbb{P} T\X: ([f], q) \in X_9, ~~\nabla^4 f(v,v, v, v) =0 
\quad \forall v \in l_q\}. 
%%, \\
%\mathcal{P}\D_k &:= \{ ([f], l_q) \in \DD \times \mathbb{P} T\X: ([f], q) \in \D_k, ~~\nabla^3 f(v,v, \cdot) =0 
%\quad \forall v \in l_q\} \qquad \textnormal{if ~~$k \geq 5$}, \\
%\mathcal{P}\E_k &:= \{ ([f], l_q) \in \DD \times \mathbb{P} T\X: ([f], q) \in \E_k, ~~\nabla^3 f(v,v, \cdot) =0 
%\quad \forall v \in l_q\} \qquad \textnormal{if ~~$k=6,7,8$}.
%\PP\D_k^{\vee} &:= \{ ([f], q) \in \D \times \P T\P^2: ([f], q) \in \D_k, 
%~~\nabla^3 f|_q(v,v,v) =0, ~~\nabla^3 f|_q(v,v,w) \neq 0,  \\
%& \qquad \qquad \forall ~~v  \in  l_{q} -0 ~~\textup{and} ~~w \in (T_{q}\X)/l_{q} -0 \}, \qquad \textnormal{if $~k>4$}. 
\end{align*}
For example, $\PP A_1$ is the space of curves with a marked point and a marked direction, such that 
the curve has a node at the marked point and the second derivative along the marked direction vanishes. 
Note that there are two such distinguished directions. Hence, the projection map  
$\pi : \PP A_1 \lra A_1$ 
is two to one. 
Similarly, the maps 
\[ \pi: \PP D_4 \lra D_4 \qquad \textnormal{and} \qquad  \pi: \PP X_9 \lra X_9  \]
are three to one and four to one respectively. Let us also define 
\begin{align*}
\hat{\mathfrak{X}} &:= \{ ([f], l_q) \in \DD \times \mathbb{P}T\X: ([f], q) \in \mathfrak{X} \} = \pi^{-1}(\mathfrak{X}).
\end{align*}
Finally, if $S_1, \ldots S_n$ are subsets of $\DD \times X$ and 
$T$ is a subset of $\DD \times \mathbb{P}TX$, then we define 
\begin{align*}
S_1 \circ S_2 \circ \ldots \circ S_n\circ T &:= \{ ([f], q_1, \ldots, q_n, l_{q_{n+1}}) \in 
\DD \times X^n \times \mathbb{P} TX:  
([f], q_1) \in S_1, \ldots, ([f], q_n) \in S_n, \\ 
%& ,  \\ %\qquad \forall i = 1 ~~\textnormal{to} ~~n, \\ 
& \qquad \qquad \qquad \qquad \qquad \qquad \qquad \qquad \qquad \qquad ([f], l_{q_{n+1}}) \in T \qquad \textnormal{and}\\ 
& \qquad \qquad \qquad \qquad \qquad \qquad  \qquad \qquad \qquad \qquad q_1, \ldots, q_n, q_{n+1} 
\qquad \textnormal{are all distinct}\}.
\end{align*}
For example, $A_1^{2} \circ \PP A_2$ is the space of curves with 
three distinct ordered points, where the curve has a simple node at the first two points and a cusp at the 
last point and a distinguished direction at the last marked point, 
such that the Hessian vanishes along that direction.   
%and all the three points are distinct. 
%For example, $A_1^2 \circ \PP A_2$ is the 
%we define $A_1^{\delta} \circ \PP \X$ and $A_1^{\delta} \circ \hat{\X}$  to be as follows 
%\begin{align*}
%A_1^{\delta} \circ \PP \X &:= \{ ([f], q_1, \ldots, q_{\delta}, l_{q_{\delta+1}}) 
%\in \DD \times (\X)^{\delta} \times \P T\X: f ~~\textnormal{has a signularity of type $A_1$ at 
%$q_1, \ldots q_{\delta}$ and} \\ 
%& \qquad \qquad \qquad \qquad \qquad \qquad \qquad  ([f], l_{q_{\delta +1}}) \in \PP \mathfrak{X} 
%~~\textnormal{and $q_1 \ldots q_{\delta+1}$ are all distinct} \}. 
%\end{align*}
We are now in a position to define a few numbers. Let us define 
%Equipped with this notation, we can define the following numbers:
\begin{align}
\Num(A_1^{\delta} \mathfrak{\X} , n_1, m_1, m_2) &:= \langle y^{\delta_L-k-n_1-m_1-2m_2-\delta} c_1^{m_1} x_1^{m_1} x_2^{m_2}, 
~~[\overline{A_1^{\delta} \circ \mathfrak{X}}] \rangle,\\
\Num(A_1^{\delta} \PP \mathfrak{\X} , n_1, m_1, m_2, \theta) &:= 
\langle y^{\delta_L-k-n_1-m_1-2m_2-\delta - \theta} c_1^{m_1} x_1^{m_1} x_2^{m_2} \lambda^{\theta}, 
~~[\overline{A_1^{\delta} \circ \PP \mathfrak{X}}] \rangle \qquad \textnormal{and}  \\ 
\Num(A_1^{\delta} \hat{\mathfrak{\X}} , n_1, m_1, m_2, \theta) &:= 
\langle y^{\delta_L-k-n_1-m_1-2m_2-\delta - \theta} c_1^{m_1} x_1^{m_1} x_2^{m_2} \lambda^{\theta}, 
~~[\overline{A_1^{\delta} \circ \hat{\mathfrak{X}}}] \rangle, 
%\qquad \textnormal{where}
%\qquad \textnormal{and} \qquad 
\end{align}
where 
\begin{align*}
c_1 &: = c_1(L), ~~x_i := c_i(T^*\X), ~~\lambda := c_1(\hat{\gamma}^*), ~~ y:= c_1(\gamma_{\D}^*) 
\end{align*}
and $\gamma_{\D} \lra \D$ and $\hat{\gamma} \lra \mathbb{P} T\X$ are 
the tautological line bundles. \\
\hf \hf Let us now make a few observations. First, let us abbreviate 
$\Num(A_1^{\delta} \mathfrak{\X} , 0, 0, 0)$ as $\Num(A_1^{\delta} \mathfrak{\X})$; 
this quantity  
is the number of curves passing through the 
right number of points that have $\delta$ distinct nodes and one singularity of type $\mathfrak{X}$. 
%Hence, let us abbreviate $\Num(A_1^{\delta} \mathfrak{\X} , 0, 0, 0)$ as $\Num(A_1^{\delta} \mathfrak{\X})$.
Next, we note that 
\begin{align}
\Num(\A_1^{\delta} \mathfrak{X}, n_1, m_1, m_2) & = \frac{1}{\textnormal{deg}(\pi)} 
\Num(\A_1^{\delta} \PP \mathfrak{X}, n_1, m_1, m_2, 0),   \label{deg_to_one_up_to_down}
\end{align}
where $\textnormal{deg}(\pi)$ is the degree of the projection map $\pi:\PP \mathfrak{X}\lra\mathfrak{X}$. 
We remind the reader that 
the degree is always one, except in the cases when $\mathfrak{X} = A_1, D_4$ or $X_9$, in which 
case the degree is $2$, $3$ and $4$ respectively. \\ 
\hf \hf Next, we note that using the ring structure of  
$H^*(\mathbb{P} T\X, \mathbb{Z})$ (cf \cite{BoTu}, pp. 270), 
we conclude that 
\begin{align}
\Num(A_1^{\delta} \PP \mathfrak{X}, n_1, m_1, m_2, \theta) & = 
\Num(A_1^{\delta} \PP \mathfrak{X}, n_1, m_1+1, m_2, \theta-1) 
-\Num(A_1^{\delta} \PP \mathfrak{X}, n_1, m_1, m_2+1, \theta-2) \nonumber \\ 
& \qquad \qquad \textnormal{if} \qquad \theta>1. \label{proj_tang_space_ring_eqn}
\end{align}
Finally, we note that 
\begin{align}
\Num(A_1^{\delta} \hat{\mathfrak{X}}, n_1, m_1, m_2, \theta) & = 0 \qquad \textnormal{if} \qquad \theta =0, \nonumber \\ 
\Num(A_1^{\delta} \hat{\mathfrak{X}}, n_1, m_1, m_2, \theta) & = \Num(A_1^{\delta} \mathfrak{X}, n_1, m_1, m_2) 
\qquad \textnormal{if} \qquad \theta =1 \qquad 
\textnormal{and} \nonumber \\ 
\Num(A_1^{\delta} \hat{\mathfrak{X}}, n_1, m_1, m_2, \theta) & = 
\Num(A_1^{\delta} \hat{\mathfrak{X}}, n_1, m_1+1, m_2, \theta-1) 
-\Num(A_1^{\delta} \hat{\mathfrak{X}}, n_1, m_1, m_2+1, \theta-2) \nonumber \\ 
& \qquad \qquad \textnormal{if} \qquad \theta>1. \label{deg_to_one_up_to_down_hat}
\end{align}
This again follows from \cite{BoTu}, pp. 270.

\section{Recursive Formulas}
\label{recursive_formulas}
%\ni The following  are the main theorems of this paper. 
%In particular, they 
%give us a formula for $\Num(A_1^{\delta} \mathfrak{X})$ if $\delta+k \leq 8$. 
%Before stating the main results of this paper, we 
%note that for any given singularity, we have the following equality: 
%\begin{align}
%\Num(A_1^{\delta} \PP \mathfrak{X}, n_1, m_1, m_2, \theta) & = 
%\Num(A_1^{\delta} \PP \mathfrak{X}, n_1, m_1+1, m_2, \theta-1) 
%-\Num(A_1^{\delta} \PP \mathfrak{X}, n_1, m_1, m_2+1, \theta-2). \label{proj_tang_space_ring_eqn}
%\end{align}
%This follows from using the ring structure of $H^*(\mathbb{P} T\X, \mathbb{Z})$ 
%(cf \cite{BoTu}, pp. 270).
%We also note that 
%\begin{align}
%\Num(\A_1^{\delta} \mathfrak{X}, n_1, m_1, m_2) & = \frac{1}{\textnormal{deg}(\pi)} 
%\Num(\A_1^{\delta} \PP \mathfrak{X}, n_1, m_1, m_2, 0),   \label{deg_to_one_up_to_down}
%\end{align}
%where $\textnormal{deg}(\pi)$ is the degree of the projection map $\pi:\PP \mathfrak{X}\lra\mathfrak{X}$. 
%We remind the reader that 
%the degree is always one, except in the cases when $\mathfrak{X} = A_1, D_4$ or $X_9$, in which 
%case the degree is $2$, $3$ and $4$ respectively. Let us also abbreviate 
%$\Num(A_1^{\delta} \mathfrak{X}, 0,0,0)$ as $\Num(A_1^{\delta} \mathfrak{X})$.\\  
\hf \hf We are now ready to state the recursive formulas.  
We are going to give a recursive formula for curves with 
$\delta$-nodes and one singularity of type $\mathfrak{X}$ in terms of 
number of curves with $\delta-1$ nodes and one singularity more degenerate than 
$\mathfrak{X}$. Combined with the base case of the recursion (Theorem \ref{base_case_recursion}), 
equations \eqref{deg_to_one_up_to_down}, \eqref{proj_tang_space_ring_eqn}  
and \eqref{deg_to_one_up_to_down_hat}, 
Theorems \ref{a1_delta_a1} to \ref{a1_delta_pa8}  
enable us to obtain an explicit formula for 
$\Num(A_1^{\delta} \mathfrak{X})$ if $\delta+k \leq 8$. \\
\hf \hf We have written a mathematica program to implement these recursive formulas and 
obtain the final formulas. The program is available on the second author's homepage 
\[ \textnormal{\url{https://www.sites.google.com/site/ritwik371/home}} \]
For the convenience of the reader, we 
have generated an output of all the numbers 
$\Num(A_1^{\delta} \mathfrak{X})$
and have appended it at the end of this paper 
(after the Bibliography). \\
%include in the appendix the output of the 
%program that includes a list of all the explicit formulas. 
%We are now ready to state the main results of this paper.  
\hf \hf We remind the reader that 
in the formulas for $N(A_1^{\delta} \mathfrak{X})$, the nodes are always \textit{ordered}. 
To obtain the corresponding number of curves where the nodes are unordered, we have to divide by 
$(\delta+1)!$ or $\delta!$, depending on whether the last singularity ($\mathfrak{X}$) is $A_1$ 
or something different, respectively.

\begin{thm}
\label{base_case_recursion}
The number $\N(\A_1, n_1, m_1, m_2)$ is given by 
\begin{align}
\label{a1_chern_eqn}
\N(\A_1, n_1, m_1, m_2) & = \begin{cases} 3 c_1^2 + 2 c_1 x_1 + x_2 & \mbox{if} ~~(n_1, m_1, m_2) =(0,0,0), \\ 
   3c_1^2 + c_1 x_1  & \mbox{if }  ~~(n_1, m_1, m_2) = (1,0,0),\\ 
   c_1^2 & \mbox{if } ~~(n_1, m_1, m_2) = (2,0,0), \\ 
   3 c_1 x_1 + x_1^2 &  \mbox{if} ~~~(n_1, m_1, m_2) = (0,1,0), \\ 
   c_1 x_1 & \mbox{if} ~~~(n_1, m_1, m_2) = (1,1,0), \\ 
   x_1^2 & \mbox{if} ~~~(n_1, m_1, m_2) = (0,2,0), \\ 
   x_2 & \mbox{if} ~~~(n_1, m_1, m_2) = (0,0,1), \\ 
   0 & \mbox{otherwise}.\end{cases} 
\end{align}
\end{thm}

\begin{rem}
Theorem \ref{base_case_recursion} will serve as the base case of our recursion.  
\end{rem}

\begin{thm}
\label{a1_delta_a1}
If $1 \leq \delta \leq 7$, then 
\begin{align*}
\N(\A_1^{\delta} \A_1 , n_1, m_1, m_2) & = \N(\A_1^{\delta}, 0, 0, 0) \times \N(\A_1, n_1,m_1, m_2)  \\ 
                                       & - \Big(\binom{\delta}{1}\Big( \N(\A_1^{\delta -1}, n_1,m_1, m_2) + 
                                       \N(\A_1^{\delta -1}, n_1+1, m_1, m_2)  + 
                                       3 \N(\A_1^{\delta -1} \A_2, n_1,m_1, m_2) \Big) \\ 
                                        & + 4\binom{\delta}{2} \N(\A_1^{\delta -2} \A_3, n_1,m_1, m_2) +  
                                          18 \binom{\delta}{3}\N(\A_1^{\delta -3} \D_4, n_1,m_1, m_2) \Big),
\end{align*}
provided the line bundle is sufficiently $(2 \delta +1)$-ample.  
\end{thm}

\begin{thm}
\label{a1_delta_pa1}
If $0 \leq \delta \leq 7$, then $\N(\A_1^{\delta} \PP \A_1 , n_1, m_1, m_2, \theta)$ is as follows:
\begin{align*}
\theta=0: &\hf 2 \N(\A_1^{\delta} \A_1, n_1, m_1, m_2) \\
\theta=1: &\hf 2 \N(\A_1^{\delta} \A_1, n_1, m_1+1, m_2) + \N(\A_1^{\delta} \A_1, n_1, m_1, m_2) +\N(\A_1^{\delta} \A_1, n_1+1, m_1, m_2) \\ 
                &\hf 
                -2 \binom{\delta}{2} \N(\A_1^{\delta-2} \PP \D_4, n_1, m_1, m_2, \theta-1),
                %,\\
                %- j_1 \binom{\delta}{5}\N(A_1^{\delta -5}X_9,n_1,m_1,m_2),\\
%\theta\geq 2:&\hf \N(\A_1^{\delta} \PP \A_1 , n_1, m_1+1, m_2, \theta-1) -\N(\A_1^{\delta} \PP \A_1 , n_1, m_1, m_2+1, \theta-2) ,
\end{align*}
provided the line bundle is sufficiently $(2 \delta +1)$-ample.  
\end{thm}

\begin{thm}
\label{a1_delta_pa2}
If $0 \leq \delta \leq 6$, then 
%$\N(\A_1^{\delta} \PP \A_2 , n_1, m_1, m_2, \theta)$ is given by 
%as follows:
\begin{align*}
\N(\A_1^{\delta} \PP \A_2 , n_1, m_1, m_2, \theta) 
&=%2\N(\A_1^{\delta} \A_1 , n_1, m_1, m_2) + 2\N(\A_1^{\delta} \A_1 , n_1, m_1+1, m_2 ) 
%+ 2\N(\A_1^{\delta} \A_1 , n_1+1, m_1, m_2),
%, \\ 
%&\hf - 2 \binom{\delta}{1}\N(\A_1^{\delta-1} \PP \A_3, n_1, m_1, m_2,\theta) 
%  -4 \binom{\delta}{2}\N(\A_1^{\delta-2} \PP \D_4, n_1, m_1, m_2, \theta)  \\ 
%\theta=1:
\N(\A_1^{\delta} \PP \A_1 , n_1, m_1, m_2, \theta) + \N(\A_1^{\delta} \PP \A_1 , n_1, m_1+1, m_2, \theta ) \\
               &\hf + \N(\A_1^{\delta} \PP \A_1 , n_1+1, m_1, m_2, \theta) - 2 \binom{\delta}{1}\N(\A_1^{\delta-1} \PP \A_3, n_1, m_1, m_2, \theta)\\
               &\hf - 3\binom{\delta}{1}\N(\A_1^{\delta-1} \hat{\D}_4, n_1, m_1, m_2, \theta) 
               -4 \binom{\delta}{2}\N(\A_1^{\delta-2} \PP \D_4, n_1, m_1, m_2, \theta) \\
               &\hf -4\binom{\delta}{2}\N(\A_1^{\delta-2} \hat{\D}_5, n_1, m_1, m_2, \theta), 
               %-j_2\binom{\delta}{5}\N(A_1^{\delta -5}X_9,n_1,m_1,m_2)\\
%\theta\geq 2:&\hf \N(\A_1^{\delta} \PP \A_2 , n_1, m_1+1, m_2, \theta-1) 
%-\N(\A_1^{\delta} \PP \A_2 , n_1, m_1, m_2+1, \theta-2) ,
\end{align*}
provided the line bundle is sufficiently $(2 \delta +2)$-ample.  
\end{thm}

%\begin{rem}
%Note that 
%\begin{align}
%\Num(\A_1^{\delta} \hat{\mathfrak{X}}, n_1, m_1, m_2,\theta) & = 0 \qquad \textnormal{if} \qquad 
%\theta =0 \qquad \textnormal{and} \nonumber \\
%\Num(\A_1^{\delta} \hat{\mathfrak{X}}, n_1, m_1, m_2,\theta) & = 
%\Num(\A_1^{\delta} \mathfrak{X}, n_1, m_1, m_2) \qquad \textnormal{if} \qquad \theta =1. \label{hat_eqn}
%\end{align}
%************
%In order to compute $\Num(\A_1^{\delta} \PP \A_2 , n_1, m_1, m_2, \theta)$ using 
%Theorem \ref{a1_delta_pa2}, 
%we use \eqref{deg_to_one_up_to_down_hat} for $\theta =0$ or $1$ and \eqref{proj_tang_space_ring_eqn}, 
%when $\theta>1$. 
%\end{rem}

\begin{thm}
\label{a1_delta_pa3}
If $0 \leq \delta \leq 5$, then 
\begin{align*}
\N(\A_1^{\delta} \PP \A_3 , n_1, m_1, m_2, \theta)  & = 3 \N(\A_1^{\delta} \PP \A_2 , n_1, m_1, m_2, \theta+1) 
                                                      + \N(\A_1^{\delta} \PP \A_2 , n_1+1, m_1, m_2, \theta) \\
                                                     & +\N(\A_1^{\delta} \PP \A_2 , n_1, m_1, m_2, \theta) -2 \binom{\delta}{1}\N(\A_1^{\delta-1} 
                                                      \PP \A_4, n_1, m_1, m_2, \theta) \\
                                                     & -2\binom{\delta}{2}\N(\A_1^{\delta-2} 
                                                      \PP \D_5, n_1, m_1, m_2, \theta),
\end{align*}
provided the line bundle is sufficiently $(2 \delta +3)$-ample.
\end{thm}

\begin{thm}
\label{a1_delta_pa4}
If $0 \leq \delta \leq 4$, then 
\begin{align*}
\N(\A_1^{\delta} \PP \A_4 , n_1, m_1, m_2, \theta) & = 2\N(\A_1^{\delta} \PP \A_3 , n_1, m_1, m_2, \theta+1)
                                                                                      +2\N(\A_1^{\delta} \PP \A_3 , n_1, m_1+1, m_2, \theta)\\
                                                                                      &+2\N(\A_1^{\delta} \PP \A_3 , n_1, m_1, m_2, \theta)
                                                                                      +2\N(\A_1^{\delta} \PP \A_3 , n_1+1, m_1, m_2, \theta)\\
                                                                                      &-2\binom{\delta}{1}\N(\A_1^{\delta-1} \PP \A_5 , n_1, m_1, m_2, \theta)
                                                                                      -4\binom{\delta}{2}\N(\A_1^{\delta-2} \PP \D_6 , n_1, m_1, m_2, \theta),
\end{align*}
provided the line bundle is sufficiently $(2 \delta +4)$-ample.
\end{thm}

\begin{thm}
\label{a1_delta_pd4}
If $0 \leq \delta \leq 4$, then $\N(\A_1^{\delta} \PP \D_4 , n_1, m_1, m_2, \theta)$ is as follows:
\begin{align*}
\theta=0: &\hf 2\N(\A_1^{\delta} \PP \A_3 , n_1, m_1+1, m_2, \theta) 
              -2\N(\A_1^{\delta} \PP \A_3 , n_1, m_1, m_2, \theta+1)\\
                &\hf +\N(\A_1^{\delta} \PP \A_3 , n_1, m_1, m_2, \theta)+
                \N(\A_1^{\delta} \PP \A_3 , n_1+1, m_1, m_2, \theta)\\
                &\hf -2\binom{\delta}{1}\N(\A_1^{\delta-1} \D_5 , n_1, m_1, m_2)  
                -2\binom{\delta}{2}\N(\A_1^{\delta-2} \PP \D_6 , n_1, m_1, m_2, 0)\\
\theta=1: &\hf \textstyle{\frac{1}{3}}\big(\N(\A_1^{\delta} \PP \D_4 , n_1, m_1, m_2, \theta-1) 
+3\N(\A_1^{\delta} \PP \D_4 , n_1, m_1+1, m_2,\theta-1)\\
                &\hf +\N(\A_1^{\delta} \PP \D_4 , n_1+1, m_1, m_2, \theta-1)\big) 
                -24\binom{\delta}{3}\Num(X_9, n_1, m_1, m_2)\\
\theta\geq 2: &\hf \N(\A_1^{\delta} \PP \D_4 , n_1, m_1+1, m_2, \theta-1) -\N(\A_1^{\delta} \PP \D_4 , n_1, m_1, m_2+1, \theta-2),
\end{align*}
provided the line bundle is sufficiently $(2 \delta +4)$-ample.
\end{thm}

\begin{rem}
Note that the structure of this formula is a little bit different; we have something different when 
$\theta=0$ and when $\theta=1$. The reason for this will be clear when we prove the Theorem 
(subsection \ref{proof_of_a1_delta_pd4}).
\end{rem}

\begin{thm}
\label{a1_delta_pd5}
If $0 \leq \delta \leq 3$, then $\N(\A_1^{\delta} \PP \D_5 , n_1, m_1, m_2, \theta)$ is as follows:
\begin{align*}
\N(\A_1^{\delta} \PP \D_5 , n_1, m_1, m_2, \theta) &= \hf \N(\A_1^{\delta} \PP \D_4 , n_1, m_1, m_2, \theta+1) +\N(\A_1^{\delta} \PP \D_4 , n_1, m_1+1, m_2, \theta)\\
& +\N(\A_1^{\delta} \PP \D_4 , n_1, m_1, m_2, \theta)+\N(\A_1^{\delta} \PP \D_4 , n_1+1, m_1, m_2, \theta)\\
& -2\binom{\delta}{1} \N(\A_1^{\delta-1} \PP \D_6 , n_1, m_1, m_2, \theta) 
 -12\binom{\delta}{1} \N(\A_1^{\delta-1} \hat{X}_9 , n_1, m_1, m_2, \theta) \\
& -18\binom{\delta}{3}\N(\PP X_9, n_1, m_1, m_2, \theta),
%\theta=1: &\hf \N(\A_1^{\delta} \PP \D_4 , n_1, m_1, m_2, \theta+1) +\N(\A_1^{\delta} \PP \D_4 , n_1, m_1+1, m_2, \theta)\\
%& +\N(\A_1^{\delta} \PP \D_4 , n_1, m_1, m_2, \theta)+\N(\A_1^{\delta} \PP \D_4 , n_1+1, m_1, m_2, \theta)\\
%& -2\binom{\delta}{1} \N(\A_1^{\delta-1} \PP \D_6 , n_1, m_1, m_2, \theta) - q_6\binom{\delta}{3}\N(A_1^{\delta -3}X_9,n_1,m_1,m_2)\\
%& - \delta_{\delta 2}q_2\N(X_9,n_1,m_1,m_2)\\
%\theta\geq 2: &\hf \N(\A_1^{\delta} \PP \D_5 , n_1, m_1+1, m_2, \theta-1) -\N(\A_1^{\delta} \PP \D_5 , n_1, m_1, m_2+1, \theta-2),
\end{align*}
provided the line bundle is sufficiently $(2 \delta +5)$-ample.
\end{thm}

\begin{thm}
\label{a1_delta_pa5}
If $0 \leq \delta \leq 3$, then 
\begin{align*}
\N(\A_1^{\delta} \PP \A_5 , n_1, m_1, m_2, \theta) &= 3\N(\A_1^{\delta} \PP \A_4 , n_1, m_1, m_2, \theta+1) +2 \N(\A_1^{\delta} \PP \A_4 , n_1, m_1+1, m_2, \theta)\\
& +2\N(\A_1^{\delta} \PP \A_4 , n_1, m_1, m_2, \theta)+2\N(\A_1^{\delta} \PP \A_4 , n_1+1, m_1, m_2, \theta)\\
& - 2\binom{\delta}{1} \N(\A_1^{\delta-1} \PP \A_6 , n_1, m_1, m_2, \theta) -\binom{\delta}{1} \N(\A_1^{\delta-1} \PP \E_6 , n_1, m_1, m_2, \theta)\\
& -4\binom{\delta}{2}\N(\A_1^{\delta-2} \PP \D_7 , n_1, m_1, m_2, \theta),
%- \delta_{\delta 2}\delta_{\theta 1}v_6 \N(X_9,n_1,m_1,m_2),
\end{align*}
provided the line bundle is sufficiently $(2 \delta +5)$-ample. 
%The above formula is valid even when $\delta=3$ and $\theta=0$.
\end{thm}

\begin{thm}
\label{a1_delta_pd6}
If $0 \leq \delta \leq 2$, then 
\begin{align*}
\N(\A_1^{\delta} \PP \D_6 , n_1, m_1, m_2, \theta) &= 4\N(\A_1^{\delta} \PP \D_5 , n_1, m_1, m_2, \theta+1) + \N(\A_1^{\delta} \PP \D_5 , n_1, m_1, m_2, \theta)\\
& +\N(\A_1^{\delta} \PP \D_5 , n_1+1, m_1, m_2, \theta) - 2\binom{\delta}{1}\N(\A_1^{\delta-1} \PP \D_7 , n_1, m_1, m_2, \theta)\\
& -\binom{\delta}{1}\N(\A_1^{\delta-1} \PP \E_7 , n_1, m_1, m_2, \theta), 
%- \delta_{\delta 2}\delta_{\theta 0}q_4 \N(\X_9,n_1,m_1,m_2),
\end{align*}
provided the line bundle is sufficiently $(2 \delta +6)$-ample.  
%The above formula is valid even when $\delta=2$ and $\theta=0$.
\end{thm}

\begin{thm}
\label{a1_delta_pe6}
If $0 \leq \delta \leq 2$, then $\N(\A_1^{\delta} \PP \E_6 , n_1, m_1, m_2, \theta)$ is as follows:
\begin{align*}
\N(\A_1^{\delta} \PP \E_6 , n_1, m_1, m_2, \theta) &= 2\N(\A_1^{\delta} \PP \D_5 , n_1, m_1+1, m_2, \theta) - \N(\A_1^{\delta} \PP \D_5 , n_1, m_1, m_2, \theta+1)\\
& +\N(\A_1^{\delta} \PP \D_5 , n_1, m_1, m_2, \theta)+\N(\A_1^{\delta} \PP \D_5 , n_1+1, m_1, m_2, \theta)\\
& - \binom{\delta }{1}\N(\A_1^{\delta-1} \PP \E_7 , n_1, m_1, m_2, \theta) 
  -4 \binom{\delta}{1}\N(\hat{X}_9,n_1,m_1,m_2,\theta), 
%\theta=1: &\hf 2\N(\A_1^{\delta} \PP \D_5 , n_1, m_1+1, m_2, \theta) - \N(\A_1^{\delta} \PP \D_5 , n_1, m_1, m_2, \theta+1)\\
%& +\N(\A_1^{\delta} \PP \D_5 , n_1, m_1, m_2, \theta)+\N(\A_1^{\delta} \PP \D_5 , n_1+1, m_1, m_2, \theta)\\
%& - \binom{\delta }{1}\N(\A_1^{\delta-1} \PP \E_7 , n_1, m_1, m_2, \theta) - \delta_{\delta 1}v_2 \N(X_9,n_1,m_1,m_2)\\
%\theta\geq 2: &\hf \N(\A_1^{\delta} \PP \E_6 , n_1, m_1+1, m_2, \theta-1) -\N(\A_1^{\delta} \PP \E_6 , n_1, m_1, m_2+1, \theta-2),
\end{align*}
provided the line bundle is sufficiently $(2 \delta +6)$-ample.  
%The above formula is valid even when $\delta=2$ and $\theta=0$.
\end{thm}

\begin{thm}
\label{a1_delta_pa6}
If $0 \leq \delta \leq 2$, then $\N(\A_1^{\delta} \PP \A_6 , n_1, m_1, m_2, \theta)$ is as follows:
\begin{align*}
\N(\A_1^{\delta} \PP \A_6 , n_1, m_1, m_2, \theta) & =  
2 \N(\A_1^{\delta} \PP \A_5 , n_1, m_1, m_2, \theta+1) +4\N(\A_1^{\delta} \PP \A_5 , n_1, m_1+1, m_2, \theta)\\
& +3\N(\A_1^{\delta} \PP \A_5 , n_1, m_1, m_2, \theta) +3 \N(\A_1^{\delta} \PP \A_5 , n_1+1, m_1, m_2, \theta)\\
& - 2\N(\A_1^{\delta} \PP \D_6 , n_1, m_1, m_2, \theta) - \N(\A_1^{\delta} \PP \E_6 , n_1, m_1, m_2, \theta)\\
& - 2\binom{\delta}{1}\N(\A_1^{\delta-1} \PP \A_7 , n_1, m_1, m_2, \theta) - 3\binom{\delta}{1}\N(\A_1^{\delta-1} \PP \E_7 , n_1, m_1, m_2, \theta)\\
& -12\binom{\delta}{1}\N(\A_1^{\delta-1} \hat{X}_9 , n_1, m_1, m_2, \theta)-6 \binom{\delta}{2}  \N(\PP \D_8,n_1,m_1,m_2), 
%\\
%\theta=1: &\hf 2 \N(\A_1^{\delta} \PP \A_5 , n_1, m_1, m_2, \theta+1) +4\N(\A_1^{\delta} \PP \A_5 , n_1, m_1+1, m_2, \theta)\\
%& +3\N(\A_1^{\delta} \PP \A_5 , n_1, m_1, m_2, \theta) +3 \N(\A_1^{\delta} \PP \A_5 , n_1+1, m_1, m_2, \theta)\\
%& - 2\N(\A_1^{\delta} \PP \D_6 , n_1, m_1, m_2, \theta) - \N(\A_1^{\delta} \PP \E_6 , n_1, m_1, m_2, \theta)\\
%& - 2\binom{\delta}{1}\N(\A_1^{\delta-1} \PP \A_7 , n_1, m_1, m_2, \theta) - 3\binom{\delta}{1}\N(\A_1^{\delta-1} \PP \E_7 , n_1, m_1, m_2, \theta)\\
%& -\delta_{\delta 1}v_3\N(\X_9, n_1, m_1, m_2)\\
%\theta\geq 2: &\hf \N(\A_1^{\delta} \PP \A_6 , n_1, m_1+1, m_2, \theta-1) -\N(\A_1^{\delta} \PP \A_6 , n_1, m_1, m_2+1, \theta-2),
\end{align*}
provided the line bundle is sufficiently $(2 \delta +6)$-ample. 
%The above formula is valid even when $\delta=2$ and $\theta=0$.
\end{thm}

\begin{thm}
\label{a1_delta_pd7}
If $0\leq \delta\leq 1$, then 
\begin{align*}
\N(\A_1^{\delta} \PP \D_7 , n_1, m_1, m_2, \theta) &= 4 \N(\A_1^{\delta} \PP \D_6 , n_1, m_1, m_2, \theta+1) +2\N(\A_1^{\delta} \PP \D_6 , n_1, m_1+1, m_2, \theta)\\
& +2\N(\A_1^{\delta} \PP \D_6 , n_1, m_1, m_2, \theta) +2 \N(\A_1^{\delta} \PP \D_6 , n_1+1, m_1, m_2, \theta)\\
& -2 \binom{\delta}{1} \N(\PP \D_8, n_1, m_1, m_2,\theta),
\end{align*}
provided the line bundle is sufficiently $(2 \delta +7)$-ample. 
\end{thm}

\begin{thm}
\label{a1_delta_pe7}
If $0\leq \delta\leq 1$, then 
\begin{align*}
\N(\A_1^{\delta} \PP \E_7 , n_1, m_1, m_2, \theta) &= \N(\A_1^{\delta} \PP \D_6 , n_1, m_1, m_2, \theta) 
-\N(\A_1^{\delta} \PP \D_6 , n_1, m_1, m_2, \theta+1)\\
& 2\N(\A_1^{\delta} \PP \D_6 , n_1, m_1+1, m_2, \theta)+\N(\A_1^{\delta} \PP \D_6 , n_1+1, m_1, m_2, \theta)\\
& -3 \binom{\delta}{1} \N(\PP X_9, n_1,m_1,m_2, \theta), \\
\end{align*}
provided the line bundle is sufficiently $(2 \delta +3)$-ample. 
\end{thm}

\begin{thm}
\label{a1_delta_pa7}
If $0\leq \delta\leq 1$, then 
%$\N(\A_1^{\delta} \PP \A_7 , n_1, m_1, m_2, \theta)$ is as follows:
\begin{align*}
\N(\A_1^{\delta} \PP \A_7 , n_1, m_1, m_2, \theta) & = 
-\N(\A_1^{\delta} \PP \A_6 , n_1, m_1, m_2, \theta+1) +8\N(\A_1^{\delta} \PP \A_6 , n_1, m_1+1, m_2, \theta)\\
& +5\N(\A_1^{\delta} \PP \A_6 , n_1, m_1, m_2, \theta) +5 \N(\A_1^{\delta} \PP \A_6 , n_1+1, m_1, m_2, \theta) \\
& -5\N(A_1^{\delta} \hat{\X}_9, n_1,m_1,m_2, \theta) \\
& - 6\N(\A_1^{\delta} \PP \D_7 , n_1, m_1, m_2, \theta)- 7\N(\A_1^{\delta} \PP \E_7 , n_1, m_1, m_2, \theta)\\
& -2\binom{\delta}{1}\N(\PP \A_8, n_1,m_1,m_2,\theta)-14 \binom{\delta}{1} \N(\PP \E_8, n_1,m_1,m_2,\theta), 
%& *, \\
%\theta=1: & \hf - \N(\A_1^{\delta} \PP \A_6 , n_1, m_1, m_2, \theta+1) +8\N(\A_1^{\delta} \PP \A_6 , n_1, m_1+1, m_2, \theta)\\
%& +5\N(\A_1^{\delta} \PP \A_6 , n_1, m_1, m_2, \theta) +5 \N(\A_1^{\delta} \PP \A_6 , n_1+1, m_1, m_2, \theta)\\
%& - 6\N(\A_1^{\delta} \PP \D_7 , n_1, m_1, m_2, \theta)- 7\N(\A_1^{\delta} \PP \E_7 , n_1, m_1, m_2, \theta)\\
%& -5 \N(\X_9, n_1,m_1,m_2)\\
%\theta\geq 2:& \hf \N(\A_1^{\delta} \PP \A_7 , n_1, m_1+1, m_2, \theta-1) -\N(\A_1^{\delta} \PP \A_7 , n_1, m_1, m_2+1, \theta-2), 
\end{align*}
provided the line bundle is sufficiently $(2 \delta +7)$-ample. 
%The above formula is valid even when $\delta=1$ and $\theta=0$.
\end{thm}

\begin{thm}
\label{a1_delta_pe8}
The number $\N(\PP \E_8 , n_1, m_1, m_2, \theta)$ is given by
\begin{align*}
\N(\PP \E_8 , n_1, m_1, m_2, \theta) &= \N(\PP \E_7, n_1, m_1, m_2,\theta)+\N(\PP \E_7, n_1, m_1+1, m_2,\theta)+3\N(\PP \E_7, n_1+1, m_1, m_2,\theta)\\
& +2\N(\PP \E_7, n_1, m_1, m_2,\theta+1)
\end{align*}
provided the line bundle is sufficiently $8$-ample. 
\end{thm}

\begin{thm}
\label{a1_delta_pd8}
The number $\N(\PP \D_8 , n_1, m_1, m_2, \theta)$ is given by
\begin{align*}
\N(\PP \D_8 , n_1, m_1, m_2, \theta) &= 4\N(\PP \D_7, n_1, m_1, m_2,\theta)+6\N(\PP \D_7, n_1, m_1+1, m_2,\theta)+3\N(\PP \A_7, n_1, m_1, m_2,\theta+1)\\
& 4\N(\PP \D_7, n_1+1, m_1, m_2,\theta)-3\N(\PP \E_8, n_1, m_1, m_2,\theta)
\end{align*}
provided the line bundle is sufficiently $8$-ample. 
\end{thm}

\begin{thm}
\label{a1_delta_px9}
The number $\N(\PP X_9 , n_1, m_1, m_2, \theta)$ is given by
\begin{align*}
\N(\PP X_9 , n_1, m_1, m_2, \theta) &=  \N(\PP \E_7, n_1, m_1, m_2,\theta)
+4\N(\PP \E_7, n_1, m_1, m_2,\theta+1)+\N(\PP \E_7, n_1+1, m_1, m_2,\theta)
\end{align*}
provided the line bundle is sufficiently $4$-ample. 
\end{thm}

It remains to give a formula for $\Num(\PP A_8, n_1, m_1, m_2, \theta)$. 
%and $\Num(X_9)$. 
Before we do that, we need to define 
one more space. 
%In this paper, we will directly define the closure  
\begin{align*}
\PP X_9^{\vee} &:= \{ ([f], l_{q}) \in \DD \times \mathbb{P} TX:  ([f], l_{q}) \in \hat{X}_9, 
          ~~\Big(- \frac{f_{31}^3}{8 } + \frac{3 f_{22} f_{31} f_{40}}{16 } - \frac{f_{13} f_{40}^2}{16}  \Big) = 0\}. 
\end{align*}

\begin{thm}
\label{a1_delta_px9vee}
The number $\N(\PP X_9^{\vee} , n_1, m_1, m_2, \theta)$ is given by
\begin{align*}
\N(\PP X_9^{\vee} , n_1, m_1, m_2, \theta) &= 3 \N(\hat{X}_9, n_1, m_1, m_2, \theta) 
                                            + 3\N(\hat{X}_9, n_1, m_1+1, m_2, \theta)\\ 
                                           & +6\N(\hat{X}_9, n_1, m_1, m_2, \theta+1)+
                                            3\N(\hat{X}_9, n_1+1, m_1, m_2, \theta),
\end{align*}
provided the line bundle is sufficiently $4$-ample. 
\end{thm}

\begin{thm}
\label{a1_delta_pa8}
The number $\N(\PP \A_8 , n_1, m_1, m_2)$ is given by
\begin{align*}
\N(\PP \A_8 , n_1, m_1, m_2, \theta)&= 6 \N(\PP \A_7, n_1, m_1, m_2) -2 \N(\PP \A_7, n_1, m_1, m_2, 1)+ 10 \N(\PP \A_7, n_1, m_1+ 1, m_2)\\
& +6 \N(\PP \A_7, n_1+ 1, m_1, m_2) -8 \N (\PP \D_8, n_1, m_1, m_2) - 16 \N (\PP \E_8, n_1, m_1, m_2)\\
& - 6\N(\PP X_9^{\vee}, n_1, m_1, m_2, \theta),
\end{align*}
provided the line bundle is sufficiently $8$-ample. 
\end{thm}

%-------------------------------------------------------------- Euler class -----------------------------------------------------------------------------------
\section{Proof of the Recursive Formulas via Euler class computation} 
\hf\hf We will now give a proof of the recursive formulas we stated in 
section \ref{recursive_formulas}. The relevant claims on 
closure, multiplicity and transversality are proved in our earlier 
papers \cite{BM13_one_singular_point_published}, \cite{BM13_2pt_published}, 
\cite{RM_Hypersurfaces_published} and \cite{BM_linear_system_codim_seven}, when 
$\delta =0$ or $1$ and $\delta+k \leq 7$. For the remaining cases 
(namely $\delta>1$ or $\delta+k>7$), we will prove the relevant 
closure, multiplicity and transversality claims in 
\cite{BM_closure_of_seven_points}. 
%the results stated in the introduction.

\subsection{Proof of Theorem \ref{base_case_recursion}: computing $\N(\A_1, n_1, m_1, m_2)$} 
\hf \hf This is proved in \cite{BM_linear_system_codim_seven}, Proposition 5.2. \qed

\subsection{Proof of Theorem \ref{a1_delta_a1}: computing $\N(\A_1^{\delta} \A_1, n_1, m_1, m_2)$}
\hf\hf Recall that in section \ref{notation} we have defined the space  
\begin{align*}
\A^{\delta}_1 \circ (\DD \times X) := \{ ([f], q_1, \ldots, q_{\delta}, q_{\delta+1}) \in \DD \times (X)^{\delta +1}: 
&\textnormal{$f$ has a singularity of type $\A_1$ at $q_1, \ldots, q_{\delta}$}, \\ 
                   & \textnormal{$q_1, \ldots, q_{\delta+1}$ all distinct}\}. 
\end{align*}
Let $\mu$ be a generic pseudocycle representing the homology class 
Poincar\'{e} dual 
to 
\[c_1^{n_1} x_1^{m_1} x_2^{m_2} y^{\delta_L - (n_1 + m_1 + 2 m_2 + \delta+1)}.\]
We now consider sections of the following two bundles 
that are induced by the evaluation map and the vertical derivative, namely: 
\begin{align*}
\psi_{\A_0}: \A^{\delta}_1 \circ (\DD \times X) \cap \mu \lra \mathcal{L}_{\A_0} & := 
\gamma_{\DD}^* \otimes L, 
\qquad \qquad  
\{\psi_{\A_0}([f], q_1, \ldots, q_{\delta+1})\}(f) := f(q_{\delta+1}) \qquad  \nonumber \\
\psi_{\A_1} : \psi_{\A_0}^{-1}(0) \lra \mathcal{V}_{\A_1} &:= \gamma_{\D}^*\otimes T^*\X \otimes L, 
\qquad 
\{\psi_{\A_1}([f], q_1, \ldots, q_{\delta+1})\}(f) := \nabla f|_{q_{\delta+1}}. 
%\label{psi_a0_a1_section_defn}
\end{align*}
%In \cite{RM_Hypersurfaces} we show that 
In \cite{BM_closure_of_seven_points}, we show that if 
$L$ is sufficiently $(2\delta+1)$-ample, then 
these sections are transverse to the zero set. 
Next, let us define 
\begin{align*}
\mathcal{B} &:= \overline{\A^{\delta}_1 \circ (\DD \times X)}- \A^{\delta}_1 \circ (\DD \times X). 
\end{align*}
Hence 
\begin{align*}
\lan e(\mathcal{L}_{\A_0}) e(\mathcal{V}_{\A_1}), 
~~[\overline{\A^{\delta}_1 \circ (\DD \times X)}] \cap [\mu] \ran & = \N(\A_1^{\delta}\A_1, n_1, m_1, m_2)  
+ \mathcal{C}_{\mathcal{B}\cap \mu},
%\label{na1_pseudocycle}
\end{align*}
where $\mathcal{C}_{\mathcal{B}\cap \mu}$ denotes the contribution of the section to the Euler 
class from the points of $\mathcal{B}\cap \mu$.
We now give an explicit description of the boundary $\mathcal{B}\cap \mu$; 
the proof will be given in \cite{BM_closure_of_seven_points}.  
First of all we note that the whole of 
$\mathcal{B}$ is not relevant while computing the contribution to the Euler class; only the points 
at which the section vanishes is relevant. Hence, we only need to consider the component of 
$\mathcal{B}\cap \mu$ where one (or more) of the $q_i$ become equal to the last point $q_{\delta+1}$. 
Let us define 
\begin{align*}
\B (q_1, q_{\delta +1 }) &:= \{ ([f], q_1, \ldots, q_{\delta +1}) \in \B: q_1 = q_{\delta +1} \}.
\end{align*}
%Geometrically, $\B (q_1, q_{\delta +1 })$ denotes the component of the boundary when the 
%point $q_1$ collides with the last marked point $q_{\delta+1}$. 
The spaces $\B (q_2, q_{\delta +1 }), \ldots, \B (q_{\delta}, q_{\delta +1 })$ are defined 
similarly. 
%Each of these stratum can be identified with a copy of $\overline{\A_1^{\delta}}$. 
It is shown in \cite[Corollary 6.6]{BM13_2pt_published} 
%(Corollary 6.6) 
and \cite[Proposition 6.1]{BM_linear_system_codim_seven} 
%(Proposition 6.1) 
that the contribution from 
$\B (q_1, q_{\delta +1 }) \cap \mu$ is 
\begin{align*}
\N(\A_1^{\delta -1}, n_1,m_1, m_2) + 
                                       \N(\A_1^{\delta -1}, n_1+1, m_1, m_2)  + 
                                       3 \N(\A_1^{\delta -1} \A_2, n_1,m_1, m_2) .
\end{align*}
The same contribution arises from $\B (q_i, q_{\delta +1 }) \cap \mu$ for 
all $i = 1$ to $\delta$. 
Hence the total contribution from $\B (q_1, q_{\delta +1 }) \cap \mu, \ldots, \B (q_{\delta}, q_{\delta +1 }) \cap \mu$ 
is 
\begin{align*}
\delta \Big(\N(\A_1^{\delta -1}, n_1,m_1, m_2) + 
                                       \N(\A_1^{\delta -1}, n_1+1, m_1, m_2)  + 
                                       3 \N(\A_1^{\delta -1} \A_2, n_1,m_1, m_2) \Big).
\end{align*}
Next, let us define 
\begin{align*}
\B (q_1, q_2, q_{\delta +1 }) &:= \{ ([f], q_1, \ldots, q_{\delta +1}) \in \B: q_1 = q_2 = q_{\delta +1} \}.
\end{align*}
The spaces $\B (q_{i_1}, q_{i_2}, q_{\delta +1 })$, and more generally 
$\B (q_{i_1}, q_{i_2},\ldots,q_{i_j}, q_{\delta +1 })$, are defined similarly. 
Let us now focus on $\B (q_1, q_2, q_{\delta +1 })$. 
We claim that 
\begin{align}
\B (q_1, q_2, q_{\delta +1 }) & \approx \overline{\A_1^{\delta-2} \circ \A}_3.  \label{two_nodes_tacnode_claim}
\end{align}
Let us give an intuitive justification as to why this is a believable claim.  
Geometrically, $\B (q_1, q_2, q_{\delta +1 })$ denotes the component of the boundary when the 
nodal points $q_1$ and $q_2$ collide with each other. Geometrically, we expect the following thing 
to happen: 
\newpage 
\begin{figure}[h!]
\vspace*{0.2cm}
\begin{center}\includegraphics[scale = 0.5]{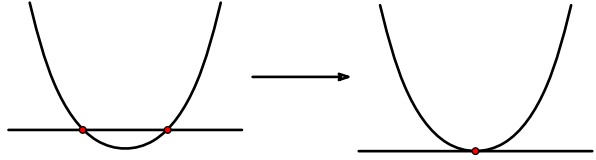}\vspace*{-0.2cm}\end{center}
\caption{Two nodes colliding into a tacnode}
\end{figure} 
%\newpage 
That is indeed the case. We give a rigorous proof of \eqref{two_nodes_tacnode_claim} 
in \cite[Lemma 6.3(2)]{BM13_2pt_published}.
We will show in \cite{BM_closure_of_seven_points} 
that 
%and 
%\begin{figure}[h!]
%\label{fig_3_nodes_triple_pt}
%\vspace*{0.2cm}
%\begin{center}\includegraphics[scale = 0.5]{3nodestriple.eps}\vspace*{-0.2cm}\end{center}
%\caption{Three nodes colliding into a triple point}
%\end{figure}
%and 
%\begin{figure}[h!]
%\label{fig_3_nodes_a5}
%\vspace*{0.2cm}
%\begin{center}\includegraphics[scale = 0.5]{3nodesA5.eps}\vspace*{-0.2cm}\end{center}
%\caption{Three nodes colliding into an $\A_5$-node}
%\end{figure}
the contribution from each of the tacnodal points, namely  $\B (q_1, q_2, q_{\delta +1 }) \cap \mu$ 
is $4$. Hence the total contribution from all the components 
of type $\B (q_{i_1}, q_{i_2}, q_{\delta +1 })$ is 
\begin{align*}
4 \binom{\delta}{2} \N(\A_1^{\delta} \A_3, n_1, m_1, m_2).   
\end{align*}
Next, let us focus on $\B (q_1, q_2, q_3, q_{\delta +1 })$.
We claim that 
\begin{align}
\B (q_1, q_2, q_3, q_{\delta +1 })\cap \mu & \approx \overline{\A_1^{\delta-3} \circ \D_4}\cap \mu.  
\label{three_nodes_collide_to_triple_point_claim}
\end{align}
Again, let us give an intuitive justification as to why this is a believable claim.  
Geometrically, $\B (q_1, q_2, q_3, q_{\delta +1 })$ denotes the component of the boundary when three of  
the nodal points $q_1$, $q_2$ and $q_3$ collide with each other. 
Geometrically, we expect the following thing 
to happen: 
\newpage 
\begin{figure}[h!]
\label{fig_3_nodes_triple_pt}
\vspace*{0.2cm}
\begin{center}\includegraphics[scale = 0.5]{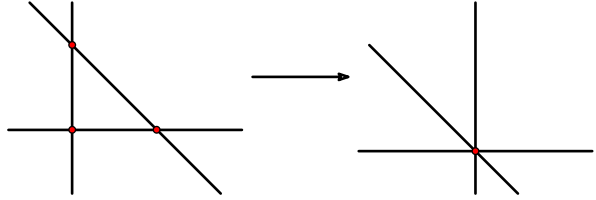}\vspace*{-0.2cm}\end{center}
\caption{Three nodes colliding into a triple point}
\end{figure}
%and 
%\begin{figure}[h!]
%\label{fig_3_nodes_a5}
%\vspace*{0.2cm}
%\begin{center}\includegraphics[scale = 0.5]{3nodesA5.eps}\vspace*{-0.2cm}\end{center}
%\caption{Three nodes colliding into an $\A_5$-node}
%\end{figure}
That is indeed the case. We will give a rigorous proof of \eqref{three_nodes_collide_to_triple_point_claim} 
in \cite{BM_closure_of_seven_points}. We will also show that  
the contribution from each of the points of $\B (q_1, q_2, q_3, q_{\delta +1 }) \cap \mu$ 
is $18$. Hence the total contribution from all the components 
of type $\B (q_{i_1}, q_{i_2}, q_{i_3}, q_{\delta +1 })$ is 
\begin{align*}
18 \binom{\delta}{3} \N(\A_1^{\delta} \D_4, n_1, m_1, m_2).   
\end{align*}
Finally, we will also show in \cite{BM_closure_of_seven_points} that if we consider 
a one dimensional family of $r$-nodal curves, then all the $r$-nodes 
can not collide together if $4 \leq  r \leq 7$. 
%when more than 
%three nodes collide, the 
In other words, 
%claim that 
\begin{align}
\B(q_1, q_2, \ldots, q_r, q_{\delta+1}) \cap \mu & = \varnothing \qquad \forall ~ 4 \leq r \leq 7.  
\label{r_nodal_come_together}
\end{align}
Hence, the remaining stratum of $\B$ does not contribute to the Euler class,
%, since its intersection with $\mu$ 
%is empty. 
giving us Theorem \ref{a1_delta_a1}.\qed 

\begin{rem}
Note that \eqref{r_nodal_come_together} is true because we intersected the left hand side 
%and right hand side 
with $\mu$. It is \textit{not} true that 
\begin{align*}
\B(q_1, q_2, \ldots, q_r, q_{\delta+1}) & = \varnothing \qquad \forall ~ 4 \leq r \leq 7.  
%\label{r_nodal_come_together}
\end{align*}
Equation \eqref{r_nodal_come_together} is saying that when more than three nodes 
come together, the resulting singularity will have such a high codimension that it will not 
intersect a generic codimension one class (i.e. if we have 
a one dimensional family of curves, then more than three nodes can not come together). 
Similarly, \eqref{three_nodes_collide_to_triple_point_claim} is true only with the $\mu$ 
present; it is \textit{not} true that 
\begin{align}
\B (q_1, q_2, q_3, q_{\delta +1 })\cap \mu & \approx \overline{\A_1^{\delta-3} \circ \D_4}.  
%\label{three_nodes_collide_to_triple_point_claim}
\end{align}
Geometrically, the following thing can also occur: 
\begin{figure}[h!]
\label{fig_3_nodes_a5}
\vspace*{0.2cm}
\begin{center}\includegraphics[scale = 0.5]{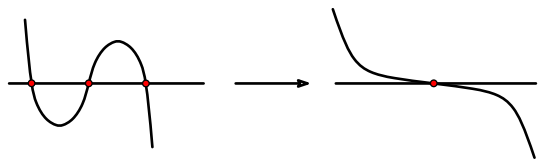}\vspace*{-0.2cm}\end{center}
\caption{Three nodes colliding into an $\A_5$-singularity}
\end{figure}
\newpage 
Hence, three nodes can collide to form an $A_5$-singularity as well. However, since 
this is a codimension $5$ singularity, its intersection with $\mu$ will be empty; only 
codiemnsion four singularities will survive. All this will be proven rigorously in 
\cite{BM_closure_of_seven_points}. 
\end{rem}

\subsection{Proof of Theorem \ref{a1_delta_pa1}: computing $\N(\A_1^{\delta} \PP \A_1, n_1, m_1, m_2, \theta)$}
When $\theta =0$, the formula 
follows from \eqref{deg_to_one_up_to_down}. 
%holds because the map $\pi: A_1^{\delta} \circ \PP A_1 \lra A_1^{\delta+1}$ 
%is two to one. 
Let us now assume $\theta>0$.\\ 
\hf \hf Recall that in section \ref{notation} we have defined the space  
\begin{align*}
\A^{\delta}_1 \circ \overline{\hat{\A}}_1 := \{ ([f], q_1, \ldots, q_{\delta}, 
l_{q_{\delta+1}}) \in \DD \times (X)^{\delta} \times \mathbb{P} T \X: 
&\textnormal{$f$ has a singularity of type $\A_1$ at $q_1, \ldots, q_{\delta}$}, \\ 
                   & ([f], l_{q_{\delta+1}}) \in \overline{\hat{A}}_1, ~~\textnormal{$q_1, \ldots, q_{\delta+1}$ all distinct}\}. 
\end{align*}
Let $\mu$ be a generic pseudocycle representing the homology class 
Poincar\'{e} dual 
to 
\[c_1^{n_1} x_1^{m_1} x_2^{m_2} \lambda^{\theta}y^{\delta_L - (n_1 + m_1 + 2 m_2 +\theta+\delta+2)}.\]
We now define a section of the following bundle 
\begin{align*}
\psi_{\PP \A_1}:  \A^{\delta}_1 \circ \overline{\hat{\A}}_1 
\lra \mathcal{L}_{\PP \A_1} & := \gamma_{\D}^\ast\otimes\hat{\gamma}^{\ast 2}\otimes L, \qquad \textnormal{given by}\\
\{\psi_{\PP \A_1}([f], q_1,\ldots,q_\delta,l_{q_{\delta+1}})\}(f\otimes v^{\otimes 2}) &:= \nabla^2 f|_{q_{\delta+1}}(v, v).  \nonumber
%\psi_{\A_1} : \psi_{\A_0}^{-1}(0) \lra \mathcal{V}_{\A_1} &:= \gamma_{\D}^*\otimes T^*\X \otimes L, 
%\qquad \textnormal{given by} \qquad 
%\{\psi_{\A_1}([f], q)\}(f) := \nabla f|_q. \label{psi_a0_a1_section_defn}
\end{align*}
Here $\hat{\gamma} \lra \P TX$ is the tautological line bundle over $\P TX$. 
%In \cite{RM_Hypersurfaces} we show that 
We will show in \cite{BM_closure_of_seven_points} that 
if $L$ is sufficiently $(2\delta+2)$-ample, then 
this section is transverse to the zero set. 
Next, let us define 
\begin{align*}
\mathcal{B} &:= \overline{\A^{\delta}_1 \circ \overline{\hat{\A}}}_1- \A^{\delta}_1 \circ \overline{\hat{\A}}_1. 
\end{align*}
Hence 
\begin{align}
\lan e(\mathcal{L}_{\PP \A_1}), 
~~[\overline{\A^{\delta}_1 \circ \overline{\hat{\A}}}_1] \cap [\mu] \ran & = 
\N(\A_1^{\delta}\PP \A_1, n_1, m_1, m_2, \theta)  
+ \mathcal{C}_{\mathcal{B}\cap \mu},
%\label{na1_pseudocycle}
\end{align}
where as before, $\mathcal{C}_{\mathcal{B} \cap \mu}$ denotes the contribution of the section 
to the Euler class from $\mathcal{B} \cap \mu$. 
We now give an explicit description of $\mathcal{B}$. 
%First of all we note that the whole of 
%$\mathcal{B}$ is not relevant while computing the contribution to the Euler class; only the points 
%at which the section vanishes is relevant. 
As before, we only need to consider the component of 
$\mathcal{B}$ where one (or more) of the $q_i$ become equal to the last point $q_{\delta+1}$. 
Let us define 
\begin{align*}
\B (q_1, l_{q_{\delta +1 }}) &:= \{ ([f], q_1, \ldots, q_{\delta}, l_{q_{\delta +1}}) \in \B: q_1 = q_{\delta +1} \}.
\end{align*}
The spaces $\B (q_2, q_{\delta +1 }), \ldots, \B (q_{\delta}, q_{\delta +1 })$ are defined 
similarly. 
%Each of these stratum can be identified with a copy of $\overline{\A_1^{\delta}}$. 
%Next, let us define 
%\begin{align*}
%\B (q_1, q_2, q_{\delta +1 }) &:= \{ ([f], q_1, \ldots, q_{\delta +1}) \in \B: q_1 = q_2 = q_{\delta +1} \}.
%\end{align*}
%The spaces $\B (q_{i_1}, q_{i_2}, q_{\delta +1 })$ are defined similarly. 
Following the same argument  as in \eqref{two_nodes_tacnode_claim}, we conclude that 
%The claim is that 
\begin{align*}
\B (q_1, l_{q_{\delta +1 }}) & \approx \overline{\A_1^{\delta-1} \circ \hat{\A}}_3.  
\end{align*}
Let us now intersect $\overline{\A_1^{\delta-1} \circ \hat{\A}}_3$ with $\mu$. This will 
be an isolated set of finite points. Hence, the section $\psi_{\PP \A_1}$ 
will not vanish on $\overline{\A_1^{\delta-1} \circ \hat{\A}}_3 \cap \mu$.  
%The section 
%The section does not vanish on $\B (q_1, l_{q_{\delta +1 }}) \cap \mu$. 
Hence it does not 
contribute to the Euler class. 
Next, let us define 
\begin{align*}
\B (q_1, q_2, l_{q_{\delta +1 }}) &:= 
\{ ([f], q_1, \ldots, q_{\delta}, l_{q_{\delta +1}}) \in \B: q_1=q_2 = l_{q_{\delta +1}} \}.
\end{align*}
The spaces $\B (q_{i_1}, q_{i_2}, \ldots,q_{i_j},l_{q_{\delta +1 }})$ are defined similarly. 
Following the same argument  as in \eqref{three_nodes_collide_to_triple_point_claim}, we conclude that 
%We claim that 
\begin{align*}
\B (q_1, q_2, l_{q_{\delta +1 }}) & \approx \overline{\A_1^{\delta-2} \circ \hat{\D}}_4.  
\end{align*}
The section $\psi_{\PP \A_1}$ vanishes everywhere on 
$\overline{\A_1^{\delta-2} \circ \hat{\D}}_4$; hence it also vanishes on 
$\overline{\A_1^{\delta-2} \circ \hat{\D}}_4 \cap \mu$. We will show in \cite{BM_closure_of_seven_points} 
that the 
contribution from each of the points of $\B (q_1, q_2, l_{q_{\delta +1 }}) \cap \mu$  
is $6$. Hence the total contribution from all the components 
of type $\B (q_{i_1}, q_{i_2}, l_{q_{\delta +1 }})$ is 
%non-zero only when $\theta=1$ and equals
\begin{align*}
6 \binom{\delta}{2} \N(\A_1^{\delta-2} \hat{\D}_4, n_1, m_1, m_2, \theta).   
\end{align*}
Finally, using \eqref{r_nodal_come_together} 
we conclude that 
\begin{align*}
\B(q_1, q_2, \ldots, q_r, l_{q_{\delta+1}}) \cap \mu & = \varnothing \qquad \forall ~ 3 \leq r \leq 7.  
\end{align*}
Hence the remaining stratum of $\B \cap \mu$ does not contribute 
to the Euler class,
%, since its intersection with $\mu$ 
%is empty. 
giving Theorem \eqref{a1_delta_pa1}. \qed 

\begin{rem}
The reason we need to compute $\N(\A_1^{\delta}\PP A_1, n_1, m_1, m_2, \theta)$ is because 
it arises in the recursive formula for 
$\N(\A_1^{\delta} \PP \A_2, n_1, m_1, m_2, \theta)$, as will be clear shortly. 
\end{rem}

\subsection{Proof of Theorem \ref{a1_delta_pa2}: computing $\N(\A_1^{\delta} \PP \A_2, n_1, m_1, m_2, \theta)$}
\hf \hf Recall that in section \ref{notation}
we have defined the space  
\begin{align*}
\A^{\delta}_1 \circ \overline{\PP \A}_1 := \{ ([f], q_1, \ldots, q_{\delta}, 
l_{q_{\delta+1}}) \in \DD \times (X)^{\delta} \times \mathbb{P} T \X: 
&\textnormal{$f$ has a singularity of type $\A_1$ at $q_1, \ldots, q_{\delta}$}, \\ 
                   & ([f], l_{q_{\delta+1}}) \in \overline{\PP \A}_1, ~~\textnormal{$q_1, \ldots, q_{\delta+1}$ all distinct}\}. 
\end{align*}
Let $\mu$ be a generic pseudocycle representing the homology class 
Poincar\'{e} dual 
to 
\[c_1^{n_1} x_1^{m_1} x_2^{m_2} \lambda^{\theta}y^{\delta_L - (n_1 + m_1 + 2 m_2 +\theta+\delta+2)}.\]
We now define a section of the following line bundle 
\begin{align}
\Psi_{\PP \A_2}:  \A^{\delta}_1 \circ \overline{\PP \A}_1 
\lra \mathbb{L}_{\PP \A_2} & := \gamma_{\D}^\ast\otimes \hat{\gamma}^\ast \otimes (\pi^*T X/\hat{\gamma})^\ast\otimes L, 
\qquad \textnormal{given by} \nonumber \\
\{\Psi_{\PP \A_2}([f], q_1,\ldots,q_\delta,l_{q_{\delta+1}})\}(f\otimes v \otimes w) &:= 
\nabla^2 f|_{q_{\delta+1}}(v,w). \nonumber
%\psi_{\A_1} : \psi_{\A_0}^{-1}(0) \lra \mathcal{V}_{\A_1} &:= \gamma_{\D}^*\otimes T^*\X \otimes L, 
%\qquad \textnormal{given by} \qquad 
%\{\psi_{\A_1}([f], q)\}(f) := \nabla f|_q. \label{psi_a0_a1_section_defn}
\end{align}
%In \cite{RM_Hypersurfaces} we show that 
If $L$ is sufficiently $(2\delta+2)$-ample, then this section is transverse to zero. 
Next, let us define 
\begin{align*}
\mathcal{B} &:= \overline{\A^{\delta}_1 \circ \overline{\PP \A}}_1- \A^{\delta}_1 \circ \overline{\PP\A}_1. 
\end{align*}
Hence 
\begin{align*}
\lan e(\mathbb{L}_{\PP \A_2}), 
~~[\overline{\A^{\delta}_1 \circ \overline{\PP\A}}_1] \cap [\mu] \ran & = \N(\A_1^{\delta}\PP \A_2, n_1, m_1, m_2, \theta)  
+ \mathcal{C}_{\mathcal{B}\cap \mu}.
%\label{na1_pseudocycle}
\end{align*}
We now give an explicit description of $\mathcal{B}$. 
%First of all we note that the whole of 
%$\mathcal{B}$ is not relevant while computing the contribution to the Euler class; only the points 
%at which the section vanishes is relevant. 
As before, we only need to consider the component of 
$\mathcal{B}$ where one (or more) of the $q_i$ become equal to the last point $q_{\delta+1}$. 
Define $\mathcal{B}(q_{i_1}, \ldots q_{i_k}, l_{q_{\delta}})$ as before. 
%Let us define 
%\begin{align*}
%\B (q_1, l_{q_{\delta +1 }}) &:= \{ ([f], q_1, \ldots, q_{\delta}, l_{q_{\delta +1}}) \in \B: q_1 = q_{\delta +1} \}.
%\end{align*}
%The spaces $\B (q_2, q_{\delta +1 }), \ldots, \B (q_{\delta}, q_{\delta +1 })$ are defined 
%similarly. 
%Each of these stratum can be identified with a copy of $\overline{\A_1^{\delta}}$. 
%Next, let us define 
%\begin{align*}
%\B (q_1, q_2, q_{\delta +1 }) &:= \{ ([f], q_1, \ldots, q_{\delta +1}) \in \B: q_1 = q_2 = q_{\delta +1} \}.
%\end{align*}
%The spaces $\B (q_{i_1}, q_{i_2}, q_{\delta +1 })$ are defined similarly. 
We will show in \cite{BM_closure_of_seven_points} that   
\begin{align*}
\B (q_1, l_{q_{\delta +1 }}) & \approx \overline{\A_1^{\delta-1} \circ \PP \A_3}\cup 
\overline{\A_1^{\delta-1} \circ \hat{\D}_4}.
\end{align*}
%The argument is similar to how we prove *. 
Furthermore, we show that the contribution from 
$\overline{\A_1^{\delta-1} \circ \PP \A_3} \cap \mu$ is $2$, while the contribution from 
$\overline{\A_1^{\delta-1} \circ \hat{\D}_4} \cap \mu$ is $3$. 
%The points of the second type do not contribute when $\theta=0$. 
%For $\theta=1$ the contribution from each of the points of $\B (q_1, l_{q_{\delta +1 }}) \cap \mu$ of the 
%first type and second type are $2$ and $3$ respectively. 
Hence the total contribution from all 
the components of type $\B (q_{i_1}, l_{q_{\delta +1 }})$ equals
\bgd
2\binom{\delta}{1}N(\A_1^{\delta-1}\PP\A_3,n_1,m_1,m_2,\theta)
+3\binom{\delta}{1}N(\A_1^{\delta-1}\hat{\D}_4,n_1,m_1,m_2, \theta).
\edd
%when $\theta=1$. 
%Next, let us define 
%\begin{align*}
%\B (q_1, q_2, l_{q_{\delta +1 }}) &:= 
%\{ ([f], q_1, \ldots, q_{\delta}, l_{q_{\delta +1}}) \in \B: q_1=q_2 = l_{q_{\delta +1}} \}.
%\end{align*}
Next, we show in \cite{BM_closure_of_seven_points} that
\begin{align*}
\B (q_1, q_2, l_{q_{\delta +1 }}) & \approx 
\overline{\A_1^{\delta-2} \circ \PP \D}_4 \cup\overline{\A_1^{\delta-2}\circ \hat{\D}}_5.  
\end{align*}
%The points of the first type do not contribute when $\theta=0$. 
We also show that 
%claim that 
the contribution from 
$\overline{\A_1^{\delta-2} \circ \hat{\D}_4} \cap \mu$ is $2$, while the contribution from 
$\overline{\A_1^{\delta-2} \circ \PP \D_5}$ is $3$. 
%For $\theta=1$ 
%The contribution from each of the points of $\B (q_1, q_2, l_{q_{\delta +1 }}) \cap \mu$ of 
%the first and second type are $4$ in both cases.
Hence the total contribution from all the components of type $\B (q_{i_1}, q_{i_2}, l_{q_{\delta +1 }})$ 
%when $\theta=1$ 
equals
\begin{align*}
4 \binom{\delta}{2} \N(\A_1^{\delta-2} \PP \D_4, n_1, m_1, m_2, \theta)+
4 \binom{\delta}{2}  \N(\A_1^{\delta-2} \hat{\D}_5, n_1, m_1, m_2,\theta).   
\end{align*}
Finally, we show in \cite{BM_closure_of_seven_points} that 
\begin{align}
\B(q_1, q_2, \ldots, q_r, l_{q_{\delta+1}}) \cap \mu & = \varnothing \qquad \forall ~ 3 \leq r \leq 5.  
\label{five_nodes_one_directional_node_special_quadruple_point}
\end{align}
Hence the remaining stratum of $\B$ does not contribute to the Euler class, 
giving us Theorem \ref{a1_delta_pa2}. \qed 

\begin{rem}
Although we have devoted \cite{BM_closure_of_seven_points} for 
the rigorous proofs of our closure claims, there is very  
subtle point about 
equation \eqref{five_nodes_one_directional_node_special_quadruple_point} 
that we should mention here. 
Let us look at the claim for $r=3$. This is essentially saying that four nodes can not 
come together, once we intersect the variety with $\mu$. Let us see if this claim is 
believable; consider the following  picture: 
%aa \\ 
%\newline PICTURE. \\
\newpage 
\begin{figure}[h!]
\vspace*{0.2cm}
\begin{center}\includegraphics[scale = 0.5]{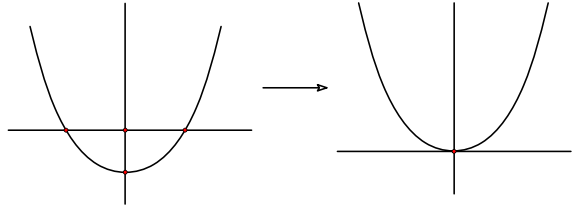}\vspace*{-0.2cm}\end{center}
\caption{Four nodes collapsing to a D6-singularity}
\end{figure} 
%\subsection{TT}
\noindent Now it seems naively that one of the components of    
$\B(q_1, q_2, q_3, l_{q_4})$ is  $\hat{D}_6$. If that were 
indeed the case then the intersection with $\mu$ would not be zero when 
$\theta =1$ (it would be zero when $\theta =0$). Hence, what we are 
claiming here is that 
one of the components of    
$\B(q_1, q_2, \ldots, q_3, l_{q_4})$ is  \textit{not} $\hat{D}_6$, but is in fact  
$\PP D_6$. Let us look at it a bit more carefully. We are \textit{not} considering 
here the space of curves with four distinct nodes. We are in fact considering the 
space of curves with three nodes and one node with a distinguished branch of the 
node; namely the space $\A_1^3 \circ \PP A_1$ as opposed to $\A_1^4$. Hence, any 
random curve with a $D_6$-singularity will not be in the closure; i.e. all elements of 
$\hat{D}_6$ are not in the closure. Only those elements will be in the closure where a 
a specific directional derivative vanishes, namely only elements of $\PP D_6$ will be 
in the closure. Hence, \eqref{five_nodes_one_directional_node_special_quadruple_point} 
is consistent with the picture we have drawn. \\ 
\hf \hf Next, let us look at the claim for $r=5$. 
This is essentially saying that six nodes can not 
come together, once we intersect the variety with $\mu$.
Again, let us see if this claim is 
believable by looking at the following  picture: \\ 
\newpage 
\begin{figure}[h!]
\vspace*{0.2cm}
\begin{center}\includegraphics[scale = 0.5]{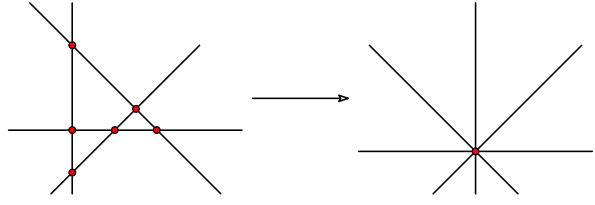}\vspace*{-0.2cm}\end{center}
\caption{Six nodes collapsing to an ordinary quadruple point}
\end{figure} 
%\newline PICTURE. \\
\noindent Now it seems naively that one of the components of    
$\B(q_1, q_2, q_3, q_4, q_5, l_{q_6})$ is  $\hat{X}_9$. If that were 
indeed the case then the intersection with $\mu$ would not be zero when 
$\theta =1$ (it would be zero when $\theta =0$). Hence, what we are 
claiming here is that 
one of the components of    
$\B(q_1, q_2, \ldots, q_5, l_{q_6})$ is  \textit{not} $\hat{X}_9$, but is in fact  
$\PP X_9$. Let us look at it a bit more carefully. 
Again, we are \textit{not} considering 
here the space of curves with six distinct nodes. We are in fact considering the 
space of curve with five nodes and one node with a distinguished branch of the 
node; namely the space $\A_1^5 \circ \PP A_1$ as opposed to $\A_1^6$. Hence, any 
random curve with an $X_9$ singularity will not be in the closure; i.e. all elements of 
$\hat{X}_9$ are not in the closure. Only those elements will be in the closure where a 
a specific directional derivative vanishes, namely only elements of $\PP X_9$ will be 
in the closure. Hence, \eqref{five_nodes_one_directional_node_special_quadruple_point} 
is consistent with the picture we have drawn.   
%\begin{align}
%\B(q_1, q_2, \ldots, q_3, l_{q_4}) & = \hat{D}_6 \cup \ldots   
%\label{five_nodes_one_directional_node_special_quadruple_point}
%\end{align}
%Till $r =4$, 
%\eqref{five_nodes_one_directional_node_special_quadruple_point} is geometrically very believable; 
%it is saying that if we have a one dimensional family of curves with upto five nodes, then all the 
%five nodes can not come together. 
%There is a very 
%till $\delta =4$. 
%Finally, we show in 
%\cite{BM_closure_of_seven_points} that 
%\begin{align*}
%\B(q_1, q_2, \ldots, q_5, l_{q_{6}}) & = \PP X_9. %\qquad \forall ~ 3 \leq r \leq 4.  
%\end{align*}
%Since $\PP X_9 \cap \mu = \varnothing$, we conclude that 
%since its intersection with $\mu$ 
%is empty. That proves 
\end{rem}

\subsection{Proof of Theorem \ref{a1_delta_pa3}: computing $\N(\A_1^{\delta} \PP \A_3, n_1, m_1, m_2, \theta)$}
\hf \hf Recall that 
in section \ref{notation}
we have defined the space  
\begin{align*}
\A^{\delta}_1 \circ \overline{\PP \A}_2 := \{ ([f], q_1, \ldots, q_{\delta}, 
l_{q_{\delta+1}}) \in \DD \times (X)^{\delta} \times \mathbb{P} T \X: 
&\textnormal{$f$ has a singularity of type $\A_1$ at $q_1, \ldots, q_{\delta}$}, \\ 
                   & ([f], l_{q_{\delta+1}}) \in \overline{\PP \A}_2, 
                   ~~\textnormal{$q_1, \ldots, q_{\delta+1}$ all distinct}\}. 
\end{align*}
Let $\mu$ be a generic pseudocycle representing the homology class 
Poincar\'{e} dual 
to 
\[c_1^{n_1} x_1^{m_1} x_2^{m_2} \lambda^{\theta}y^{\delta_L - (n_1 + m_1 + 2 m_2 +\theta+\delta+3)}.\]
We now define a section of the following bundle 
\begin{align}
\Psi_{\PP \A_3}:  \A^{\delta}_1 \circ \overline{\PP \A}_2 
\lra \mathbb{L}_{\PP \A_3} & := \gamma_{\D}^\ast\otimes\hat{\gamma}^{\ast 3}\otimes L, \qquad \textnormal{given by} 
\nonumber \\
\{\Psi_{\PP \A_3}([f], q_1,\ldots,q_\delta,l_{q_{\delta+1}})\}(f\otimes v^{\otimes 3}) &:= 
\nabla^3f|_{q_{\delta+1}}(v,v,v).  \nonumber
%\psi_{\A_1} : \psi_{\A_0}^{-1}(0) \lra \mathcal{V}_{\A_1} &:= \gamma_{\D}^*\otimes T^*\X \otimes L, 
%\qquad \textnormal{given by} \qquad 
%\{\psi_{\A_1}([f], q)\}(f) := \nabla f|_q. \label{psi_a0_a1_section_defn}
\end{align}
%In \cite{RM_Hypersurfaces} we show that 
%The variety $\A^{\delta}_1 \circ \overline{\PP \A}_2$ is smooth; this follows from the 
%transversality of the section $\Psi_{\PP \A_2}$ defined in the proof of Theorem  \ref{a1_delta_pa3}. 
%Let us define 
%\begin{align*}
%M &:= \{ ([f], q_1, \ldots, q_{\delta}, l_{q_{\delta+1}}) \in A_1^{\delta} \circ \overline{\PP A}_2: 
%\nabla^2 f|_{q_{\delta+1}}(w,w) \neq 0\}.  
%\end{align*}
%Note that 
By 
%Lemma $7.1$, statement $8$ of the arXiv version 
\cite[Lemma 7.1, statement 8]{BM13}, we have that 
\begin{align}
\overline{\PP A}_2 &= \PP A_2 \cup \overline{\PP A}_3 \cup \overline{\hat{D}}_4. \label{pa2_closure}
\end{align}
Hence, let us define 
\begin{align*}
\mathcal{B} &:= \overline{A_1^{\delta} \circ \overline{\PP A}}_2 - \A^{\delta}_1 \circ (\PP\A_2\cup \overline{\PP A}_3). 
\end{align*}
If $L$ is sufficiently $(2\delta+3)$-ample, then the section 
$\Psi_{\PP \A_3}$ vanishes on the points of  
$\A_1^{\delta}\circ \PP A_3$ transversally 
(all the points of $\A_1^{\delta}\circ \PP A_3$ are smooth points of the 
variety $\A_1^{\delta}\circ \overline{\PP A}_2$).
%$M$ is transverse to the zero set. 
%Hence, the number of zeros of the section $\Psi_{\PP A_3}$ restricted to 
%$\A_1^{\delta}\circ \overline{\PP \A}_2 $
%$M \cap \mu$ is 
%the number $\N(\A_1^{\delta}\PP \A_3, n_1, m_1, m_2, \theta)$.
Hence,  
\begin{align*}
\lan e(\mathbb{L}_{\PP \A_3}), 
~~[\overline{\A^{\delta}_1 \circ \overline{\PP\A}}_2] \cap [\mu] \ran & = \N(\A_1^{\delta}\PP \A_3, n_1, m_1, m_2, \theta)  
+ \mathcal{C}_{\mathcal{B}\cap \mu}. \nonumber 
%\label{na1_pseudocycle}
\end{align*}
We now give an explicit description of $\mathcal{B}$. 
Let us first define 
\begin{align*}
\mathcal{B}_0 &:= \{ ([f], q_1, \ldots q_{\delta}, l_{q_{\delta+1}}) \in \mathcal{B}: q_1, q_2 \ldots q_{\delta+1} 
~~\textnormal{are all distinct}\}. 
\end{align*}
In other words, $\mathcal{B}_0$ is that component of the boundary, where all the points are still distinct.   
By \eqref{pa2_closure}, we conclude that 
\begin{align*}
\mathcal{B}_0 &= \overline{A_1^{\delta}} \circ \overline{\hat{D}}_4. 
\end{align*}
If we intersect $\mathcal{B}_0$ with $\mu$ then we will get a finite set of points. 
Since the representative $\mu$ is generic, we conclude that the third derivative along 
$v$ will not vanish, i.e. the section $\Psi_{\PP A_3}$ will not vanish on those points. 
Hence, $\mathcal{B}_0\cap \mu$ does not contribute to the Euler class. \\ 
\hf\hf Next, let us 
%First of all we note that the whole of 
%$\mathcal{B}$ is not relevant while computing the contribution to the Euler class; only the points 
%at which the section vanishes is relevant. Hence, 
%As before, we only need to 
consider the components of 
$\mathcal{B}$ where one (or more) of the $q_i$ become equal to the last point $q_{\delta+1}$. 
Define $\mathcal{B}(q_{i_1}, \ldots q_{i_k}, l_{q_{\delta}})$ as before. 
%Let us define 
%\begin{align*}
%\B (q_1, l_{q_{\delta +1 }}) &:= \{ ([f], q_1, \ldots, q_{\delta}, l_{q_{\delta +1}}) \in \B: q_1 = q_{\delta +1} \}.
%\end{align*}
%The spaces $\B (q_2, q_{\delta +1 }), \ldots, \B (q_{\delta}, q_{\delta +1 })$ are defined 
%similarly. 
We show in \cite[Lemma 3.1]{BM13_2pt_published} 
%(Lemma 6.3 (3)) 
that 
%\begin{align*}
%\B 
%\end{align*}
%The claim is that 
\begin{align*}
\B (q_1, l_{q_{\delta +1 }}) & \approx \overline{\A_1^{\delta-1} \circ \PP \A}_4 \cup 
\overline{\A_1^{\delta-1} \circ \hat{D}}_5.
\end{align*}
We also show in \cite{BM13_2pt_published} (Corollary $6.13$) that 
the contribution to the Euler class from each of the points of 
$\overline{\A_1^{\delta-1} \circ \PP \A}_4 \cap \mu$
%$\B (q_1, l_{q_{\delta +1 }}) \cap \mu$ 
is $2$. 
Furthermore, 
%after we intersect with $\mu$, 
the section $\Psi_{\PP A_3}$ does not vanish on 
$\overline{\A_1^{\delta-1} \circ \hat{D}}_5 \cap \mu$, since $\mu$ is generic.  
Hence, the total contribution from all the components of type 
$\B (q_{i_1}, l_{q_{\delta +1 }})$ equals
\bgd
2\binom{\delta}{1}N(\A_1^{\delta-1}\PP\A_4,n_1,m_1,m_2,\theta).
\edd
%Next, let us define 
%\begin{align*}
%\B (q_1, q_2, l_{q_{\delta +1 }}) &:= 
%\{ ([f], q_1, \ldots, q_{\delta}, l_{q_{\delta +1}}) \in \B: q_1=q_2 = l_{q_{\delta +1}} \}.
%\end{align*}
Next, we analyze what happens when 
two nodes and one cusp come together. 
We show in \cite{BM_closure_of_seven_points} that 
%we claim that 
\begin{align*}
\B (q_1, q_2, l_{q_{\delta +1 }}) & \approx \overline{\A_1^{\delta-2}\circ \PP \D}_5.  
\end{align*}
Geometrically this is a believable statement due to the following picture: 
\newpage 
\begin{figure}[h!]
\vspace*{0.2cm}
\begin{center}\includegraphics[scale = 0.5]{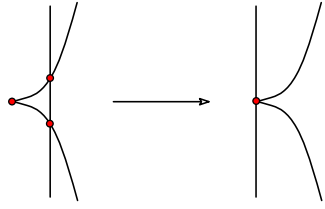}\vspace*{-0.2cm}\end{center}
\caption{Two nodes and one cusp collapsing to a $D_5$-singularity}
\end{figure} 
%\subsection{TT}
%\[ \textnormal{PICTURE}. \]
We also show that the contribution from each of the points of
$\overline{\A_1^{\delta-2}\circ \PP \D}_5 \cap \mu$
%$\B (q_1, q_2, l_{q_{\delta +1 }}) \cap \mu$ 
is $2$.
Hence the total contribution from all the components of type $\B (q_{i_1}, q_{i_2}, l_{q_{\delta +1 }})$ equals
\begin{align*}
2 \binom{\delta}{2} \N(\A_1^{\delta-2} \PP\D_5, n_1, m_1, m_2,\theta).   
\end{align*}
Finally, we claim that 
\begin{align*}
\B(q_1, q_2, \ldots, q_r, l_{q_{\delta+1}}) \cap \mu & = \varnothing \qquad \forall ~ 3 \leq r \leq 5.  
\end{align*}
In other words, a one dimensional family of curves with $r$-nodes and one cusp can not come together 
if $3\leq r \leq 5$.  
Hence the remaining stratum of $\B$ does not contribute to the Euler class, giving us Theorem \ref{a1_delta_pa3}.\qed

\subsection{Proof of Theorem \ref{a1_delta_pa4}: computing $\N(\A_1^{\delta} \PP \A_4, n_1, m_1, m_2, \theta)$}
\hf \hf Recall that 
in section \ref{notation}
we have defined the space  
\begin{align*}
\A^{\delta}_1 \circ \overline{\PP \A}_3 := \{ ([f], q_1, \ldots, q_{\delta}, 
l_{q_{\delta+1}}) \in \DD \times (X)^{\delta} \times \mathbb{P} T \X: 
&\textnormal{$f$ has a singularity of type $\A_1$ at $q_1, \ldots, q_{\delta}$}, \\ 
                   & ([f], l_{q_{\delta+1}}) \in \overline{\PP \A}_3, ~~\textnormal{$q_1, \ldots, q_{\delta+1}$ all distinct}\}. 
\end{align*}
Let $\mu$ be a generic pseudocycle representing the homology class 
Poincar\'{e} dual 
to 
\[c_1^{n_1} x_1^{m_1} x_2^{m_2} \lambda^{\theta}y^{\delta_L - (n_1 + m_1 + 2 m_2 +\theta+\delta+4)}.\]
Let $v \in \hat{\gamma}$ and $w \in \pi^*TX/\hat{\gamma}$ be two fixed non zero vectors. 
Let us introduce the following abbreviation:  
\begin{align*}
%\label{abbreviation}
f_{ij} & := \nabla^{i+j} f|_{q}
(\underbrace{v,\cdots v}_{\textnormal{$i$ times}}, \underbrace{w,\cdots w}_{\textnormal{$j$ times}}).
\end{align*}
We now define a section of the following bundle 
\begin{align*}
\psi_{\PP \A_4}:  \A^{\delta}_1 \circ \overline{\PP \A}_3 
\lra \mathbb{L}_{\PP \A_4} & := 
\gamma_{\D}^{\ast 2}\otimes\hat{\gamma}^{\ast 4}\otimes(\pi^* TX/\hat{\gamma})^{\ast 2}\otimes L^{\ast 2}, \\
\{\psi_{\PP \A_4}([f], q_1,\ldots,q_\delta,l_{q_{\delta+1}})\}(f^{\otimes 2}
\otimes v^{\otimes 4}\otimes w^{\otimes 2}) &:= f_{02}A_4^f, \qquad \textnormal{where} \qquad 
\A^{f}_4 := f_{40}-\frac{3 f_{21}^2}{f_{02}}.\nonumber
%\psi_{\A_1} : \psi_{\A_0}^{-1}(0) \lra \mathcal{V}_{\A_1} &:= \gamma_{\D}^*\otimes T^*\X \otimes L, 
%\qquad \textnormal{given by} \qquad 
%\{\psi_{\A_1}([f], q)\}(f) := \nabla f|_q. \label{psi_a0_a1_section_defn}
\end{align*}
%where 
%\begin{align*}
%\A^{f}_4 := f_{40}-\frac{3 f_{21}^2}{f_{02}}. 
%\end{align*}
%In \cite{RM_Hypersurfaces} we show that 
%Let us define 
%\begin{align*}
%M &:= \{ ([f], q_1, \ldots, q_{\delta}, l_{q_{\delta+1}}) \in A_1^{\delta} \circ \overline{\PP A}_3: 
%\nabla^2 f|_{q_{\delta+1}}(w,w) \neq 0\}.  
%\end{align*}
%\subsection{9AA}
By 
%Lemma $7.1$ statement $9$ of the arXiv version 
\cite[Lemma 7.1, statement 9]{BM13}, we have that 
\begin{align}
\overline{\PP A}_3 & = \PP A_3 \cup \overline{\PP A}_4 \cup \overline{\PP D}_4.  \label{pa3_closure}
\end{align}
%We prove this  in the Proof Theorem $4.20$ in \cite{BM13_one_singular_point_published}.
Hence, let us define 
\begin{align*}
\mathcal{B} &:= \overline{A_1^{\delta} \circ \overline{\PP A}}_3 - \A^{\delta}_1 \circ (\PP\A_3\cup \overline{\PP A}_4). 
\end{align*}
If $L$ is sufficiently $(2\delta+4)$-ample, then the section 
$\Psi_{\PP \A_4}$ vanishes on the points of  
$\A_1^{\delta}\circ \PP A_4$ transversally 
(all the points of $\A_1^{\delta}\circ \PP A_4$ are smooth points of the 
variety $\A_1^{\delta}\circ \overline{\PP A}_3$).
%***********
%If $L$ is sufficiently $(2\delta+4)$-ample, then the section $\Psi_{\PP A_4}$ restricted to 
%$M$ is transverse to the zero set. 
%Hence, the number of zeros of the section $\Psi_{\PP A_4}$ restricted to $M \cap \mu$ is 
%the number $\N(\A_1^{\delta}\PP \A_4, n_1, m_1, m_2, \theta)$.
Next, let us define 
\begin{align*}
\mathcal{B} &:= \overline{A_1^{\delta} \circ \PP A}_3 - \A^{\delta}_1 \circ \PP\A_3\cup \overline{\PP A}_4. 
\end{align*}
Hence,  
\begin{align*}
\lan e(\mathbb{L}_{\PP \A_4}), 
~~[\overline{\A^{\delta}_1 \circ \overline{\PP\A}}_3] \cap [\mu] \ran & = \N(\A_1^{\delta}\PP \A_4, n_1, m_1, m_2, \theta)  
+ \mathcal{C}_{\mathcal{B}\cap \mu}.
%\label{na1_pseudocycle}
\end{align*}
We now give an explicit description of $\mathcal{B}$. 
As before, let us first define $\mathcal{B}_0$ as  
\begin{align*}
\mathcal{B}_0 &:= \{ ([f], q_1, \ldots q_{\delta}, l_{q_{\delta+1}}) \in \mathcal{B}: q_1, q_2 \ldots q_{\delta+1} 
~~\textnormal{are all distinct}\}. 
\end{align*}
Hence, $\mathcal{B}_0$ is that component of the boundary, where all the points are still distinct.   
By \eqref{pa3_closure}, we conclude that 
\begin{align*}
\mathcal{B}_0 &= \overline{A^{\delta}_1} \circ \overline{\PP D}_4. 
\end{align*}
If we intersect $\mathcal{B}_0$ with $\mu$ then we will get a finite set of points. 
Since the representative $\mu$ is generic, we conclude that the directional derivative 
$f_{21}$ will not vanish on those points. Since $f_{02} =0$ on $\mathcal{B}_0$ 
we conclude that 
\[ f_{02}A^f_4 = f_{02}f_{40} - 3f_{21}^2 \neq 0 \] 
if $f_{21} \neq 0$.
%, i.e. $f_{02}A^f_4$ does not vanish on 
%$\mathcal{B}_0 \cap \mu$. 
Hence the section $\Psi_{\PP A_4}$ will not vanish on 
$\mathcal{B}_0\cap \mu$. \\
%those 
%points. \\ 
%the third derivative along 
%$v$ will not vanish, i.e. the section $\psi_{\PP A_3}$ will not vanish on those points. 
%Hence, $\mathcal{B}_0\cap \mu$ does not contribute to the Euler class. \\ 
\hf\hf Next, let us 
%First of all we note that the whole of 
%$\mathcal{B}$ is not relevant while computing the contribution to the Euler class; only the points 
%at which the section vanishes is relevant. Hence, 
%As before, we only need to 
consider the components of 
$\mathcal{B}$ where one (or more) of the $q_i$ become equal to the last point $q_{\delta+1}$. 
Define $\mathcal{B}(q_{i_1}, \ldots q_{i_k}, l_{q_{\delta}})$ as before. 
%Let us define 
%\begin{align*}
%\B (q_1, l_{q_{\delta +1 }}) &:= \{ ([f], q_1, \ldots, q_{\delta}, l_{q_{\delta +1}}) \in \B: q_1 = q_{\delta +1} \}.
%\end{align*}
%The spaces $\B (q_2, q_{\delta +1 }), \ldots, \B (q_{\delta}, q_{\delta +1 })$ are defined 
%similarly. 
%*************
%Next, let us define 
%\begin{align*}
%\mathcal{B} &:= \overline{\A^{\delta}_1 \circ \PP \A_3}- \A^{\delta}_1 \circ \PP\A_3. 
%\end{align*}
%Hence 
%\begin{align}
%\lan e(\mathbb{L}_{\PP \A_4}), 
%~~[\overline{\A^{\delta}_1 \circ \PP\A_3}] \cap [\mu] \ran & = \N(\A_1^{\delta}\PP \A_4, n_1, m_1, m_2, \theta)  
%+ \mathcal{C}_{\mathcal{B}\cap \mu}.
%\label{na1_pseudocycle}
%\end{align}
%We now give an explicit description of $\mathcal{B}$. First of all we note that the whole of 
%$\mathcal{B}$ is not relevant while computing the contribution to the Euler class; only the points 
%at which the section vanishes is relevant. 
%As before, we only need to consider the component of 
%$\mathcal{B}$ where one (or more) of the $q_i$ become equal to the last point $q_{\delta+1}$. 
%Define $\mathcal{B}(q_{i_1}, \ldots q_{i_k}, l_{q_{\delta}})$ as before.
%Let us define 
%\begin{align*}
%\B (q_1, l_{q_{\delta +1 }}) &:= \{ ([f], q_1, \ldots, q_{\delta}, l_{q_{\delta +1}}) \in \B: q_1 = q_{\delta +1} \}
%\end{align*}
%The spaces $\B (q_2, q_{\delta +1 }), \ldots, \B (q_{\delta}, q_{\delta +1 })$ are defined 
%similarly. 
%The claim is 
We show in \cite[Lemma 6.3(4)]{BM13_2pt_published}, 
%(Lemma 6.3 (4)) 
that 
%that 
\begin{align*}
\B (q_1, l_{q_{\delta +1 }}) & \approx \overline{\A_1^{\delta-1} \circ \PP \A}_5 \cup 
\overline{\A_1^{\delta-1} \circ \PP \D^{\vee}_5},
\end{align*}
where  we define $\PP D_k^{\vee}$ to be the space 
\begin{align}
\PP D_k^{\vee} &:= \{ ([f], l_{q}) \in \DD \times \mathbb{P} TX: ([f], q) \in \DD_k, ~~ f_{30}=0  
~~ \textnormal{and} ~~ f_{21} \neq 0\}. \label{pdk_dual} 
\end{align}
The section $\Psi_{\PP A_4}$ will not vanish on $\PP D_5^{\vee}$, since 
by definition on $\PP D_k^{\vee}$ we have $f_{02} =0$ but 
$f_{21} \neq 0$. Hence, the section will not vanish on 
$\overline{\A_1^{\delta-1} \circ \PP \D^{\vee}_5} \cap \mu$. 
We also show in \cite[Corollary 6.13]{BM13_2pt_published},
%(Corollary $6.13$) 
that the contribution 
of the section $\Psi_{\PP A_4}$ to the Euler class from 
%The contribution from 
each of the 
points of $\overline{\A_1^{\delta-1} \circ \PP \A}_5 \cap \mu$ is $2$. 
Hence the total contribution 
from 
all the components of type $\B (q_{i_1}, l_{q_{\delta +1 }})$ equals
\bgd
2\binom{\delta}{1}N(\A_1^{\delta-1}\PP\A_5,n_1,m_1,m_2,\theta).
\edd
Next, 
%let us consider the boundary when two nodes and one tacnode ****
we analyze what happens when 
two nodes and one tacnode come together. 
We show in \cite{BM_closure_of_seven_points} that
%\section{yy}
%define 
%We show in \cite{BM_closure_of_seven_points} that 
%we claim that 
\begin{align*}
\B (q_1, q_2, l_{q_{\delta +1 }}) \cap \mu  & \approx \overline{\A_1^{\delta-2}\circ \PP \D}_6 \cap \mu.  
\end{align*}
Geometrically this is a believable statement due to the following picture:
\newpage 
\begin{figure}[h!]
\label{figure_two_nodes_one_tacnode_to_d6}
\vspace*{0.2cm}
\begin{center}\includegraphics[scale = 0.5]{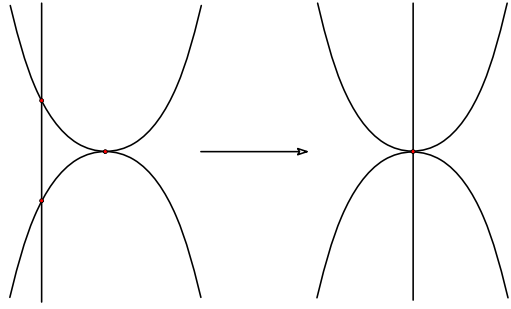}\vspace*{-0.2cm}\end{center}
\caption{Two nodes and one tacnode collapsing to a $D_6$-singularity}
\end{figure} 
We also show that the contribution from each of the points of 
$\overline{\A_1^{\delta-2}\circ \PP \D}_6 \cap \mu$
%$\B (q_1, q_2, l_{q_{\delta +1 }}) \cap \mu$ 
is $4$.
Hence the total contribution from all the components of type $\B (q_{i_1}, q_{i_2}, l_{q_{\delta +1 }})$ equals
\begin{align*}
4 \binom{\delta}{2} \N(\A_1^{\delta-2} \PP\D_6, n_1, m_1, m_2,\theta).   
\end{align*}
Finally, we show in \cite{BM_closure_of_seven_points} that 
\begin{align*}
\B(q_1, q_2, \ldots, q_r, l_{q_{\delta+1}}) \cap \mu & = \varnothing \qquad \forall ~ 3 \leq r \leq 4.  
\end{align*}
In other words, a one dimensional family of curves with $r$-nodes and one tacnode can not come together 
if $3\leq r \leq 4$.  
Hence the remaining stratum of $\B$ does not contribute to the Euler class, giving us Theorem \ref{a1_delta_pa4}.\qed 
%********
%\begin{align*}
%\B (q_1, q_2, l_{q_{\delta +1 }}) &:= 
%\{ ([f], q_1, \ldots, q_{\delta}, l_{q_{\delta +1}}) \in \B: q_1=q_2 = l_{q_{\delta +1}}\}.
%\end{align*}
%The claim is that 
%\begin{align*}
%\B (q_1, q_2, l_{q_{\delta +1 }}) & \approx \overline{\A_1^{\delta-2}\circ \PP \D_6}.  
%\end{align*}
%The contribution from each of the points of $\B (q_1, q_2, l_{q_{\delta +1 }}) \cap \mu$ is $4$.
%Hence the total contribution from all the components of type $\B (q_{i_1}, q_{i_2}, l_{q_{\delta +1 }})$ equals
%\begin{align*}
%4 \binom{\delta}{2} \N(\A_1^{\delta-2} \PP\D_6, n_1, m_1, m_2,\theta).   
%\end{align*}
%Finally, we claim that 
%\begin{align*}
%\B(q_1, q_2, \ldots, q_r, l_{q_{\delta+1}}) \cap \mu & = \varnothing \qquad \forall ~ 3 \leq r \leq 7.  
%\end{align*}
%Hence the remaining stratum of $\B$ does not contribute, since its intersection with $\mu$ 
%is empty. That gives us Theorem \ref{a1_delta_pa4}. \qed 

\subsection{Proof of Theorem \ref{a1_delta_pd4}: computing $\N(\A_1^{\delta} \PP \D_4, n_1, m_1, m_2, \theta)$}
\label{proof_of_a1_delta_pd4}
\hf \hf Let us first do the case $\theta=0$. 
%We continue with the setup of the proof of 
%Theorem \ref{a1_delta_pa4}; the spaces $M$ and $\mathcal{B}$ (with different decorations) 
%and the pseudocycle $\mu$ mean the same thing in this proof. 
%Recall that we have defined the space  
%\begin{align*}
%\A^{\delta}_1 \circ \PP \A_3 := \{ ([f], q_1, \ldots, q_{\delta}, 
%l_{q_{\delta+1}}) \in \DD \times (X)^{\delta} \times \mathbb{P} T \X: 
%&\textnormal{$f$ has a singularity of type $\A_1$ at $q_1, \ldots, q_{\delta}$}, \\ 
%                   & ([f], l_{q_{\delta+1}}) \in \PP \A_3, ~~\textnormal{$q_1, \ldots, q_{\delta+1}$ all distinct}\}. 
%\end{align*}
%Let $\mu$ be a generic pseudocycle representing the homology class 
%Poincar\'{e} dual 
%to 
%\[c_1^{n_1} x_1^{m_1} x_2^{m_2} \lambda^{\theta}y^{\delta_L - (n_1 + m_1 + 2 m_2 +\theta+\delta+4)}.\]
Define a section of the following bundle 
\begin{align}
\Psi_{\PP \D_4}:  \A^{\delta}_1 \circ \overline{\PP \A}_3 
\lra \mathbb{L}_{\PP \D_4} & := \gamma_{\D}^{\ast}\otimes(\pi^*TX/\hat{\gamma})^{\ast 2}\otimes L, \qquad 
\textnormal{given by} \nonumber \\
\{\Psi_{\PP \D_4}([f], q_1,\ldots,q_\delta,l_{q_{\delta+1}})\}(f\otimes w^{\otimes 2}) &:= \nabla^2f|_q(w,w). \nonumber
%\psi_{\A_1} : \psi_{\A_0}^{-1}(0) \lra \mathcal{V}_{\A_1} &:= \gamma_{\D}^*\otimes T^*\X \otimes L, 
%\qquad \textnormal{given by} \qquad 
%\{\psi_{\A_1}([f], q)\}(f) := \nabla f|_q. \label{psi_a0_a1_section_defn}
\end{align}
%In \cite{RM_Hypersurfaces} we show that 
%As stated earlier, 
By 
%Lemma $7.1$ statement $9$ of the arXiv version 
\cite[Lemma 7.1, statement 9]{BM13}, we conclude that 
%As stated earlier ()
\begin{align}
\overline{\PP A}_3 & = \PP A_3 \cup \overline{\PP A}_4 \cup \overline{\PP D}_4.  \label{pa3_closure_again}
\end{align}
Hence, let us define 
\begin{align*}
\mathcal{B} &:= \overline{A_1^{\delta} \circ \overline{\PP A}}_3 - \A^{\delta}_1 \circ (\PP\A_3\cup \overline{\PP D}_4). 
\end{align*}
If $L$ is sufficiently $(2\delta+4)$-ample, then the section 
$\Psi_{\PP \D_4}$ vanishes on the points of  
$\A_1^{\delta}\circ \PP D_4$ transversally 
(all the points of $\A_1^{\delta}\circ \PP D_4$ are smooth points of the 
variety $\A_1^{\delta}\circ \overline{\PP A}_3$).
%***********
%If $L$ is sufficiently $(2\delta+4)$-ample, then the section $\Psi_{\PP A_4}$ restricted to 
%$M$ is transverse to the zero set. 
%Hence, the number of zeros of the section $\Psi_{\PP A_4}$ restricted to $M \cap \mu$ is 
%the number $\N(\A_1^{\delta}\PP \A_4, n_1, m_1, m_2, \theta)$.
%Next, let us define 
%\begin{align*}
%\mathcal{B} &:= \overline{A_1^{\delta} \circ \overline{\PP A}}_3 - \A^{\delta}_1 \circ \PP\A_3\cup \overline{\PP A}_4. 
%\end{align*}
%Hence,  
%*********************
%Hence, the number of zeros of the section $\Psi_{\PP \D_4}$, restricted to 
%$\A^{\delta}_1 \circ \overline{\PP \A}_3 \cap \mu$ is the number 
%$\N(\A_1^{\delta}\PP \D_4, n_1, m_1, m_2, \theta)$ 
%(since the section $\Psi_{\PP \D_4}$ does not vanish on $\PP A_4$).
%If $L$ is sufficiently $(2\delta+4)$-ample, then the section $\Psi_{\PP \D_4}$ 
%is transverse to the zero set. 
%Next, let us define 
%\begin{align*}
%\mathcal{B} &:= \overline{\A^{\delta}_1 \circ \overline{\PP \A}}_3- \A^{\delta}_1 \circ \PP\A_3. 
%\end{align*}
Hence, 
\begin{align*}
\lan e(\mathbb{L}_{\PP \D_4}), 
~~[\overline{\A^{\delta}_1 \circ \PP\A}_3] \cap [\mu] \ran & = \N(\A_1^{\delta}\PP \D_4, n_1, m_1, m_2, \theta)  
+ \mathcal{C}_{\mathcal{B}\cap \mu}.
%\label{na1_pseudocycle}
\end{align*}
We now give an explicit description of $\mathcal{B}$.
Let us define $\mathcal{B}_0$ as before; it denotes the component of the 
boundary where all the points are distinct. The section 
$\Psi_{\PP D_4}$ does not vanish on $\overline{A_1^{\delta}} \circ \PP A_4$; hence 
\eqref{pa3_closure_again}, the section does not vanish on $\mathcal{B}_0$. \\ 
\hf \hf Next, let us analyze the component of 
$\mathcal{B}$ where one (or more) of the $q_i$ become equal to the last point $q_{\delta+1}$. 
Define $\mathcal{B}(q_{i_1}, \ldots q_{i_k}, l_{q_{\delta}})$ as before.
%First of all we note that the whole of 
%$\mathcal{B}$ is not relevant while computing the contribution to the Euler class; only the points 
%at which the section vanishes is relevant. Hence, we only need to consider the component of 
%$\mathcal{B}$ where one (or more) of the $q_i$ become equal to the last point $q_{\delta+1}$. 
%Let us define 
%\begin{align*}
%\B (q_1, l_{q_{\delta +1 }}) &:= \{ ([f], q_1, \ldots, q_{\delta}, l_{q_{\delta +1}}) \in \B: q_1 = q_{\delta +1} \}
%\end{align*}
%The spaces $\B (q_2, q_{\delta +1 }), \ldots, \B (q_{\delta}, q_{\delta +1 })$ are defined 
%similarly. 
We show in \cite[Lemma 6.3(4)]{BM13_2pt_published} 
%(Lemma 6.3 (4)) 
%that 
that 
\begin{align*}
\B (q_1, l_{q_{\delta +1 }}) & \approx \overline{\A_1^{\delta-1} \circ \PP \A}_5 \cup 
\overline{\A_1^{\delta-1} \circ \PP \D^{\vee}_5},
\end{align*}
where $\PP D_5^{\vee}$ is as defined in \eqref{pdk_dual}. 
We also show in \cite[Corollary 6.17]{BM13_2pt_published} 
%(Corollary $6.17$), 
that 
the section 
$\Psi_{\PP D_4}$ vanishes on $\A_1^{\delta-1} \circ \PP D_5^{\vee} \cap \mu$ with a 
multiplicity of $2$. Furthermore, we notice that the section does not vanish 
on $\A_1^{\delta-1}\circ \PP A_5 \cap \mu$. 
%The claim is that 
%\begin{align*}
%\B (q_1, l_{q_{\delta +1 }}) & \approx \overline{\A_1^{\delta-1} \circ \PP \A_5} \cup 
%\overline{\A_1^{\delta-1} \circ \PP \D_5^{\vee}}.
%\end{align*}
%\begin{align*}
%\B (q_1, l_{q_{\delta +1 }})\cap \{\Psi_{\PP \D_4}=0\} & \approx \overline{\A_1^{\delta-1} \circ \PP \D_5}.
%\end{align*}
%Not that we are \textit{not} claiming that 
%\begin{align*}
%\B (q_1, l_{q_{\delta +1 }}) & \approx \overline{\A_1^{\delta-1} \circ \PP \D_5}.
%\end{align*}
%There will 
%The contribution from each of the points of $\B (q_1, l_{q_{\delta +1 }}) \cap \mu$ is $2$. 
Hence the total contribution from all the components of type $\B (q_{i_1}, l_{q_{\delta +1 }})$ equals
\bgd
2\binom{\delta}{1}N(\A_1^{\delta-1}\PP\D_5^{\vee},n_1,m_1,m_2,0).
\edd
Next, we note that 
the projection map 
%$\pi: \A_1^{\delta}\PP \D_5 \lra \D_5 $ and 
$\pi: \A_1^{\delta} \circ \PP \D^{\vee}_5 \lra \A_1^{\delta} \circ  \D_5$ 
is one to one. 
Hence, 
\[ N(\A_1^{\delta-1}\PP\D_5^{\vee},n_1,m_1,m_2,0) = N(\A_1^{\delta-1} \D_5, n_1, m_1, m_2). \] 
%= N(\A_1^{\delta-1}\PP\D_5,n_1,m_1,m_2,0) = 
Hence the total contribution from all the components of type $\B (q_{i_1}, l_{q_{\delta +1 }})$ equals
\bgd
2\binom{\delta}{1}N(\A_1^{\delta-1} \D_5,n_1,m_1,m_2).
\edd
%Next, let us define 
%\begin{align*}
%\B (q_1, q_2, l_{q_{\delta +1 }}) &:= 
%\{ ([f], q_1, \ldots, q_{\delta}, l_{q_{\delta +1}}) \in \B: q_1=q_2 = l_{q_{\delta +1}}. \}
%\end{align*}
Next, let us describe the boundary when three of the points come together. 
As mentioned earlier in the Proof of 
Theorem \ref{a1_delta_pa4}, 
we show in \cite{BM_closure_of_seven_points} that 
%The claim is that 
\begin{align*}
\B (q_1, q_2, l_{q_{\delta +1 }}) & \approx \overline{\A_1^{\delta-2}\circ \PP \D}_6.  
\end{align*}
%The intuitive justification for this claim is the picture drawn in 
%\ref{figure_two_nodes_one_tacnode_to_d6}. 
We also show that the contribution from each of the points of 
$\overline{\A_1^{\delta-2}\circ \PP \D}_6 \cap \mu$ is $2$.
Hence the total contribution from all the components of type $\B (q_{i_1}, q_{i_2}, l_{q_{\delta +1 }})$ equals
\begin{align*}
2 \binom{\delta}{2} \N(\A_1^{\delta-2} \PP\D_6, n_1, m_1, m_2,\theta).   
\end{align*}
Finally, we show in 
\cite{BM_closure_of_seven_points} that 
\begin{align*}
\B(q_1, q_2, \ldots, q_r, l_{q_{\delta+1}}) \cap \mu & = \varnothing \qquad \forall ~ 3 \leq r \leq 4.  
\end{align*}
In other words, if we have a one dimensional family of curves with $r$-nodes and one tacnode, then 
all the points can not come together, if $r$ is $3$ or $4$. 
Hence the remaining stratum of $\B$ does not contribute to the Euler class, 
giving us Theorem \ref{a1_delta_pd4} for $\theta=0$. \\
\hf\hf Next, we consider the case $\theta=1$. 
%Recall that in section \ref{notation} we 
%have defined 
Let us consider the following space  
in section \ref{notation}, namely 
\begin{align*}
\A^{\delta}_1 \circ \overline{\hat{D}}_4 := \{ ([f], q_1, \ldots, q_{\delta}, 
l_{q_{\delta+1}}) \in \DD \times (X)^{\delta} \times \mathbb{P} T \X: 
&\textnormal{$f$ has a singularity of type $\A_1$ at $q_1, \ldots, q_{\delta}$}, \\ 
                   & ([f], l_{q_{\delta+1}}) \in \overline{\hat{D}}_4, ~~\textnormal{$q_1, \ldots, q_{\delta+1}$ all distinct}\}. 
\end{align*}
Let $\mu$ be a generic pseudocycle representing the homology class 
Poincar\'{e} dual 
to 
\[c_1^{n_1} x_1^{m_1} x_2^{m_2} \lambda^{\theta}y^{\delta_L - (n_1 + m_1 + 2 m_2 +\theta+\delta+4)}.\]
We again consider the  section of the following bundle 
\begin{align}
\Psi_{\PP \A_3}:  \A^{\delta}_1 \circ \overline{\hat{\D}}_4 
\lra \mathbb{L}_{\PP \A_3} & := \gamma_{\D}^\ast \otimes \hat{\gamma}^{\ast 3} \otimes L, \nonumber \\
\{\Psi_{\PP \A_3}([f], q_1,\ldots,q_\delta,l_{q_{\delta+1}})\}(f\otimes v^{\otimes 3}) &:= \nabla^3f|_q(v,v,v).  \nonumber
%\psi_{\A_1} : \psi_{\A_0}^{-1}(0) \lra \mathcal{V}_{\A_1} &:= \gamma_{\D}^*\otimes T^*\X \otimes L, 
%\qquad \textnormal{given by} \qquad 
%\{\psi_{\A_1}([f], q)\}(f) := \nabla f|_q. \label{psi_a0_a1_section_defn}
\end{align}
%Note that 
%\begin{align}
%\overline{\hat{D}}_4 & = \hat{D}_4 \cup \overline{\hat{D}}_5. 
%\end{align}
If $L$ is sufficiently $(2\delta+4)$-ample, then this section is transverse to the zero set. 
Let, 
\begin{align*}
\mathcal{B} &:= \overline{\A^{\delta}_1 \circ \overline{\hat{\D}}}_4- \A^{\delta}_1 \circ \hat{\D}_4. 
\end{align*}
Hence 
\begin{align*}
\lan e(\mathbb{L}_{\PP \A_3}), 
~~[\overline{\A^{\delta}_1 \circ \hat{\D}_4}] \cap [\mu] \ran & = \N(\A_1^{\delta}\PP \D_4, n_1, m_1, m_2, \theta)  
+ \mathcal{C}_{\mathcal{B}\cap \mu}.
%\label{na1_pseudocycle}
\end{align*}
Let us denote $\B_0$ to be the component of the boundary where all the points are distinct. 
%from the last point. 
The section $\Psi_{\PP \A_3}$ will not vanish on $\B_0 \cap \mu$, 
since the last marked point will comprise of a generic $D_5$-singularity. \\ 
\hf \hf Next, let us consider the component when one or more points come together.  
Define $\B(q_1, \ldots, q_r, l_{q_{\delta+1}})$ as before. 
We show in \cite{BM_closure_of_seven_points} that 
\begin{align*}
\B(q_1, l_{q_{\delta+1}}) & = \overline{\hat{D}}_6.  
\end{align*}
Again, when we intersect $\overline{\hat{D}}_6$, the section 
$\Psi_{\PP A_3}$ will not vanish. Hence, it does not contribute to the Euler class. 
Next, we show in \cite{BM_closure_of_seven_points} that
%However, *.
\begin{align*}
\B(q_1, q_2, l_{q_{\delta+1}}) \cap \mu & = \varnothing.  
\end{align*}
Hence, $\B(q_1, q_2, l_{q_{\delta+1}})$ does not contribute to the Euler class. 
Finally, we show that when three nodes and one triple point come together, we get a 
triple point. In other words, we show in \cite{BM_closure_of_seven_points} that: 
\begin{align}
\B(q_1, q_2, q_3, l_{q_{4}}) & = \overline{\hat{X}}_9.  \label{three_nodes_d4_hat_x9}
\end{align}
Geometrically this is a believable statement due to the following picture:
\newpage 
\begin{figure}[h!]
\vspace*{0.2cm}
\begin{center}\includegraphics[scale = 0.5]{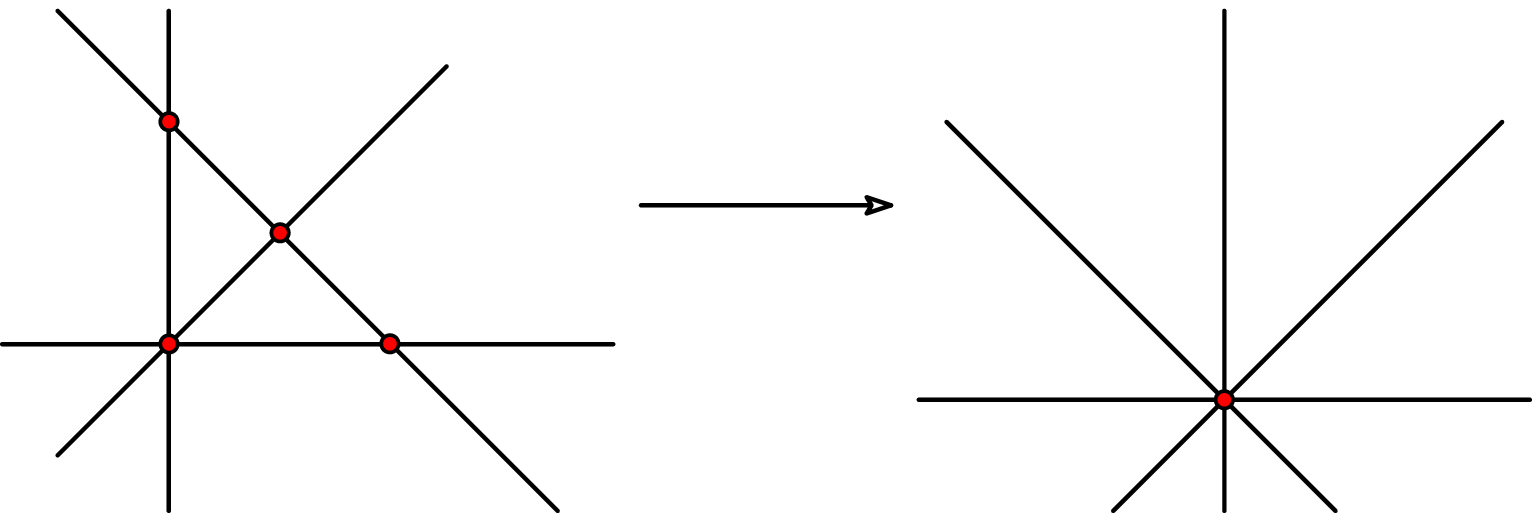}\vspace*{-0.2cm}\end{center}
\caption{Three nodes and one $D_4$-singularity collapsing to a quadruple point}
\end{figure} 
We also show in \cite{BM_closure_of_seven_points} that 
the contribution of the section from 
each of the points of $\overline{\hat{X}}_9 \cap \mu$ 
%$\B(q_1, q_2, q_3, l_{q_{4}}) \cap \mu$ 
is $24$. 
That gives us 
Theorem \ref{a1_delta_pd4} for $\theta=1$. \qed 

\subsection{Proof of Theorem \ref{a1_delta_pd5}: computing $\N(\A_1^{\delta} \PP \D_5, n_1, m_1, m_2, \theta)$}
\hf \hf Recall that in section \ref{notation}, we have defined the following space  
\begin{align*}
\A^{\delta}_1 \circ \overline{\PP \D}_4 := \{ ([f], q_1, \ldots, q_{\delta}, 
l_{q_{\delta+1}}) \in \DD \times (X)^{\delta} \times \mathbb{P} T \X: 
&\textnormal{$f$ has a singularity of type $\A_1$ at $q_1, \ldots, q_{\delta}$}, \\ 
                   & ([f], l_{q_{\delta+1}}) \in \overline{\PP \D}_4, 
                   ~~\textnormal{$q_1, \ldots, q_{\delta+1}$ all distinct}\}. 
\end{align*}
Let $\mu$ be a generic pseudocycle representing the homology class 
Poincar\'{e} dual 
to 
\[c_1^{n_1} x_1^{m_1} x_2^{m_2} \lambda^{\theta}y^{\delta_L - (n_1 + m_1 + 2 m_2 +\theta+\delta+5)}.\]
As before, we let $v \in \hat{\gamma}$ and $w \in \pi^*TX/\hat{\gamma}$ be two fixed non zero vectors 
and let $f_{ij}$ be defined as before. 
%\begin{align*}
%\label{abbreviation}
%f_{ij} & := \nabla^{i+j} f|_{q}
%(\underbrace{v,\cdots v}_{\textnormal{$i$ times}}, \underbrace{w,\cdots w}_{\textnormal{$j$ times}}).
%\end{align*}
We now define a section of the following bundle 
\begin{align}
\Psi_{\PP \D_5}:  \A^{\delta}_1 \circ \overline{\PP \D}_4 
\lra \mathbb{L}_{\PP \D_5} & := \gamma_{\D}^{\ast} \otimes \hat{\gamma}^{\ast 2} \otimes 
(\pi^*TX/\hat{\gamma})^{\ast} \otimes L^{\ast} , \qquad \textnormal{given by} \nonumber \\
\{\Psi_{\PP \D_5}([f], q_1,\ldots,q_\delta,l_{q_{\delta+1}})\}(f\otimes v^{\otimes 2}\otimes w) &:= 
f_{21}. \nonumber
%\psi_{\A_1} : \psi_{\A_0}^{-1}(0) \lra \mathcal{V}_{\A_1} &:= \gamma_{\D}^*\otimes T^*\X \otimes L, 
%\qquad \textnormal{given by} \qquad 
%\{\psi_{\A_1}([f], q)\}(f) := \nabla f|_q. \label{psi_a0_a1_section_defn}
\end{align}
%In \cite{RM_Hypersurfaces} we show that
In 
%Lemma 
%$7.1$ statement $4$ of the arXiv version 
\cite[Lemma 7.1, statement 4]{BM13}, we prove that 
%As stated earlier ()
\begin{align}
\overline{\PP D}_4 & = \PP D_4 \cup \overline{\PP D}_5 \cup \overline{\PP D^{\vee}_5}.  \label{pd4_closure_again}
\end{align}
Hence, let us define 
\begin{align*}
\mathcal{B} &:= \overline{A_1^{\delta} \circ \overline{\PP D}}_4 - \A^{\delta}_1 \circ (\PP\D_4\cup \overline{\PP D}_5). 
\end{align*}
If $L$ is sufficiently $(2\delta+4)$-ample, then the section 
$\Psi_{\PP \D_5}$ vanishes on the points of  
$\A_1^{\delta}\circ \PP D_5$ transversally 
(all the points of $\A_1^{\delta}\circ \PP D_5$ are smooth points of the 
variety $\A_1^{\delta}\circ \overline{\PP D}_4$).
%************************************* 
%If $L$ is sufficiently $(2\delta+5)$-ample, then this section is transverse to the zero set. 
%Next, let us define 
%\begin{align*}
%\mathcal{B} &:= \overline{\A^{\delta}_1 \circ \PP \D_4}- \A^{\delta}_1 \circ \PP\D_4. 
%\end{align*}
Hence, 
\begin{align*}
\lan e(\mathbb{L}_{\PP \D_5}), 
~~[\overline{\A^{\delta}_1 \circ \overline{\PP\D}}_4] \cap [\mu] \ran & = \N(\A_1^{\delta}\PP \D_5, n_1, m_1, m_2, \theta)  
+ \mathcal{C}_{\mathcal{B}\cap \mu}.
%\label{na1_pseudocycle}
\end{align*}
Let us define $\B_0$ as before, namely the component of the boundary where 
all the points are distinct. 
%none of the 
%$(q_1, \ldots, q_{\delta})$ points are equal to the last marked point $q_{\delta+1}$. 
The section $\Psi_{\PP D_5}$ does not vanish on $\PP D_5^{\vee}$ by definition 
(recall that by definition of the space $\PP D_5^{\vee}$, the directional derivative 
$f_{21} \neq 0$). 
Hence, the section $\Psi_{\PP D_5}$ does not vanish 
on $\overline{A^{\delta}_1} \circ \overline{\PP D^{\vee}_5} \cap \mu$ for a generic $\mu$. 
Therefore, $\B_0$ does not contribute to the Euler class. \\ 
\hf \hf Next, let us analyze the boundary when one 
or more points come together; define 
%We now give an explicit description of $\mathcal{B}$. 
%First of all we note that the whole of 
%$\mathcal{B}$ is not relevant while computing the contribution to the Euler class; only the points 
%at which the section vanishes is relevant. Hence, we only need to consider the component of 
%$\mathcal{B}$ where one (or more) of the $q_i$ become equal to the last point $q_{\delta+1}$. 
%Let us define 
%\begin{align*}
%\B (q_1, l_{q_{\delta +1 }}) &:= \{ ([f], q_1, \ldots, q_{\delta}, l_{q_{\delta +1}}) \in \B: q_1 = q_{\delta +1} \}
%\end{align*}
%The spaces $\B (q_2, q_{\delta +1 }), \ldots, \B (q_{\delta}, q_{\delta +1 })$ are defined 
%similarly. 
%As before, we only need to consider the component of 
%$\mathcal{B}$ where one (or more) of the $q_i$ become equal to the last point $q_{\delta+1}$. 
%Define 
$\mathcal{B}(q_{i_1}, \ldots q_{i_k}, l_{q_{\delta}})$ as before.
In \cite{BM_closure_of_seven_points}, we show that 
%The claim is that 
\begin{align*}
\B (q_1, l_{q_{\delta +1 }}) & \approx 
\overline{\A_1^{\delta-1} \circ \PP \D_6} \cup \overline{\A_1^{\delta-1} \circ \PP \D_6^{\vee}}. 
\end{align*}
%\begin{align*}
%\B (q_1, l_{q_{\delta +1 }})\cap\{\Psi_{\PP \D_5}=0\} & \approx \overline{\A_1^{\delta-1} \circ \PP \D_6}.
%\end{align*}
Again, by the definition of $\PP D_k^{\vee}$, the directional derivative $f_{21}$ does not 
vanish on $\PP D_6^{\vee}$; 
hence the section $\Psi_{\PP \D_5}$ 
does not vanish on $\overline{\A_1^{\delta-1} \circ \PP \D_6^{\vee}}\cap \mu$. 
%Hence if $\mu$ is a generic 
%representative of the Poincar\'{e} dual of
%\[c_1^{n_1} x_1^{m_1} x_2^{m_2} \lambda^{\theta}y^{\delta_L - (n_1 + m_1 + 2 m_2 +\theta+\delta+5)}\]
%then it will not intersect  $\ov{\A_1^{\delta-1} \circ \PP \D^{\vee}_6}-\A_1^{\delta-1} \circ \PP \D^{\vee}_6$ anywhere. 
We also show in \cite{BM_closure_of_seven_points} that 
the section vanishes on $\overline{\A_1^{\delta-1} \circ \PP \D_6}\cap \mu$ 
%and on $\overline{\A_1^{\delta-1} \circ \hat{X}_9}\cap \mu$
%. We claim that 
%it vanishes 
with a multiplicity of $2$. 
%and $12$ respectively. 
%Hence,  
%The contribution from each of the points of 
%$\B (q_1, l_{q_{\delta +1 }}) \cap \mu$ is $2$. 
Hence the total contribution from all the 
components of type $\B (q_{i_1}, l_{q_{\delta +1 }})$ equals
\bgd
2\binom{\delta}{1}N(\A_1^{\delta-1}\PP\D_6,n_1,m_1,m_2,\theta).
\edd
Next, let us consider the boundary when two points come together. 
We show in \cite{BM_closure_of_seven_points} that 
\begin{align*}
\B (q_1, q_2, l_{q_{\delta +1 }}) & \approx 
\overline{\A_1^{\delta-2} \circ \hat{X}_9}. 
\end{align*}
Furthermore, we show that the section vanishes on 
$\overline{\A_1^{\delta-2} \circ \hat{X}_9}$. Hence, the 
contribution from all points of type $\B (q_{i_1}, q_{i_2}, l_{q_{\delta +1 }})$ 
equals 
\bgd
12\binom{\delta}{1}N(\A_1^{\delta-2}\hat{X}_9, n_1,m_1,m_2,\theta).
\edd
Finally, we 
show in \cite{BM_closure_of_seven_points} that 
%if we have a one dimensional family of 
%curves with $3$-nodes and one singularity of type $\PP D_4$, then 
%we 
%all the singular 
%points can not come together. More precisely, we show that 
%claim that 
\begin{align}
\B(q_1, q_2, q_3, l_{q_{4}}) \cap & = \overline{\PP X}_9.  \label{three_nodes_and_one_pd4}
\end{align}
This is geometrically believable due to the following picture: 
\newpage
\begin{figure}[h!]
\vspace*{0.2cm}
\begin{center}\includegraphics[scale = 0.5]{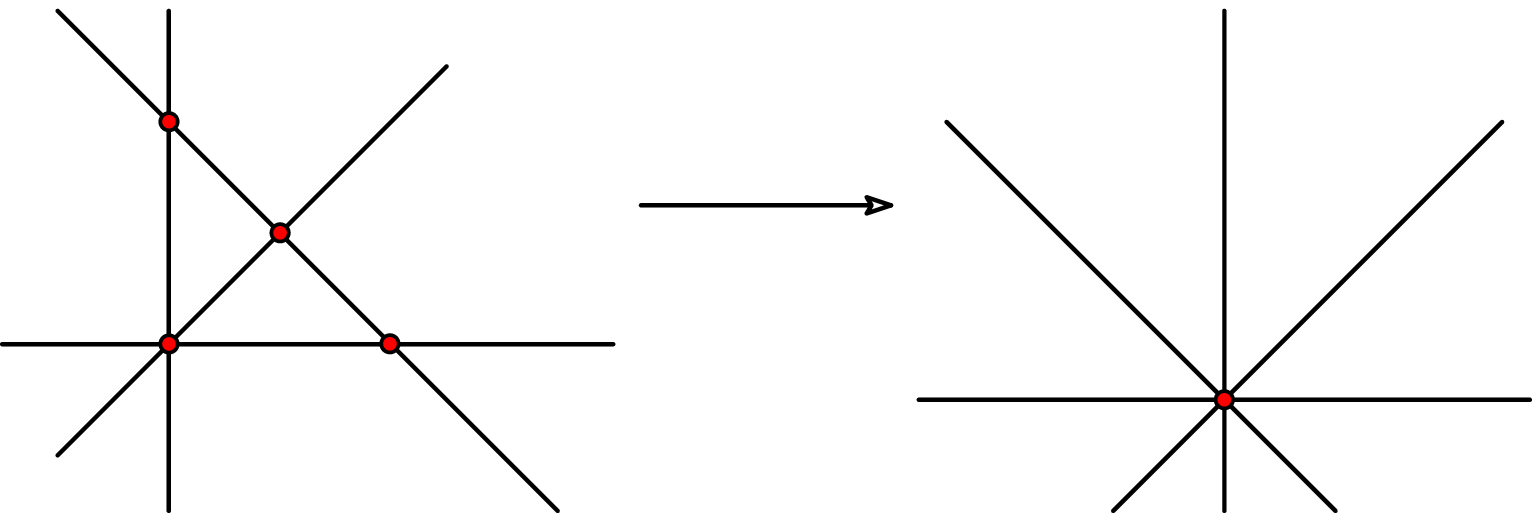}\vspace*{-0.2cm}\end{center}
\caption{Three nodes and one triple point collapsing to a quadruple point}
\end{figure} 
%\subsection{JJ}
%\[ \textnormal{PICTURE.} \] 
We also show in \cite{BM_closure_of_seven_points} that the contribution from 
each of the point of $\overline{\PP X}_9$ is $18$, thereby giving 
us Theorem \ref{a1_delta_pd5}. \qed  

\begin{rem}
We remind the reader again that there is a subtle difference between 
\eqref{three_nodes_d4_hat_x9} 
and 
\eqref{three_nodes_and_one_pd4}. Both the statements happen to be 
true, but they are not exactly the same statements. In the first statement  
we are looking at the space of curves with three nodes and one triple point; in the 
closure we get a quadruple point. In other words we were looking at the space 
$A_1^3 \circ \hat{D}_4$ and getting $\hat{X}_9$ in the closure.\\ 
\hf \hf In the second case, we are looking at the space of curves with three 
nodes and a triple point with a distinguished branch along which the third derivative 
vanishes. In other words, we are looking at the space $A_1\circ\PP D_4$ as opposed 
to $A_1^{3} \circ \hat{D}_4$. What we get in the closure is not any random quadruple point. 
We get a quadruple point with a distinguished branch such that the fourth derivative 
along that direction vanishes. In other words, we do not get $\hat{X}_9$ in the closure, but 
rather $\PP X_9$. 
\end{rem}

%Next, let us define 
%\begin{align*}
%\B (q_1, q_2, l_{q_{\delta +1 }}) &:= 
%\{ ([f], q_1, \ldots, q_{\delta}, l_{q_{\delta +1}}) \in \B: q_1=q_2 = l_{q_{\delta +1}}. \}
%\end{align*}
%The claim is that 
%\begin{align*}
%\B (q_1, q_2, l_{q_{\delta +1 }}) & \approx \overline{\A_1^{\delta-2}\circ \PP \D_6}.  
%\end{align*}
%The contribution from each of the points of $\B (q_1, q_2, l_{q_{\delta +1 }}) \cap \mu$ is $4$.
%Hence the total contribution from all the components of type $\B (q_{i_1}, q_{i_2}, l_{q_{\delta +1 }})$ equals
%\begin{align*}
%4 \binom{\delta}{2} \N(\A_1^{\delta-2} \PP\D_6, n_1, m_1, m_2,\theta).   
%\end{align*}
%Finally, we claim that 
%\begin{align*}
%\B(q_1, q_2, \ldots, q_r, l_{q_{\delta+1}}) \cap \mu & = \varnothing \qquad \forall ~ 3 \leq r \leq 7.  
%\end{align*}
%Hence the remaining stratum of $\B$ does not contribute, since its intersection with $\mu$ 
%is empty. \qed 
%{\bf The above presumably holds for $\theta=1$ as well. What we have is for $\theta=0$. }\qed

\subsection{Proof of Theorem \ref{a1_delta_pa5}: computing $\N(\A_1^{\delta} \PP \A_5, n_1, m_1, m_2, \theta)$}
\hf\hf Recall that in section \ref{notation} we have defined the space  
\begin{align*}
\A^{\delta}_1 \circ \overline{\PP \A}_4 := \{ ([f], q_1, \ldots, q_{\delta}, 
l_{q_{\delta+1}}) \in \DD \times (X)^{\delta} \times \mathbb{P} T \X: 
&\textnormal{$f$ has a singularity of type $\A_1$ at $q_1, \ldots, q_{\delta}$}, \\ 
                   & ([f], l_{q_{\delta+1}}) \in \overline{\PP \A}_4, ~~\textnormal{$q_1, \ldots, q_{\delta+1}$ all distinct}\}. 
\end{align*}
Let $\mu$ be a generic pseudocycle representing the homology class 
Poincar\'{e} dual to 
\[c_1^{n_1} x_1^{m_1} x_2^{m_2} \lambda^{\theta}y^{\delta_L - (n_1 + m_1 + 2 m_2 +\theta+\delta+5)}.\]
As before, we let $v \in \hat{\gamma}$ and $w \in \pi^*TX/\hat{\gamma}$ be two fixed non zero vectors 
and define $f_{ij}$ as before. 
We now define a section of the following bundle 
\begin{align*}
\Psi_{\PP \A_5}:  \A^{\delta}_1 \circ \overline{\PP \A}_4 
\lra \mathbb{L}_{\PP \A_5} & := \hat{\gamma}^{\ast 5}
\otimes\gamma_{\D}^{\ast 3}\otimes(\pi^*TX/\hat{\gamma})^{\ast 4} \otimes L^{\ast 3} \qquad \textnormal{given by} \\
\{\Psi_{\PP \A_5}([f], q_1,\ldots,q_\delta, l_{q_{\delta+1}})\}(f^{\otimes 2}
\otimes v^{\otimes 5}\otimes w^{\otimes 4}) &:= f_{02}^2 A_5^f, \nonumber 
%\qquad \textnormal{where} 
%\qquad A^f_5:=  f_{50} -\frac{10 f_{21} f_{31}}{f_{02}} + 
%\frac{15 f_{12} f_{21}^2}{f_{02}^2}.\nonumber
%\psi_{\A_1} : \psi_{\A_0}^{-1}(0) \lra \mathcal{V}_{\A_1} &:= \gamma_{\D}^*\otimes T^*\X \otimes L, 
%\qquad \textnormal{given by} \qquad 
%\{\psi_{\A_1}([f], q)\}(f) := \nabla f|_q. \label{psi_a0_a1_section_defn}
\end{align*}
where
\[ \A^{f}_5 := f_{50} -\frac{10 f_{21} f_{31}}{f_{02}} + 
\frac{15 f_{12} f_{21}^2}{f_{02}^2}.\]
%In \cite{RM_Hypersurfaces} we show that 
In 
%Lemma $7.1$ statement $*$ of the arXiv version 
\cite[Lemma 7.1, statement 10]{BM13}, we prove that 
%As stated earlier ()
\begin{align}
\overline{\PP A}_4 & = \PP A_4 \cup \overline{\PP A}_5 \cup \overline{\PP D}_5.  \label{pa4_closure}
\end{align}
Hence, let us define 
\begin{align*}
\mathcal{B} &:= \overline{A_1^{\delta} \circ \overline{\PP A}}_4 - \A^{\delta}_1 \circ (\PP\A_4\cup \overline{\PP A}_5). 
\end{align*}
If $L$ is sufficiently $(2\delta+5)$-ample, then the section 
$\Psi_{\PP \A_5}$ vanishes on the points of  
$\A_1^{\delta}\circ \PP A_5$ transversally 
(all the points of $\A_1^{\delta}\circ \PP A_5$ are smooth points of the 
variety $\A_1^{\delta}\circ \overline{\PP A}_4$).
%************************************* 
%If $L$ is sufficiently $(2\delta+5)$-ample, then this section is transverse to the zero set. 
%Next, let us define 
%\begin{align*}
%\mathcal{B} &:= \overline{\A^{\delta}_1 \circ \PP \D_4}- \A^{\delta}_1 \circ \PP\D_4. 
%\end{align*}
Hence, 
\begin{align*}
\lan e(\mathbb{L}_{\PP \A_5}), 
~~[\overline{\A^{\delta}_1 \circ \overline{\PP\A}}_4] \cap [\mu] \ran & = \N(\A_1^{\delta}\PP \A_5, n_1, m_1, m_2, \theta)  
+ \mathcal{C}_{\mathcal{B}\cap \mu}.
%\label{na1_pseudocycle}
\end{align*}
Let us define $\B_0$ as before, namely the component of the boundary where none of the 
$(q_1, \ldots, q_{\delta})$ points are equal to the last marked point $q_{\delta+1}$. 
We show in \cite[Corollary 7.4]{BM13}, 
%Lemma *, 
that the section $\Psi_{\PP A_5}$ 
vanishes on all points of $\overline{A_1^{\delta}}\circ \PP \D_5$ with a multiplicity 
of $2$. 
%The section $\Psi_{\PP D_5}$ does not vanish on $\PP D_5^{\vee}$ by definition 
%(recall that by definition of the space $\PP D_5^{\vee}$, the directional derivative 
%$f_{21} \neq 0$). 
%Hence, the section $\Psi_{\PP D_5}$ does not vanish 
%on $\overline{A^{\delta}_1} \circ \overline{\PP D^{\vee}_5} \cap \mu$ for a generic $\mu$. 
Hence, the contribution of $\B_0$ 
to the Euler class is 
\begin{align*}
2 \Num(\A_1^{\delta}\PP \D_5, n_1, m_1, m_2, \theta). 
\end{align*}
%does not contribute to the Euler class. \\ 
%*****************
%If $L$ is sufficiently $(2\delta+5)$-ample, then this section is transverse to the zero set. 
%Next, let us define 
%\begin{align*}
%\mathcal{B} &:= \overline{\A^{\delta}_1 \circ \PP \A}_4- \A^{\delta}_1 \circ \PP\A_4. 
%\end{align*}
%Hence 
%\begin{align}
%\lan e(\mathbb{L}_{\PP \A_5}), 
%~~[\overline{\A^{\delta}_1 \circ \PP\A}_4] \cap [\mu] \ran & = \N(\A_1^{\delta}\PP \A_5, n_1, m_1, m_2, \theta)  
%+ \mathcal{C}_{\mathcal{B}\cap \mu}.
%\label{na1_pseudocycle}
%\end{align}
%We now give an explicit description of $\mathcal{B}$.
\hf \hf Next, let us analyze the boundary when 
one (or more) of the $q_i$ become equal to the last point $q_{\delta+1}$.
Define $\mathcal{B}(q_{i_1}, \ldots q_{i_k}, l_{q_{\delta}})$ as before.
%some of the points come together.  *
%As before, we only need to consider the component of 
%$\mathcal{B}$ where one (or more) of the $q_i$ become equal to the last point $q_{\delta+1}$. 
%First of all we note that the whole of 
%$\mathcal{B}$ is not relevant while computing the contribution to the Euler class; only the points 
%at which the section vanishes is relevant. 
%First, we note that $\B$ contains $\overline{\A_1^{\delta}\circ \PP \D_5}$ and 
%the section $\Psi_{\PP \A_5}$ vanishes on it with a multiplicity of $2$. 
%Thus, the contribution from $\overline{\A_1^{\delta}\circ \PP \D_5}$ is
%\bgd
%2 \N(\A_1^{\delta}\PP \D_5,n_1,m_1,m_2,\theta).
%\edd
%Now we focus on the component of 
%$\mathcal{B}$ where one (or more) of the $q_i$ become equal to the last point $q_{\delta+1}$. 
%Let us define 
%\begin{align*}
%\B (q_1, l_{q_{\delta +1 }}) &:= \{ ([f], q_1, \ldots, q_{\delta}, l_{q_{\delta +1}}) \in \B: q_1 = q_{\delta +1} \}
%\end{align*}
%The spaces $\B (q_2, q_{\delta +1 }), \ldots, \B (q_{\delta}, q_{\delta +1 })$ are defined 
%similarly. 
%The claim is that 
%{\bf What about $\PP \D_7^s$?}\\
We show in \cite[Lemma 6.3, statement 5]{BM13_2pt_published}  
%(Lemma *) that, 
that,
\begin{align*}
\B (q_1, l_{q_{\delta +1 }}) & 
\approx \overline{\A_1^{\delta-1} \circ \PP \A}_6\cup \overline{\A_1^{\delta-1} \circ \PP \E}_6 \cup 
\overline{\A_1^{\delta-1} \circ \PP D}_7.
\end{align*}
%\subsection{RR}
The contribution from $\overline{\A_1^{\delta-1} \circ \PP \A}_6 \cap \mu$ 
and $\overline{\A_1^{\delta-1} \circ \PP \E}_6 \cap \mu$ are $2$ and $5$ 
respectively. Since the pseudocycle $\mu$ is generic, it will not intersect 
$\overline{\A_1^{\delta-1} \circ \PP \D}_7 \cap \mu$. 
Hence, the total contribution 
%from the points of $\B (q_1, l_{q_{\delta +1 }}) \cap \mu$ of the 
%first and second type are $2$ and $5$ respectively. 
%Hence the total contribution 
from all the components of type $\B (q_{i_1}, l_{q_{\delta +1 }})$ equals
\bgd
2\binom{\delta}{1}N(\A_1^{\delta-1}\PP\A_6,n_1,m_1,m_2,\theta)
+5\binom{\delta}{1}N(\A_1^{\delta-1}\PP\E_6,n_1,m_1,m_2,\theta).
\edd
Next, 
%let us define 
%\begin{align*}
%\B (q_1, q_2, l_{q_{\delta +1 }}) &:= 
%\{ ([f], q_1, \ldots, q_{\delta}, l_{q_{\delta +1}}) \in \B: q_1=q_2 = l_{q_{\delta +1}}. \}
%\end{align*}
%we claim that 
we show in \cite{BM_closure_of_seven_points} that 
\begin{align*}
\B (q_1, q_2, l_{q_{\delta +1 }}) & \approx \overline{\A_1^{\delta-2}\circ \PP \D_7}.  
\end{align*}
This is geometrically believable due to the following picture: 
\newpage
\begin{figure}[h!]
\vspace*{0.2cm}
\begin{center}\includegraphics[scale = 0.5]{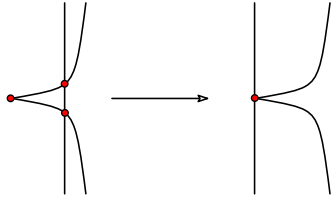}\vspace*{-0.2cm}\end{center}
\caption{Two nodes and  nodes and one $A_4$-singularity collapsing to a $D_7$-singularity}
\end{figure} 
We also show that the contribution from each of the points of 
$\overline{\A_1^{\delta-2}\circ \PP \D_7} \cap \mu$ is $6$. 
%$\B (q_1, q_2, l_{q_{\delta +1 }}) \cap \mu$ is $6$.
Hence the total contribution from all the components of type $\B (q_{i_1}, q_{i_2}, l_{q_{\delta +1 }})$ equals
\begin{align*}
6 \binom{\delta}{2} \N(\A_1^{\delta-2} \PP\D_7, n_1, m_1, m_2,\theta).   
\end{align*}
Finally, we 
%claim that
we show in \cite{BM_closure_of_seven_points} that
\begin{align*}
\B(q_1, q_2, q_3, l_{q_{4}}) \cap \mu & = \varnothing. 
%\qquad \forall ~ 0 \leq r \leq 3.  
\end{align*}
Hence the remaining stratum of $\B$ does not contribute to the Euler class, 
giving us Theorem \ref{a1_delta_pa5}. \qed

\subsection{Proof of Theorem \ref{a1_delta_pa6}: computing $\N(\A_1^{\delta} \PP \A_6, n_1, m_1, m_2, \theta)$}
\hf \hf Recall that in section \ref{notation}, we have defined the space  
\begin{align*}
\A^{\delta}_1 \circ \overline{\PP \A}_5 := \{ ([f], q_1, \ldots, q_{\delta}, 
l_{q_{\delta+1}}) \in \DD \times (X)^{\delta} \times \mathbb{P} T \X: 
&\textnormal{$f$ has a singularity of type $\A_1$ at $q_1, \ldots, q_{\delta}$}, \\ 
                   & ([f], l_{q_{\delta+1}}) \in \PP \A_5, ~~\textnormal{$q_1, \ldots, q_{\delta+1}$ all distinct}\}. 
\end{align*}
Let $\mu$ be a generic pseudocycle representing the homology class 
Poincar\'{e} dual 
to 
\[c_1^{n_1} x_1^{m_1} x_2^{m_2} \lambda^{\theta}y^{\delta_L - (n_1 + m_1 + 2 m_2 +\theta+\delta+6)}.\]
As before, we let $v \in \hat{\gamma}$ and $w \in \pi^*TX/\hat{\gamma}$ be two fixed non zero vectors 
and define $f_{ij}$ as before. 
We now define a section of the following bundle 
\begin{align}
\Psi_{\PP \A_6}:  \A^{\delta}_1 \circ \PP \A_5 
\lra \mathbb{L}_{\PP \A_6} & := \hat{\gamma}^{\ast 6}\otimes\gamma_{\D}^{\ast 4}
\otimes(\pi^*TX/\hat{\gamma})^{\ast 6} \otimes L^{\ast 4} \qquad \textnormal{given by} \nonumber \\
\{\Psi_{\PP \A_6}([f],q_1,\ldots,q_\delta, l_{q_{\delta+1}})\}
(f^{\otimes 2}\otimes v^{\otimes 6}\otimes w^{\otimes 6}) &:= f_{02}^3 A_6^f, \nonumber
%\psi_{\A_1} : \psi_{\A_0}^{-1}(0) \lra \mathcal{V}_{\A_1} &:= \gamma_{\D}^*\otimes T^*\X \otimes L, 
%\qquad \textnormal{given by} \qquad 
%\{\psi_{\A_1}([f], q)\}(f) := \nabla f|_q. \label{psi_a0_a1_section_defn}
\end{align}
where 
\[ \A^{f}_6 := f_{60}- \f{ 15 f_{21} f_{41}}{f_{02}}-\f{10 f_{31}^2}{f_{02}} + \f{60 f_{12} f_{21} f_{31}}{f_{02}^2}
   +
   \f{45 f_{21}^2 f_{22}}{f_{02}^2} - \f{15 f_{03} f_{21}^3}{f_{02}^3}
   -\f{90 f_{12}^2 f_{21}^2}{f_{02}^3} .\]
%In \cite{RM_Hypersurfaces} we show that 
In 
%Lemma $7.1$ statement $*$ of the arXiv version 
\cite[Lemma 7.1, statement 11]{BM13}, we prove that 
%As stated earlier ()
\begin{align}
\overline{\PP A}_5 & = \PP A_5 \cup \overline{\PP A}_6 \cup \overline{\PP D}_6 \cup \overline{\PP E}_6.  \label{pa5_closure}
\end{align}
Hence, let us define 
\begin{align*}
\mathcal{B} &:= \overline{A_1^{\delta} \circ \overline{\PP A}}_5 - \A^{\delta}_1 \circ (\PP\A_5\cup \overline{\PP A}_6). 
\end{align*}
If $L$ is sufficiently $(2\delta+5)$-ample, then the section 
$\Psi_{\PP \A_6}$ vanishes on the points of  
$\A_1^{\delta}\circ \PP A_6$ transversally 
(all the points of $\A_1^{\delta}\circ \PP A_6$ are smooth points of the 
variety $\A_1^{\delta}\circ \overline{\PP A}_5$).
Hence, 
\begin{align*}
\lan e(\mathbb{L}_{\PP \A_6}), 
~~[\overline{\A^{\delta}_1 \circ \overline{\PP\A}}_5] \cap [\mu] \ran & = \N(\A_1^{\delta}\PP \A_6, n_1, m_1, m_2, \theta)  
+ \mathcal{C}_{\mathcal{B}\cap \mu}.
%\label{na1_pseudocycle}
\end{align*}
Let us define $\B_0$ as before, namely the component of the boundary where 
all the points are distinct.
%none of the 
%$(q_1, \ldots, q_{\delta})$ points are equal to the last marked point $q_{\delta+1}$. 
We show in \cite[Corollary 7.7 and Corollary 7.9]{BM13}
%and \cite[Corollary 7.9]{BM13} 
%, Lemma *, 
that the section $\Psi_{\PP A_6}$ 
vanishes on all points of $\overline{A_1^{\delta}}\circ \PP \D_6 \cap \mu$ 
and $\overline{A_1^{\delta}}\circ \PP \E_6 \cap \mu$
with a multiplicity 
of $4$ and $3$ respectively. 
%We also show that 
%The section $\Psi_{\PP D_5}$ does not vanish on $\PP D_5^{\vee}$ by definition 
%(recall that by definition of the space $\PP D_5^{\vee}$, the directional derivative 
%$f_{21} \neq 0$). 
%Hence, the section $\Psi_{\PP D_5}$ does not vanish 
%on $\overline{A^{\delta}_1} \circ \overline{\PP D^{\vee}_5} \cap \mu$ for a generic $\mu$. 
Hence, the contribution of $\B_0$ 
to the Euler class is 
\bgd
4 \N(\A_1^{\delta}\PP \D_6,n_1,m_1,m_2,\theta)+3 \N(\A_1^{\delta}\PP \E_6,n_1,m_1,m_2,\theta).
\edd
%\begin{align*}
%2 \Num(\A_1^{\delta}\PP \D_5, n_1, m_1, m_2, \theta). 
%\end{align*}
%does not contribute to the Euler class. \\ 
%*****************
%If $L$ is sufficiently $(2\delta+5)$-ample, then this section is transverse to the zero set. 
%Next, let us define 
%\begin{align*}
%\mathcal{B} &:= \overline{\A^{\delta}_1 \circ \PP \A}_4- \A^{\delta}_1 \circ \PP\A_4. 
%\end{align*}
%Hence 
%\begin{align}
%\lan e(\mathbb{L}_{\PP \A_5}), 
%~~[\overline{\A^{\delta}_1 \circ \PP\A}_4] \cap [\mu] \ran & = \N(\A_1^{\delta}\PP \A_5, n_1, m_1, m_2, \theta)  
%+ \mathcal{C}_{\mathcal{B}\cap \mu}.
%\label{na1_pseudocycle}
%\end{align}
%We now give an explicit description of $\mathcal{B}$.
%\hf \hf Next, let us analyze the boundary when 
\hf \hf Next, let us analyze the boundary when 
one (or more) of the $q_i$ become equal to the last point $q_{\delta+1}$.
Define $\mathcal{B}(q_{i_1}, \ldots q_{i_k}, l_{q_{\delta}})$ as before.
We show in \cite[Lemma 6.3, statement 6]{BM13_2pt_published} 
%Lemma * 
and \cite{BM_closure_of_seven_points} that 
\begin{align*}
\B (q_1, l_{q_{\delta +1 }}) & 
\approx \overline{\A_1^{\delta-1} \circ \PP \A}_7\cup \overline{\A_1^{\delta-1} \circ \PP \E}_7 \cup 
\overline{\A_1^{\delta-1} \circ \PP D}_8 \cup \overline{\A_1^{\delta-1} \circ \hat{X}}_9.
\end{align*}
%{\bf What about $\PP \D_8^s$?}\\
We also show in  \cite[Corollary 6.13 and 6.23]{BM13_2pt_published} 
%and \cite{BM_closure_of_seven_points} 
that 
the contribution from the points of 
$\overline{\A_1^{\delta-1} \circ \PP \A}_7\cap \mu$ 
and $\overline{\A_1^{\delta-1} \circ \PP \E}_7 \cap \mu$ are $2$ and $6$ 
respectively. 
We also show in \cite{BM_closure_of_seven_points} that the contribution from the points of 
$\overline{\A_1^{\delta-1} \circ \hat{X}}_9\cap \mu$ is $12$. 
Since the pseudocycle $\mu$ is generic, 
%the intersection 
it will not intersect $\overline{\A_1^{\delta-1} \circ \PP \D}_8$.
%$\overline{\A_1^{\delta-1} \circ \PP \D}_8 \cap \mu$. 
%Hence, the total contribution
%$\B (q_1, l_{q_{\delta +1 }}) \cap \mu$ of the 
%first and second type are $2$ and $6$ respectively. 
Hence the total 
contribution from all the components of type $\B (q_{i_1}, l_{q_{\delta +1 }})$ equals
\bgd
2\binom{\delta}{1}N(\A_1^{\delta-1}\PP\A_7,n_1,m_1,m_2,\theta)+
6\binom{\delta}{1}N(\A_1^{\delta-1}\PP\E_7,n_1,m_1,m_2,\theta) + 12\binom{\delta}{1}N(\A_1^{\delta-1} \hat{X}_9, n_1, 
m_1, m_2, \theta).
\edd
Finally, we show in \cite{BM_closure_of_seven_points} that 
\begin{align*}
\B(q_1, q_2, l_{q_{3}}) \cap \mu & = \varnothing. 
%\qquad \forall ~ 2 \leq r \leq 7.  
\end{align*}
Hence the remaining stratum of $\B$ does not contribute to the Euler class, 
%, since its intersection with $\mu$ 
%is empty. This 
giving us Theorem \ref{a1_delta_pa6}. \qed

\subsection{Proof of Theorem \ref{a1_delta_pd6} and \ref{a1_delta_pe6}: 
computing $\N(\A_1^{\delta} \PP \D_6, n_1, m_1, m_2, \theta)$ and $\N(\A_1^{\delta} \PP \E_6, n_1, m_1, m_2, \theta)$}
\hf\hf The computation of 
$\N(\A_1^{\delta} \PP \D_6, n_1, m_1, m_2, \theta)$ and $\N(\A_1^{\delta} \PP \E_6, n_1, m_1, m_2, \theta)$ 
are over the same spaces; hence it is more efficient to do these to computations simultaneously. 
Recall that in section \ref{notation}, we have defined the following space  
\begin{align*}
\A^{\delta}_1 \circ \overline{\PP \D}_5 := \{ ([f], q_1, \ldots, q_{\delta}, 
l_{q_{\delta+1}}) \in \DD \times (X)^{\delta} \times \mathbb{P} T \X: 
&\textnormal{$f$ has a singularity of type $\A_1$ at $q_1, \ldots, q_{\delta}$}, \\ 
                   & ([f], l_{q_{\delta+1}}) \in \overline{\PP \D}_5, ~~\textnormal{$q_1, \ldots, q_{\delta+1}$ all distinct}\}. 
\end{align*}
Let $\mu$ be a generic pseudocycle representing the homology class 
Poincar\'{e} dual 
to 
\[c_1^{n_1} x_1^{m_1} x_2^{m_2} \lambda^{\theta}y^{\delta_L - (n_1 + m_1 + 2 m_2 +\theta+\delta+6)}.\]
As before, we let $v \in \hat{\gamma}$ and $w \in \pi^*TX/\hat{\gamma}$ be two fixed non zero vectors 
and define $f_{ij}$ as before. 
We now define sections of the following two bundles: 
\begin{align}
\Psi_{\PP \D_6}:  \A^{\delta}_1 \circ \overline{\PP \D}_5 
\lra \mathbb{L}_{\PP \D_6} & := \gamma_{\D}^*\otimes \hat{\gamma}^{\ast 4}\otimes L, \qquad \textnormal{given by} 
\nonumber \\
\{\Psi_{\PP \D_6}([f], q_1,\ldots,q_\delta,l_{q_{\delta+1}})\}(f\otimes v^{\otimes 4}) &:= f_{40} \nonumber 
\qquad \textnormal{and}\\
%\psi_{\A_1} : \psi_{\A_0}^{-1}(0) \lra \mathcal{V}_{\A_1} &:= \gamma_{\D}^*\otimes T^*\X \otimes L, 
%\qquad \textnormal{given by} \qquad 
%\{\psi_{\A_1}([f], q)\}(f) := \nabla f|_q. \label{psi_a0_a1_section_defn}
%\end{align}
%and a section of the following bundle 
%\begin{align}
\Psi_{\PP \E_6}:  \A^{\delta}_1 \circ \overline{\PP \D}_5 
\lra \mathbb{L}_{\PP \E_6} & := \hat{\gamma}^{\ast}\otimes\gamma_{\D}^{\ast }\otimes(TX/\hat{\gamma})^{\ast 2} 
\otimes L^{\ast}, \qquad \textnormal{given by} \nonumber \\
\{\Psi_{\PP \E_6}([f],q_1,\ldots,q_\delta, l_{q_{\delta+1}})\}(f\otimes v\otimes w^{\otimes 2}) &:= f_{12}.  \nonumber
%\psi_{\A_1} : \psi_{\A_0}^{-1}(0) \lra \mathcal{V}_{\A_1} &:= \gamma_{\D}^*\otimes T^*\X \otimes L, 
%\qquad \textnormal{given by} \qquad 
%\{\psi_{\A_1}([f], q)\}(f) := \nabla f|_q. \label{psi_a0_a1_section_defn}
\end{align}
%In \cite{RM_Hypersurfaces} we show that 
In 
%Lemma $7.1$ statement $*$ of the arXiv version 
\cite[Lemma 7.1, statement 6]{BM13}, we prove that 
%As stated earlier ()
\begin{align}
\overline{\PP D}_5 & = \PP D_5 \cup \overline{\PP D}_6 \cup \overline{\PP E}_6.  \label{pd5_closure}
\end{align}
Hence, let us define 
\begin{align*}
\mathcal{B} &:= \overline{A_1^{\delta} \circ \overline{\PP D}}_5 - \A^{\delta}_1 \circ (\PP\D_5\cup \overline{\PP D}_6) 
\qquad \textnormal{and} \\
\B^{\prime} &:= 
\overline{A_1^{\delta} \circ \overline{\PP D}}_5 - \A^{\delta}_1 \circ (\PP\D_5\cup \overline{\PP E}_6)
\end{align*}
If $L$ is sufficiently $(2\delta+4)$-ample and $(2\delta+3)$-ample, then the sections 
$\Psi_{\PP \D_6}$ and $\Psi_{\PP \E_6}$ vanish on the points of 
$\A_1^{\delta}\circ \PP D_6$  and $\A_1^{\delta}\circ \PP E_6$ transversally, respectively 
(all the points of $\A_1^{\delta}\circ \PP D_6$ and $\A_1^{\delta}\circ \PP E_6$ are smooth points of the 
variety $\A_1^{\delta}\circ \overline{\PP D}_5$).
Hence 
\begin{align*}
\lan e(\mathbb{L}_{\PP \D_6}), 
~~[\overline{\A^{\delta}_1 \circ \PP\D_5}] \cap [\mu] \ran & = \N(\A_1^{\delta}\PP \D_6, n_1, m_1, m_2, \theta)  
+ \mathcal{C}_{\mathcal{B}\cap \mu}(\Psi_{\PP \D_6}) \qquad \textnormal{and} \\
%\label{na1_pseudocycle}
%\end{align}
%Hence 
%\begin{align}
\lan e(\mathbb{L}_{\PP \E_6}), 
~~[\overline{\A^{\delta}_1 \circ \PP\D_5}] \cap [\mu] \ran & = \N(\A_1^{\delta}\PP \E_6, n_1, m_1, m_2, \theta)  
+ \mathcal{C}_{\mathcal{B}^{\prime}\cap \mu}(\Psi_{\PP \E_6}).
%\label{na1_pseudocycle}
\end{align*}
%Hence, 
%\begin{align}
%\lan e(\mathbb{L}_{\PP \A_6}), 
%~~[\overline{\A^{\delta}_1 \circ \overline{\PP\A}}_5] \cap [\mu] \ran & = \N(\A_1^{\delta}\PP \A_6, n_1, m_1, m_2, \theta)  
%+ \mathcal{C}_{\mathcal{B}\cap \mu}.
%\label{na1_pseudocycle}
%\end{align}
Here $\mathcal{C}_{\mathcal{B}\cap \mu}(\Psi_{\PP \D_6})$ 
and $\mathcal{C}_{\mathcal{B}^{\prime}\cap \mu}(\Psi_{\PP \E_6})$ 
denote the contribution from the boundary to the Euler class from the sections 
$\Psi_{\PP \D_6}$ and $\Psi_{\PP \E_6}$ respectively. \\ 
\hf \hf Let us define $\B_0$ and $\B_0^{\prime}$ as before, 
namely the component of the boundary where 
all the points are distinct. 
%none of the 
%$(q_1, \ldots, q_{\delta})$ points are equal to the last marked point $q_{\delta+1}$. 
It is easy to see that the section $\Psi_{\PP D_6}$ 
does not vanish on $\B^{\prime}_0 \cap \mu$ 
and the section $\Psi_{\PP E_6}$ does not vanish on $\B_0\cap \mu$. 
Hence, $\B_0$ and $\B^{\prime}_0$ do not contribute to the Euler class. \\ 
\hf \hf Next, let us analyze the boundary when 
one (or more) of the $q_i$ become equal to the last point $q_{\delta+1}$.
Define $\mathcal{B}(q_{i_1}, \ldots q_{i_k}, l_{q_{\delta}})$ and 
$\mathcal{B}^{\prime}(q_{i_1}, \ldots q_{i_k}, l_{q_{\delta}})$
as before. We show in 
\cite[Lemma 6.3, statement 9]{BM13_2pt_published} 
and \cite{BM_closure_of_seven_points} 
%Lemma *, 
that 
\begin{align*}
\B (q_1, l_{q_{\delta +1 }}) \approx \B^{\prime} (q_1, l_{q_{\delta +1 }}) & \approx 
\overline{\A_1^{\delta-1} \circ \PP \D_7}\cup \overline{\A_1^{\delta-1} \circ \PP \E_7}
\cup \overline{\A_1^{\delta-1} \circ \hat{X}}_9.
\end{align*}
Let us first focus on the section $\Psi_{\PP \D_6}$. 
We show in \cite[Corollary 6.28 and Corollary 6.30]{BM13_2pt_published}, 
%Corollary *, 
that 
the contribution of the section $\Psi_{\PP \D_6}$  
from the points of 
$\overline{\A_1^{\delta-1} \circ \PP \D_7}\cap \mu$ 
and $\overline{\A_1^{\delta-1} \circ \PP \E_7}\cap \mu$ 
are $2$ and $1$ respectively. 
Let us now consider 
$\overline{\A_1^{\delta-1} \circ \hat{X}}_9\cap \mu$. 
Suppose 
\[ ([f], q_1,\ldots,q_\delta,l_{q_{\delta+1}}) \in  \overline{\A_1^{\delta-1} \circ \hat{X}}_9\cap \mu.\]
Since $\mu$ is a generic pseudocycle, we conclude that $f_{40} \neq 0$. Hence, the section 
$\Psi_{\PP \D_6}$ does not vanish on $\overline{\A_1^{\delta-1} \circ \hat{X}}_9\cap \mu$. 
%$\B (q_1, l_{q_{\delta +1 }}) \cap \mu$ of the 
%first and second type are $2$ and $1$ respectively. 
Hence the total contribution of the section $\Psi_{\PP \D_6}$ from all the 
components of type $\B (q_{i_1}, l_{q_{\delta +1 }})$ equals
\bgd
2\binom{\delta}{1}N(\A_1^{\delta-1}\PP\D_7,n_1,m_1,m_2,\theta)+
\binom{\delta}{1}N(\A_1^{\delta-1}\PP\E_7,n_1,m_1,m_2,\theta).
\edd
Next, we show in \cite{BM_closure_of_seven_points} that  
\begin{align*}
\B(q_1, q_2, l_{q_{3}}) \cap \mu & = \varnothing. 
%\qquad \forall ~ 2 \leq r \leq 7.  
\end{align*}
Hence the remaining stratum of $\B$ does not contribute to the Euler class, 
giving us Theorem \ref{a1_delta_pd6}. \\ 
\hf \hf Next, let us consider the section $\Psi_{\PP \E_6}$.
We show in \cite[Corollary 6.30]{BM13_2pt_published}, 
%Corollary * 
and in \cite{BM_closure_of_seven_points}  
that the contribution of the section $\Psi_{\PP \E_6}$  
from the points of 
%$\overline{\A_1^{\delta-1} \circ \PP \D_7}\cap \mu$, 
$\overline{\A_1^{\delta-1} \circ \PP \E_7}\cap \mu$ and 
$\overline{\A_1^{\delta-1} \circ \hat{X}}_9 \cap \mu$
are $1$  and $4$ respectively. The section 
$\Psi_{\PP \E_6}$ does not vanish on 
$\overline{\A_1^{\delta-1} \circ \PP \D}_7$; hence it does not vanish on 
$\overline{\A_1^{\delta-1} \circ \overline{\PP \D}}_7 \cap \mu$. 
%Let us now consider 
%$\overline{\A_1^{\delta-1} \circ \PP \D}_7\cap \mu$. 
%Suppose 
%\[ ([f], q_1,\ldots,q_\delta,l_{q_{\delta+1}}) \in  \overline{\A_1^{\delta-1} \circ \PP \D}_7\cap \mu.\]
%Since $\mu$ is a generic pseudocycle, we conclude that 
%\[ ([f], q_1,\ldots,q_\delta,l_{q_{\delta+1}}) \in  \A_1^{\delta-1} \circ \PP \D_7\cap \mu.\]
%Since $f_{12} \neq 0$ on $\PP \D_7$, the section $\Psi_{\PP \E_6}$ does not vanish on 
%$\overline{\A_1^{\delta-1} \circ \PP \D}_7\cap \mu$. 
Hence the total contribution from all the components of type $\B (q_{i_1}, l_{q_{\delta +1 }})$ equals
\bgd
1\times \binom{\delta}{1}N(\A_1^{\delta-1}\PP\E_7,n_1,m_1,m_2,\theta) + 
4\binom{\delta}{1}N(\A_1^{\delta-1}\hat{X}_9,n_1,m_1,m_2,\theta).
\edd
Finally, as stated earlier, we prove in \cite{BM_closure_of_seven_points} that 
\begin{align*}
\B(q_1, q_2, l_{q_{3}}) \cap \mu & = \varnothing.
%\qquad \forall ~ 2 \leq r \leq 7.  
\end{align*}
Hence, the remaining stratum of $\B$ does not contribute to Euler class, 
giving us Theorem \ref{a1_delta_pe6}. \qed 

\subsection{Proof of Theorem \ref{a1_delta_pd7} and \ref{a1_delta_pe7}: 
computing $\N(\A_1^{\delta} \PP \D_7, n_1, m_1, m_2, \theta)$ and 
$\N(\A_1^{\delta} \PP \E_7, n_1, m_1, m_2, \theta)$}
\hf\hf The computation of 
$\N(\A_1^{\delta} \PP \D_7, n_1, m_1, m_2, \theta)$ and $\N(\A_1^{\delta} \PP \E_7, n_1, m_1, m_2, \theta)$ 
are over the same spaces; hence it is more efficient to prove these two computations simultaneously. 
Recall that in section \ref{notation}, we have defined the following space  
%Let us define the space  
\begin{align*}
\A^{\delta}_1 \circ \overline{\PP \D}_6 := \{ ([f], q_1, \ldots, q_{\delta}, 
l_{q_{\delta+1}}) \in \DD \times (X)^{\delta} \times \mathbb{P} T \X: 
&\textnormal{$f$ has a singularity of type $\A_1$ at $q_1, \ldots, q_{\delta}$}, \\ 
                   & ([f], l_{q_{\delta+1}}) \in \overline{\PP \D}_6, 
                   ~~\textnormal{$q_1, \ldots, q_{\delta+1}$ all distinct}\}. 
\end{align*}
Let $\mu$ be a generic pseudocycle representing the homology class 
Poincar\'{e} dual 
to 
\[c_1^{n_1} x_1^{m_1} x_2^{m_2} \lambda^{\theta}y^{\delta_L - (n_1 + m_1 + 2 m_2 +\theta+\delta+7)}.\]
As before, we let $v \in \hat{\gamma}$ and $w \in \pi^*TX/\hat{\gamma}$ be two fixed non zero vectors 
and define $f_{ij}$ as before. 
We now define sections of the following two bundles 
\begin{align}
\Psi_{\PP \D_7}:  \A^{\delta}_1 \circ \overline{\PP \D}_6 
\lra \mathbb{L}_{\PP \D_7} & := 
\hat{\gamma}^{\ast 6}\otimes\gamma_{\D}^{\ast 2}\otimes(\pi^*TX/\hat{\gamma})^{\ast 2} \otimes L^{\ast 2} , 
\qquad \textnormal{given by} \nonumber \\
\{\Psi_{\PP \D_7}([f], q_1,\ldots,q_\delta, l_{q_{\delta+1}})\}
(f^{\otimes 2}\otimes v^{\otimes 6}\otimes w^{\otimes 2}) &:= f_{12} D_7^f \qquad \textnormal{and}  \nonumber \\
%\psi_{\A_1} : \psi_{\A_0}^{-1}(0) \lra \mathcal{V}_{\A_1} &:= \gamma_{\D}^*\otimes T^*\X \otimes L, 
%\qquad \textnormal{given by} \qquad 
%\{\psi_{\A_1}([f], q)\}(f) := \nabla f|_q. \label{psi_a0_a1_section_defn}
%\end{align}
%and 
%\begin{align}
\Psi_{\PP \E_7}:  \A^{\delta}_1 \circ \overline{\PP \D}_6 
\lra \mathbb{L}_{\PP \E_7} & := \gamma_{\D}^{\ast}\otimes \hat{\gamma}^{\ast} 
\otimes(\pi^*TX/\hat{\gamma})^{\ast 2} \otimes L , \qquad 
\textnormal{given by} \nonumber \\
\{\Psi_{\PP \E_7}([f], q_1,\ldots,q_\delta, l_{q_{\delta+1}})\}(f\otimes v \otimes w^{\otimes 2}) &:= f_{12},  \nonumber
%\psi_{\A_1} : \psi_{\A_0}^{-1}(0) \lra \mathcal{V}_{\A_1} &:= \gamma_{\D}^*\otimes T^*\X \otimes L, 
%\qquad \textnormal{given by} \qquad 
%\{\psi_{\A_1}([f], q)\}(f) := \nabla f|_q. \label{psi_a0_a1_section_defn}
\end{align}
where $~\D^{f}_7 := f_{50} - \frac{5 f_{31}^2}{3 f_{12}}.$\\
\hf \hf In 
%Lemma $7.1$ statement $*$ of the arXiv version 
\cite[Lemma 7.1, statement 7]{BM13}, we prove that 
%As stated earlier ()
\begin{align}
\overline{\PP D}_6 & = \PP D_6 \cup \overline{\PP D}_7 \cup \overline{\PP E}_7.  \label{pd6_closure}
\end{align}
Hence, let us define 
\begin{align*}
\mathcal{B} &:= \overline{A_1^{\delta} \circ \overline{\PP D}}_6 - \A^{\delta}_1 \circ (\PP\D_6\cup \overline{\PP D}_7) 
\qquad \textnormal{and} \\
\B^{\prime} &:= 
\overline{A_1^{\delta} \circ \overline{\PP D}}_5 - \A^{\delta}_1 \circ (\PP\D_6\cup \overline{\PP E}_7)
\end{align*}
If $L$ is sufficiently $(2\delta+5)$-ample and $(2\delta+4)$-ample, then the sections 
$\Psi_{\PP \D_7}$ and $\Psi_{\PP \E_7}$ vanish on the points of 
$\A_1^{\delta}\circ \PP D_7$  and $\A_1^{\delta}\circ \PP E_7$ transversally, respectively 
(all the points of $\A_1^{\delta}\circ \PP D_7$ and $\A_1^{\delta}\circ \PP E_7$ are smooth points of the 
variety $\A_1^{\delta}\circ \overline{\PP D}_6$).
Hence 
\begin{align*}
\lan e(\mathbb{L}_{\PP \D_7}), 
~~[\overline{\A^{\delta}_1 \circ \PP\D_6}] \cap [\mu] \ran & = \N(\A_1^{\delta}\PP \D_7, n_1, m_1, m_2, \theta)  
+ \mathcal{C}_{\mathcal{B}\cap \mu}(\Psi_{\PP \D_7}) \qquad \textnormal{and} \\
%\label{na1_pseudocycle}
%\end{align}
%Hence 
%\begin{align}
\lan e(\mathbb{L}_{\PP \E_7}), 
~~[\overline{\A^{\delta}_1 \circ \PP\D_6}] \cap [\mu] \ran & = \N(\A_1^{\delta}\PP \E_7, n_1, m_1, m_2, \theta)  
+ \mathcal{C}_{\mathcal{B}^{\prime} \cap \mu}(\Psi_{\PP \E_7}).
%\label{na1_pseudocycle}
\end{align*}
%Hence, 
%\begin{align}
%\lan e(\mathbb{L}_{\PP \A_6}), 
%~~[\overline{\A^{\delta}_1 \circ \overline{\PP\A}}_5] \cap [\mu] \ran & = \N(\A_1^{\delta}\PP \A_6, n_1, m_1, m_2, \theta)  
%+ \mathcal{C}_{\mathcal{B}\cap \mu}.
%\label{na1_pseudocycle}
%\end{align}
%Here $\mathcal{C}_{\mathcal{B}\cap \mu}(\Psi_{\PP \D_6})$ 
%and $\mathcal{C}_{\mathcal{B}^{\prime}\cap \mu}(\Psi_{\PP \E_6})$ 
%denote the contribution from the boundary to the Euler class from the sections 
%$\Psi_{\PP \D_6}$ and $\Psi_{\PP \E_6}$ respectively. \\ 
\hf \hf Let us define $\B_0$ and $\B_0^{\prime}$ as before, 
namely the component of the boundary where 
all the points are distinct.
%none of the 
%$(q_1, \ldots, q_{\delta})$ points are equal to the last marked point $q_{\delta+1}$. 
Note that the section $\Psi_{\PP D_7}$ 
does not vanish on $\overline{A^{\delta}_1}\circ \PP E_7$. To see why that is 
so, we note that on $\PP E_7$, the directional derivative $f_{12} =0$ but 
$f_{03}$ and $f_{31} \neq 0$; this follows from Lemma \ref{fstr_prp_E7}. 
%\subsection{RR}
Hence, 
\[ f_{12}D^f_7 = f_{12}f_{50} - \frac{5}{3} f_{31}^2 = -\frac{5}{3} f_{31}^2 \neq 0. \]
Hence, the section $\Psi_{\PP D_7}$ does not vanish on 
$\B_0 \cap \mu$. It is also immediate that the section $\Psi_{\PP E_7}$ does not 
vanish on $\overline{A^{\delta}_1}\circ \PP D_7$; hence it does not vanish 
on $\B^{\prime}_0 \cap \mu$. 
%and the section $\Psi_{\PP E_6}$ does not vanish on $\B_0\cap \mu$. 
Hence, $\B_0$ and $\B^{\prime}_0$ do not contribute to the Euler class. \\ 
\hf \hf Next, let us analyze the boundary when 
one (or more) of the $q_i$ become equal to the last point $q_{\delta+1}$. 
%$$$$$$$$$$$$$$$!!!!!!!!!!!!!!!!!!!!!!
Define $\mathcal{B}(q_{i_1}, \ldots q_{i_k}, l_{q_{\delta}})$ and 
$\mathcal{B}^{\prime}(q_{i_1}, \ldots q_{i_k}, l_{q_{\delta}})$
as before. We show in \cite{BM_closure_of_seven_points}, that 
%***************************************
%*In \cite{RM_Hypersurfaces} we show that 
%If $L$ is sufficiently $(2\delta+7)$-ample, then these sections are transverse to the zero set. 
%Next, let us define 
%\begin{align*}
%\mathcal{B} &:= \overline{\A^{\delta}_1 \circ \PP \D_6}- \A^{\delta}_1 \circ \PP\D_6. 
%\end{align*}
%Hence 
%\begin{align}
%\lan e(\mathbb{L}_{\PP \D_7}), 
%~~[\overline{\A^{\delta}_1 \circ \PP\D_6}] \cap [\mu] \ran & = \N(\A_1^{\delta}\PP \D_7, n_1, m_1, m_2, \theta)  
%+ \mathcal{C}_{\mathcal{B} \cap \mu}(\Psi_{\PP \D_7}) 
%\label{na1_pseudocycle}
%\end{align}
%and 
%\begin{align}
%\lan e(\mathbb{L}_{\PP \E_7}), 
%~~[\overline{\A^{\delta}_1 \circ \PP\E_6}] \cap [\mu] \ran & = \N(\A_1^{\delta}\PP \E_7, n_1, m_1, m_2, \theta)  
%+ \mathcal{C}_{\mathcal{B} \cap \mu}(\Psi_{\PP \E_7}).
%\label{na1_pseudocycle}
%\end{align}
%We now give an explicit description of $\mathcal{B}$. 
%As before, we only need to consider the component of 
%$\mathcal{B}$ where one (or more) of the $q_i$ become equal to the last point $q_{\delta+1}$. 
%Define $\mathcal{B}(q_{i_1}, \ldots q_{i_k}, l_{q_{\delta}})$ as before. We claim 
%that 
\begin{align*}
\mathcal{B}(q_1, l_{q_{\delta+1}}) \approx \mathcal{B}^{\prime}(q_1, l_{q_{\delta+1}}) & \approx  
\overline{\A_1^{\delta-1} \circ \PP \D_8}\cup \overline{\A_1^{\delta-1} \circ \PP \X_9}.
\end{align*} 
%First, let us look at the section 
Let us first focus on the section $\Psi_{\PP \D_7}$. 
We show that the contribution of the section $\Psi_{\PP \D_7}$  
from the points of 
$\overline{\A_1^{\delta-1} \circ \PP \D_8}\cap \mu$ 
is $2$.
%and $\overline{\A_1^{\delta-1} \circ \PP \E_7}\cap \mu$ 
%are $2$ and $1$ respectively. 
Next, let us consider 
$\overline{\A_1^{\delta-1} \circ \PP X}_9\cap \mu$. 
Suppose 
\[ ([f], q_1,\ldots,q_\delta,l_{q_{\delta+1}}) \in  \overline{\A_1^{\delta-1} \circ \PP X}_9\cap \mu.\]
Since $\mu$ is a generic pseudocycle, we conclude that $f_{31} \neq 0$. Hence, the section 
$\Psi_{\PP \D_7}$ does not vanish on $\overline{\A_1^{\delta-1} \circ \PP X}_9\cap \mu$. 
%$\B (q_1, l_{q_{\delta +1 }}) \cap \mu$ of the 
%first and second type are $2$ and $1$ respectively. 
Hence the total contribution of the section $\Psi_{\PP \D_7}$ from 
$\mathcal{B} \cap \mu$ is 
%all the 
%components of type $\B (q_{i_1}, l_{q_{\delta +1 }})$ equals
\bgd
2\binom{\delta}{1}N(\A_1^{\delta-1}\PP\D_8,n_1,m_1,m_2,\theta).
%+
%\binom{\delta}{1}N(\A_1^{\delta-1}\PP\E_7,n_1,m_1,m_2,\theta).
\edd 
That gives us Theorem \ref{a1_delta_pd7}. \\ 
\hf \hf Next, let us look at the section $\Psi_{\PP \E_7}$. 
We show that the contribution of the section $\Psi_{\PP \E_7}$  
from the points of 
$\overline{\A_1^{\delta-1} \circ \PP \X}_9 \cap \mu$ 
is $3$. Furthermore, the section does not vanish on 
$\A_1^{\delta-1}\circ \PP D_8$ (and hence does not vanish on 
$\overline{\A_1^{\delta-1}\circ \PP D}_8 \cap \mu$ since the 
cycle $\mu$ is generic). Hence, the total contribution of the 
section from $\mathcal{B} \cap \mu$ is 
\bgd
3\binom{\delta}{1}N(\A_1^{\delta-1}\PP\X_9,n_1,m_1,m_2,\theta).
%+
%\binom{\delta}{1}N(\A_1^{\delta-1}\PP\E_7,n_1,m_1,m_2,\theta).
\edd 
That proves Theorem \ref{a1_delta_pd7}. \qed

\subsection{Proof of Theorem \ref{a1_delta_pa7}: computing $\N(\A_1^{\delta} \PP \A_7, n_1, m_1, m_2, \theta)$}
\hf \hf Recall that in section \ref{notation}, 
we have defined the space  
\begin{align*}
\A^{\delta}_1 \circ \overline{\PP \A}_6 := \{ ([f], q_1, \ldots, q_{\delta}, 
l_{q_{\delta+1}}) \in \DD \times (X)^{\delta} \times \mathbb{P} T \X: 
&\textnormal{$f$ has a singularity of type $\A_1$ at $q_1, \ldots, q_{\delta}$}, \\ 
                   & ([f], l_{q_{\delta+1}}) \in \overline{\PP \A}_6, ~~\textnormal{$q_1, \ldots, q_{\delta+1}$ all distinct}\}. 
\end{align*}
Let $\mu$ be a generic pseudocycle representing the homology class 
Poincar\'{e} dual 
to 
\[c_1^{n_1} x_1^{m_1} x_2^{m_2} \lambda^{\theta}y^{\delta_L - (n_1 + m_1 + 2 m_2 +\theta+\delta+7)}.\]
As before, we let $v \in \hat{\gamma}$ and $w \in \pi^*TX/\hat{\gamma}$ be two fixed non zero vectors 
and define $f_{ij}$ as before. 
We now define section of the following bundle 
\begin{align}
\Psi_{\PP \A_7}:  \A^{\delta}_1 \circ \overline{\PP \A}_6 
\lra \mathbb{L}_{\PP \A_7} & := 
\hat{\gamma}^{\ast 7}\otimes\gamma_{\D}^{\ast 5}\otimes(\pi^*TX/\hat{\gamma})^{\ast 8} \otimes L^{\ast  5},  
\qquad \textnormal{given by} \nonumber \\
\{\Psi_{\PP \A_7}([f], q_1,\ldots,q_\delta, l_{q_{\delta+1}})\}(f^{\otimes 5}\otimes 
v^{\otimes 7}\otimes w^{\otimes 8}) &:= f_{02}^4 A_7^f.  \nonumber
%\psi_{\A_1} : \psi_{\A_0}^{-1}(0) \lra \mathcal{V}_{\A_1} &:= \gamma_{\D}^*\otimes T^*\X \otimes L, 
%\qquad \textnormal{given by} \qquad 
%\{\psi_{\A_1}([f], q)\}(f) := \nabla f|_q. \label{psi_a0_a1_section_defn}
\end{align}
%In \cite{RM_Hypersurfaces} we show that 
In 
%Lemma $7.1$ statement $*$ of the arXiv version 
\cite[Lemma 7.1, statement 12]{BM13}, 
we prove that 
%As stated earlier ()
\begin{align}
\overline{\PP A}_6 & = \PP A_6 \cup \overline{\PP A}_7 \cup \overline{\PP D}_7 \cup 
\overline{\PP E}_7 \cup \overline{\hat{X}}_9.  \label{pa6_closure}
\end{align}
Hence, let us define 
\begin{align*}
\mathcal{B} &:= \overline{A_1^{\delta} \circ \overline{\PP A}}_6 - \A^{\delta}_1 \circ (\PP\A_6\cup \overline{\PP A}_7). 
%\qquad \textnormal{and} \\
%\B^{\prime} &:= 
%\overline{A_1^{\delta} \circ \overline{\PP D}}_5 - \A^{\delta}_1 \circ (\PP\D_6\cup \overline{\PP E}_7)
\end{align*}
If $L$ is sufficiently $(2\delta+7)$-ample 
%and $(2\delta+4)$-ample, 
then the section 
%$\Psi_{\PP \A_7}$ 
%and $\Psi_{\PP \E_7}$ 
vanishes on the points of 
$\A_1^{\delta}\circ \PP A_7$  
%and $\A_1^{\delta}\circ \PP E_7$ 
transversally  
(all the points of $\A_1^{\delta}\circ \PP A_7$ are smooth points of the 
variety $\A_1^{\delta}\circ \overline{\PP A}_6$).
Hence 
\begin{align}
\lan e(\mathbb{L}_{\PP \A_7}), 
~~[\overline{\A^{\delta}_1 \circ \PP\A_6}] \cap [\mu] \ran & = \N(\A_1^{\delta}\PP \A_7, n_1, m_1, m_2, \theta)  
+ \mathcal{C}_{\mathcal{B}\cap \mu}. 
\end{align}
\hf \hf Let us define $\B_0$ as before, 
namely the component of the boundary 
all the points are distinct.
%where none of the 
%$(q_1, \ldots, q_{\delta})$ points are equal to the last marked point $q_{\delta+1}$. 
%In Lemma $7.1$ statement $*$ of the arXiv version 
In \cite[Corollary 7.11, Corollary 7.13 and Equation 7.62]{BM13}, 
we prove that 
the section $\Psi_{\PP A_7}$ vanishes on $\overline{A^{\delta}_1} \circ \PP D_7 \cap \mu$, 
$\overline{A^{\delta}_1} \circ \PP E_7 \cap \mu$ 
and $\overline{A^{\delta}_1} \circ \hat{X}_9 \cap \mu$ with a multiplicity of 
$6$, $7$ and 
$5$ respectively. 
Thus, the contribution of $\mathcal{B}_0 \cap \mu$ to the Euler class is 
\begin{align*}
6\Num(\A_1^{\delta}\PP\D_7, n_1, m_1, m_2, \theta)+ 7\Num(\A_1^{\delta}\PP\E_7, n_1, m_1, m_2, \theta) + 
5\Num(\A_1^{\delta}\hat{\X}_9, n_1, m_1, m_2, \theta).
\end{align*}
\hf \hf Next, let us analyze the boundary when 
one (or more) of the $q_i$ become equal to the last point $q_{\delta+1}$. 
We show in \cite{BM_closure_of_seven_points} that 
%First of all we note that the whole of 
%$\mathcal{B}$ is not relevant while computing the contribution to the Euler class; only the points 
%at which the section vanishes is relevant. Note that $\B$ contains $\overline{\A_1^{\delta}\circ \PP \D_7}\cup\overline{\A_1^{\delta}\circ \PP \E_7}$ and the section $\Psi_{\PP \A_7}$ vanish on it with multiplicity $6$ and $7$ respectively. Now we focus on the component of 
%$\mathcal{B}$ where one (or more) of the $q_i$ become equal to the last point $q_{\delta+1}$. 
%Let us define 
%\begin{align*}
%\B (q_1, l_{q_{\delta +1 }}) &:= \{ ([f], q_1, \ldots, q_{\delta}, l_{q_{\delta +1}}) \in \B: q_1 = q_{\delta +1} \}
%\end{align*}
%The spaces $\B (q_2, q_{\delta +1 }), \ldots, \B (q_{\delta}, q_{\delta +1 })$ are defined 
%similarly. We claim that 
\begin{align*}
\B(q_1, l_{q_{\delta+1}})\cap \mu & = \Big(\overline{\PP A}_8\cap \mu \Big) \cup 
\Big(\overline{\PP E}_8 \cap \mu \Big).   
\end{align*}
We also show that the contribution to the Euler class from $\overline{\PP A}_8\cap \mu$ and 
$\overline{\PP E}_8\cap \mu$ is $2$ and $14$ respectively. 
Hence the total contribution from points of type $\B(q_{i_1}, l_{q_{\delta+1}})$ is 
\begin{align*}
2 \binom{\delta}{1}\Num(\A_1^{\delta}\PP\A_8, n_1, m_1, m_2, \theta) + 14 \binom{\delta}{1}
\Num(\A_1^{\delta} \PP \E_8, n_1, m_1, m_2, \theta). 
\end{align*}
This gives us Theorem \ref{a1_delta_pa7}. \qed \\
%Hence the remaining stratum of $\B$ does not contribute, since its intersection with $\mu$ 
%is empty. \qed 
\hf \hf The rest of the computations involve enumerating curves with just 
one codimension eight singular point.  Hence, there will be no reference to a 
$\delta$, since $\delta$ is going to be zero. 
%of codimension , but in codimension eight. 

\subsection{Proof of Theorem \ref{a1_delta_pd8}: computing $\N(\PP \D_8, n_1, m_1, m_2, \theta)$}
%The computation of 
%$\N(\PP \D_8, n_1, m_1, m_2, \theta)$ and $\N(\PP \E_8, n_1, m_1, m_2, \theta)$ 
%are over the same spaces; hence it is more efficient to prove these to computations simultaneously. 
%Recall that we have defined the space  
%Let us define the space  
%\begin{align*}
%\PP \D_7 := \{ ([f], l_{q}) \in \DD \times \mathbb{P} T \X: 
%&\textnormal{$f$ has a singularity of type $\A_1$ at $q_1, \ldots, q_{\delta}$}, \\ 
                   %& ([f], l_{q_{\delta+1}}) \in \PP \D_7\} 
                   %, ~~\textnormal{$q_1, \ldots, q_{\delta+1}$ all distinct}\}. 
%\end{align*}
\hf \hf Let $\mu$ be a generic pseudocycle representing the homology class 
Poincar\'{e} dual 
to 
\[c_1^{n_1} x_1^{m_1} x_2^{m_2} \lambda^{\theta}y^{\delta_L - (n_1 + m_1 + 2 m_2 +\theta+8)}.\]
As before, we let $v \in \hat{\gamma}$ and $w \in \pi^*TX/\hat{\gamma}$ be two fixed non zero vectors 
and define $f_{ij}$ as before. 
We now define a section of the following bundles 
\begin{align}
\Psi_{\PP \D_8}:  \overline{\PP \D}_7 
\lra \mathbb{L}_{\PP \D_8} & :=  
\gamma_{\D}^{\ast 3}\otimes \hat{\gamma}^{\ast 9}\otimes (\pi^* TX/\hat{\gamma})^{\ast 6} \otimes L^{\ast 4},  
\qquad \textnormal{given by} \nonumber \\
\{\Psi_{\PP \D_8}([f], q_1,\ldots,q_\delta, l_{q_{\delta+1}})\}
(f^{\otimes 2}\otimes v^{\otimes 6}\otimes w^{\otimes 2}) &:= f_{12}^3 D_8^f,  \nonumber
%\psi_{\A_1} : \psi_{\A_0}^{-1}(0) \lra \mathcal{V}_{\A_1} &:= \gamma_{\D}^*\otimes T^*\X \otimes L, 
%\qquad \textnormal{given by} \qquad 
%\{\psi_{\A_1}([f], q)\}(f) := \nabla f|_q. \label{psi_a0_a1_section_defn}
\end{align}
%and 
%\begin{align}
%\Psi_{\PP \E_8}:  \PP \D_7 
%\lra \mathbb{L}_{\PP \E_8} & := \hat{\gamma}^{\ast 4}\otimes\gamma_{\D}^{\ast}\otimes L , \\
%\{\Psi_{\PP \E_7}([f], q_1,\ldots,q_\delta, l_{q_{\delta+1}})\}(f\otimes v^{\otimes 4}) &:= f_{40},  \nonumber
%\psi_{\A_1} : \psi_{\A_0}^{-1}(0) \lra \mathcal{V}_{\A_1} &:= \gamma_{\D}^*\otimes T^*\X \otimes L, 
%\qquad \textnormal{given by} \qquad 
%\{\psi_{\A_1}([f], q)\}(f) := \nabla f|_q. \label{psi_a0_a1_section_defn}
%\end{align}
where
\[ \D^{\q}_8 := \q_{60} + 
\frac{5 \q_{03} \q_{31} \q_{50}}{3 \q_{12}^2} 
-\frac{5 \q_{31} \q_{41}}{\q_{12}} - \frac{10 \q_{03} \q_{31}^3}{3 \q_{12}^3} 
+ \frac{5 \q_{22} \q_{31}^2}{\q_{12}^2}.
 \]
%In \cite{RM_Hypersurfaces} we show that 
If $L$ is sufficiently $(2\delta+6)$-ample, then this section is transverse to the zero set. 
We show in \cite{BM_closure_of_seven_points} that 
\begin{align*}
\overline{\PP D}_7 & = \PP D_7 \cup \overline{\PP D}_8 \cup \overline{\PP E}_8.  
\end{align*}
Hence, let us define 
\begin{align*}
\mathcal{B} &:= \overline{\PP \D}_7- \PP\D_7. 
\end{align*}
Hence 
\begin{align*}
\lan e(\mathbb{L}_{\PP \D_8}), 
~~[\overline{\A^{\delta}_1 \circ \PP\D_7}] \cap [\mu] \ran & = \N(\A_1^{\delta}\PP \D_8, n_1, m_1, m_2, \theta)  
+ \mathcal{C}_{\mathcal{B} \cap \mu}. 
%\label{na1_pseudocycle}
\end{align*}
We also show in \cite{BM_closure_of_seven_points} that 
the contribution of $\overline{\PP E}_8 \cap \mu$ to the Euler class is $3$. 
That gives us Theorem \ref{a1_delta_pd8}. \qed  
%and 
%\begin{align}
%\lan e(\mathbb{L}_{\PP \E_7}), 
%~~[\overline{\A^{\delta}_1 \circ \PP\E_6}] \cap [\mu] \ran & = \N(\A_1^{\delta}\PP \E_7, n_1, m_1, m_2, \theta)  
%+ \mathcal{C}_{\mathcal{B} \cap \mu}(\Psi_{\PP \E_7}).
%\label{na1_pseudocycle}
%\end{align}
%We now give an explicit description of $\mathcal{B}$. 
%We claim that 
%\[ \mathcal{B} = \overline{\PP D}_8 \cup \overline{\PP E}_8.  \]
%The contribution from $\overline{\PP E}_8 \cap \mu$ is $3$. 
%That gives us Theorem \ref{a1_delta_pd8}. \\ 

\subsection{Proof of Theorem \ref{a1_delta_pe8} and \ref{a1_delta_px9}: 
computing $\N(\PP \E_8, n_1, m_1, m_2, \theta)$ and $\N(\PP \X_9, n_1, m_1, m_2, \theta)$} 
\hf \hf Since the computation of both these numbers is done on top of the same space, it is 
more efficient to compute them together. 
Let $\mu$ be a generic pseudocycle representing the homology class 
Poincar\'{e} dual 
to 
\[c_1^{n_1} x_1^{m_1} x_2^{m_2} \lambda^{\theta}y^{\delta_L - (n_1 + m_1 + 2 m_2 +\theta+8)}.\]
As before, we let $v \in \hat{\gamma}$ and $w \in \pi^*TX/\hat{\gamma}$ be two fixed non zero vectors 
and define $f_{ij}$ as before. 
We now define sections of the following bundles  
\begin{align}
\Psi_{\PP \E_8}:  \overline{\PP \E}_7 
\lra \mathbb{L}_{\PP \E_8} & := 
\gamma_{\D}^{\ast} \otimes \hat{\gamma}^{\ast 3}\otimes (\pi^*TX/\hat{\gamma})^{\ast} \otimes L^{\ast} 
\qquad \textnormal{given by} \nonumber \\
\{\Psi_{\PP \E_8}([f], q_1,\ldots,q_\delta, l_{q_{\delta+1}})\}
(f \otimes v^{\otimes 3}\otimes w) &:= f_{31}  \nonumber \qquad \textnormal{and} \\ 
%\psi_{\A_1} : \psi_{\A_0}^{-1}(0) \lra \mathcal{V}_{\A_1} &:= \gamma_{\D}^*\otimes T^*\X \otimes L, 
%\qquad \textnormal{given by} \qquad 
%\{\psi_{\A_1}([f], q)\}(f) := \nabla f|_q. \label{psi_a0_a1_section_defn}
\Psi_{\PP X_9}:  \overline{\PP \E}_7 
\lra \mathbb{L}_{\PP X_9} & :=  
\gamma_{\D}^* \otimes (\pi^*TX/\hat{\gamma})^{\ast 3}\otimes L^{\ast} 
\qquad \textnormal{given by} \nonumber \\
\{\Psi_{\PP X_9}([f], q_1,\ldots,q_\delta, l_{q_{\delta+1}})\}
(f \otimes w^{\otimes 3}) &:= f_{03}.  \nonumber 
%\nonumber \qquad \textnormal{and} \\ 
\end{align}
%and 
%\begin{align}
%\Psi_{\PP \E_8}:  \PP \D_7 
%\lra \mathbb{L}_{\PP \E_8} & := \hat{\gamma}^{\ast 4}\otimes\gamma_{\D}^{\ast}\otimes L , \\
%\{\Psi_{\PP \E_7}([f], q_1,\ldots,q_\delta, l_{q_{\delta+1}})\}(f\otimes v^{\otimes 4}) &:= f_{40},  \nonumber
%\psi_{\A_1} : \psi_{\A_0}^{-1}(0) \lra \mathcal{V}_{\A_1} &:= \gamma_{\D}^*\otimes T^*\X \otimes L, 
%\qquad \textnormal{given by} \qquad 
%\{\psi_{\A_1}([f], q)\}(f) := \nabla f|_q. \label{psi_a0_a1_section_defn}
%\end{align}
%In \cite{RM_Hypersurfaces} we show that 
If $L$ is sufficiently $(2\delta+4)$-ample, then these sections are transverse to the zero set. 
We show in \cite{BM_closure_of_seven_points} that 
\begin{align}
\overline{\PP E}_7 & = \PP E_7 \cup \overline{\PP E}_8 \cup \overline{\PP X}_9. 
\end{align}
The directional derivative $f_{31}$ does not vanish on $\PP X_9$ after we intersect it with $\mu$. 
Similarly, the directional derivative $f_{03}$ does not vanish on $\PP E_8$. 
Hence, 
%Hence, let us define 
%\begin{align*}
%\mathcal{B} &:= \overline{\PP \E}_7- \PP\E_7. 
%\end{align*}
%Hence 
\begin{align*}
\lan e(\mathbb{L}_{\PP \E_8}), 
~~[\overline{\PP\E_7}] \cap [\mu] \ran & = \N(\PP \E_8, n_1, m_1, m_2, \theta) 
\qquad \textnormal{and} \\ 
\lan e(\mathbb{L}_{\PP X_9}), 
~~[\overline{\PP\E_7}] \cap [\mu] \ran & = \N(\PP X_9, n_1, m_1, m_2, \theta). 
%\qquad \textnormal{and} \\ 
%\label{na1_pseudocycle}
\end{align*}
%and 
%\begin{align}
%\lan e(\mathbb{L}_{\PP \E_7}), 
%~~[\overline{\A^{\delta}_1 \circ \PP\E_6}] \cap [\mu] \ran & = \N(\A_1^{\delta}\PP \E_7, n_1, m_1, m_2, \theta)  
%+ \mathcal{C}_{\mathcal{B} \cap \mu}(\Psi_{\PP \E_7}).
%\label{na1_pseudocycle}
%\end{align}
%We now give an explicit description of $\mathcal{B}$. We claim that 
%\[ \mathcal{B} = \overline{\PP D}_8 \cup \overline{\PP E}_8.  \]
%The contribution from $\overline{\PP E}_8 \cap \mu$ is $3$. 
That gives us Theorem \ref{a1_delta_pe8} and \ref{a1_delta_px9}. \qed \\

%\subsection{Proof of Theorem \ref{a1_delta_px9}: computing $\N(\A_1^{\delta} \PP \X_9, n_1, m_1, m_2, \theta)$}
%Recall 

\subsection{Proof of Theorem \ref{a1_delta_px9vee}: computing $\N(\PP \X_9^{\vee}, n_1, m_1, m_2, \theta)$}
\hf \hf Let $\mu$ be a generic pseudocycle representing the homology class 
Poincar\'{e} dual 
to 
\[c_1^{n_1} x_1^{m_1} x_2^{m_2} \lambda^{\theta}y^{\delta_L - (n_1 + m_1 + 2 m_2 +\theta+8)}.\]
As before, we let $v \in \hat{\gamma}$ and $w \in \pi^*TX/\hat{\gamma}$ be two fixed non zero vectors 
and define $f_{ij}$ as before. 
We now define sections of the following bundle 
\begin{align*}
\Psi_{\PP X_9^{\vee}}: \overline{\hat{X}}_9 \lra \mathbb{L}_{\PP X_9^{\vee}} & :=  
\gamma_{\D}^{\ast 3 } \otimes \hat{\gamma}^{\ast 9}\otimes (\pi^*TX/\hat{\gamma})^{\ast 3} \otimes L^{\ast 3}, 
\qquad \textnormal{given by}\\ 
\{\Psi_{\PP X_9^{\vee}}\}(v^{\otimes 9} \otimes w^{\otimes 3}) & :=  
\Big(- \frac{f_{31}^3}{8 } + \frac{3 f_{22} f_{31} f_{40}}{16 } - \frac{f_{13} f_{40}^2}{16}  \Big).
\end{align*}
We show in \cite{BM_closure_of_seven_points}, that there is no extra boundary contribution 
while computing
$\N(\PP \X_9^{\vee}, n_1, m_1, m_2, \theta)$. Hence, 
\begin{align*}
\lan e(\mathbb{L}_{\PP \X_9^{\vee}}), 
~~[\overline{\hat{X}}_9] \cap [\mu] \ran & = \N(\PP X_9^{\vee}, n_1, m_1, m_2, \theta).  
%\qquad \textnormal{and} \\ 
%\lan e(\mathbb{L}_{\PP X_9}), 
%~~[\overline{\PP\E_7}] \cap [\mu] \ran & = \N(\PP X_9, n_1, m_1, m_2, \theta). 
%\qquad \textnormal{and} \\ 
%\label{na1_pseudocycle}
\end{align*}
This gives us Theorem $\ref{a1_delta_px9vee}$. \qed

\subsection{Proof of Theorem \ref{a1_delta_pa8}: computing $\N(\PP \A_8, n_1, m_1, m_2, \theta)$}
\hf \hf Finally, we are ready to compute $\N(\PP \A_8, n_1, m_1, m_2, \theta)$. 
Let $\mu$ be a generic pseudocycle representing the homology class 
Poincar\'{e} dual 
to 
\[c_1^{n_1} x_1^{m_1} x_2^{m_2} \lambda^{\theta}y^{\delta_L - (n_1 + m_1 + 2 m_2 +\theta+8)}.\]
As before, we let $v \in \hat{\gamma}$ and $w \in \pi^*TX/\hat{\gamma}$ be two fixed non zero vectors 
and define $f_{ij}$ as before. 
We now define sections of the following bundle 
\begin{align*}
\Psi_{\PP A_8}: \overline{\PP A}_7 \lra \mathbb{L}_{\PP A_8} & :=  
\gamma_{\D}^{\ast 6 } \otimes \hat{\gamma}^{\ast 8}\otimes (\pi^*TX/\hat{\gamma})^{\ast 10} \otimes L^{\ast 6}, 
\qquad \textnormal{given by} \\ 
\{ \Psi_{\PP A_8}([f], l_q) \} (f^{\otimes 6} \otimes v^{8} \otimes w^{\otimes 10}) &:= f_{02}^5 A^f_8,
\end{align*}
where 
\begin{align*}
\A^{\q}_8 & := \q_{80} -\frac{28 \q_{21} \q_{61}}{\q_{02}}
-\frac{56 \q_{31} \q_{51}}{\q_{02}} + \frac{210 \q_{21}^2 \q_{42}}{\q_{02}^2} 
+\frac{420 \q_{21} \q_{22} \q_{41}}{\q_{02}^2}-\frac{210 \q_{03} \q_{21}^2
   \q_{41}}{\q_{02}^3} +\frac{560 \q_{21} \q_{31} \q_{32}}{\q_{02}^2}- \nonumber \\ 
   & \frac{840 \q_{13} \q_{21}^2 \q_{31}}{\q_{02}^3}- \frac{ 420 \q_{21}^3 \q_{23}}{\q_{02}^3} 
   +\frac{1260 \q_{03} \q_{21}^3 \q_{22}}{\q_{04}^4}-\frac{35 \q_{41}^2}{\q_{02}}+ 
   \frac{280 \q_{22}
   \q_{31}^2}{\q_{02}^2}-\frac{280 \q_{03} \q_{21} \q_{31}^2}{\q_{02}^3}
   -\frac{1260 \q_{21}^2 \q_{22}^2}{\q_{02}^3}+  \nonumber \\ 
   & \frac{105 \q_{04}
   \q_{21}^4}{\q_{02}^4}-\frac{315 \q_{03}^2 \q_{21}^4}{\q_{02}^5}
   +\frac{168 \q_{21} \q_{51} \q_{12}}{\q_{02}^2}+\frac{280 \q_{31} \q_{41}
   \q_{12}}{\q_{02}^2}-\frac{1680 \q_{21}^2 \q_{32} \q_{12}}{\q_{02}^3}
   -\frac{3360 \q_{21} \q_{22} \q_{31} \q_{12}}{\q_{02}^3} + \nonumber \\
   & \frac{2520 \q_{03} \q_{21}^2 \q_{31} \q_{12}}{\q_{02}^4}+  
   \frac{2520 \q_{13} \q_{21}^3 \q_{12}}{\q_{02}^4} 
    -\frac{840 \q_{21} \q_{41}
   \q_{12}^2}{\q_{02}^3} +\frac{7560 \q_{21}^2 \q_{22} \q_{12}^2}{\q_{02}^4}
   -\frac{560 \q_{31}^2 \q_{12}^2}{\q_{02}^3}-\frac{5040 \q_{03} \q_{21}^3 \q_{12}^2}{
   \q_{02}^5} \nonumber \\ 
   & +\frac{3360 \q_{21} \q_{31} \q_{12}^3}{\q_{02}^4} 
   -\frac{5040 \q_{21}^2 \q_{12}^4}{\q_{02}^5}.
\end{align*}
We show in \cite{BM_closure_of_seven_points} that 
\begin{align*}
\overline{\PP A}_7 & = \PP A_7 \cup \overline{\PP A}_8 \cup \overline{\PP D}_8 \cup \overline{\PP E}_8 \cup 
\overline{\PP X^{\vee}_9}.   
\end{align*}
We also show in \cite{BM_closure_of_seven_points} that the contributions from 
$\PP D_8 \cap \mu$, $\PP E_8 \cap \mu$ and $\PP X_9^{\vee} \cap \mu$ are 
$8$, $16$ and $6$. Hence, 
\begin{align*}
\lan e(\mathbb{L}_{\PP A_8}), 
~~[\overline{\PP A}_7] \cap [\mu] \ran & = \N(\PP A_8, n_1, m_1, m_2, \theta) \\ 
& + 8 \Num(\PP D_8, n_1, m_1, m_2, \theta) + 16 \Num(\PP E_8, n_1, m_1, m_2, \theta) + 
 6 \Num(\PP X_9^{\vee}, n_1, m_1, m_2, \theta).   
%\qquad \textnormal{and} \\ 
%\lan e(\mathbb{L}_{\PP X_9}), 
%~~[\overline{\PP\E_7}] \cap [\mu] \ran & = \N(\PP X_9, n_1, m_1, m_2, \theta). 
%\qquad \textnormal{and} \\ 
%\label{na1_pseudocycle}
\end{align*}
That gives us Theorem \ref{a1_delta_pa8}. \qed

\section{Low degree checks in $\mathbb{P}^2$}
\label{low_degree_checks}
%We now make some low degree checks. 
%\subsection{Curves in $\mathbb{P}^2$}
\subsection{Verifying the number $N(A_1^{8})$}  
\hf \hf We will now subject our formulas to some low degree checks. 
%We will begin by 
%verifying a special case of the formula of Kleiman and Piene for $8$-nodal curves. \\  
%considering curves in 
%$\mathbb{P}^2$. 
Let start by verifying Kleiman-Piene's formula for $8$-nodal curves for quintics in $\mathbb{P}^2$. 
We need to find out the number for degree $5$ curves in $\mathbb{P}^2$, through $12$ generic 
points that have $8$ nodes (let us say unordered). The genus of a smooth quintic is $6$; hence 
the curve has to break into at least three components. The possibilities are as follows: it could 
break into a nodal cubic through $8$ points and a pair of lines through the remaining four points. 
The total number of such configurations are 
\[ \binom{12}{8} \times 12 \times \binom{4}{2} \times \binom{2}{2} \times \frac{1}{2}. \] 
Note that the lines are indistinguishable; hence we divide by two. Next, the quintic 
could break into two conics and a line. The total number of such configurations is 
\[ \binom{12}{5} \times \binom{7}{5} \times \frac{1}{2}. \] 
Again, we divide by two since the conics are indistinguishable. 
Adding up, we get that the total number is 
\[ \binom{12}{8} \times 12 \times \binom{4}{2} \times \binom{2}{2} \times \frac{1}{2} 
+ \binom{12}{5} \times \binom{7}{5} \times \frac{1}{2} = 26136. \] 
%The sum total of these two 
%indeed gives us the answer 
This is precisely the value predicted by the formula of Kleiman and Piene!  

\subsection{Verifying the number $N(A_1^{6} A_2)$} 
\hf \hf Let us directly compute the number of quintics through $12$ points in $\mathbb{P}^2$ 
that have $6$ nodes and one cusp. This has to break into at least two components. The possibilities 
are: it could break as a cuspidal cubic going through $7$ points and a conic. The total number 
of such configurations is 
\[ \binom{12}{7}\times 24 \times \binom{5}{5}. \] 
Next, it could also break into a rational quartic with two nodes and one cusp through $10$ points 
and a line through the remaining two points. 
The total number of such configurations is 
\[ \binom{12}{10} \times 2304 \times \binom{2}{2}. \] 
Adding up, we get that the total number is  
\[ \binom{12}{7}\times 24 \times \binom{5}{5} + \binom{12}{10} \times 2304 \times \binom{2}{2} = 171072. \]
That is precisely the value predicted by our formula!  

\subsection{Verifying the number $N(A_1^4 D_4)$} 
\hf \hf Let us directly compute the number of quintics in $\mathbb{P}^2$ through 
$12$ points that have $4$ nodes and one triple point. This has to break into at least 
two components. The possibilities are as follows: it could break into a quartic through 
$10$ points having a triple point and a line through the remaining two points. The total 
number of such configurations are 
\[ \binom{12}{10} \times 60 \times \binom{2}{2}. \] 
Next, it could be an irreducible quartic through $11$ points, with three nodes and a line through the 
remaining point also going through any of the three nodes. The total number of such configurations is 
\[ \binom{12}{11} \times 620 \times 3. \] 
Next, we could put a line through any of the two points. And we could put an irreducible quartic 
through $10$ points, that has a node on the given line. First, lets us see how we can compute that 
last number. Using our formula for $N(A_1, n_1, m_1, m_2)$ we can compute that the number of quartics 
through $10$ points that have a three nodes and one of them on a given line; that 
number is $843$ 
(the free nodes are unordered in the final answer). 
Of these, we have to figure out how many 
are irreducible. There could be a cubic through $9$ points and a line through the 
last point and any of the three points the cubic intersects. Or, we could place a line through 
two points and a cubic through the remaining eight points and the point of intersection of the 
two lines. The total number of reducible configurations are 
\[ \binom{10}{9} \times 3 + \binom{10}{2} = 75.  \]
Hence, the number of irreducible cubics through $10$ points with three nodes, one of them on a 
line is $768$. Hence, the number of ways we could put a line through 
two points and put an irreducible quartic 
through $10$ points, that has a node on the given line is 
\[ \binom{12}{2} \times 768. \] 
Next, the curve could split as $3+2$. We could put a nodal cubic through 
$8$ points and a conic through the remaining $4$ points and the nodal point. That gives us 
\[ \binom{12}{8} \times 12. \] 
We could also put a conic through $5$ points and a cubic through $7$ points with a node on the 
given conic. The total number of such configurations is 
\[ \binom{12}{5} \times 6 \times 2. \] 
Note that the number of cubics through $7$ points with a node on a given line is $6$; hence 
the corresponding number with a node on a given conic is $12$.\\ 
\hf \hf Finally, the curve can split as $3+1+1$. We could put a cubic through $9$ points, 
a line through two points and then a line through the remaining point and any one of the three 
points of intersection. The total such configurations is 
\[ \binom{12}{9} \times \binom{3}{2} \times 3. \] 
Note that in this case the lines are distinguishable (one of them goes through two points, the 
other one goes through one point). Hence, we do not divide by two. Finally, we could put a 
line through two point, another line through two point and a cubic through the remaining eight points and the 
point of intersection. The total number of such configurations is 
\[ \binom{12}{2} \times \binom{10}{2} \times \frac{1}{2}.  \] 
Note that in this case, the two lines are indistinguishable, hence, we divide by two. 
Adding up all these numbers, we get that the desired number is 
\begin{align*}
& \binom{12}{10} \times 60 \times \binom{2}{2} +  
\binom{12}{11} \times 620 \times 3 \\
& + \binom{12}{2} \times 768 + \binom{12}{8} \times 12 
+ \binom{12}{5} \times 6 \times 2\\ 
& + \binom{12}{9} \times \binom{3}{2} \times 3 + 
\binom{12}{2} \times \binom{10}{2} \times \frac{1}{2}  = 95877. 
\end{align*}
This is precisely the number predicted by our formula!  

\subsection{Verifying the number $N(A_1^4 D_4, [L])$} 
\hf \hf Let us denote $N(A_1^4 D_4, [L])$ to the be the number of degree $d$ curves in $\mathbb{P}^2$ 
through the appropriate number of generic points that have $4$ nodes and one triple point lying on a 
line. Let us verify this number directly for quintics. We will directly compute the 
number of quintics through $11$ points that have four nodes and one triple point lying on a 
line. \\ 
\hf \hf The curve could break as $4+1$. It could break as quartic through $9$ points with a 
triple point on a line and a line through the remaining two points. Total number of such 
configurations is 
\[ \binom{11}{9} \times 12. \] 
Note that the number of quartics through $9$ points with a triple point on a line is $12$. 
Next, we could put an irreducible quartic through ten point with a node on the given line 
and put a line through the remaining eleventh point and the nodal point of the 
quartic that lies on the line. The total number of such configurations is 
\[ \binom{11}{10} \times 768.  \] 
Note that $768$ is the number of irreducible quartics through $10$ points 
with three nodes, one of them 
on a line. Next, let us compute the number of irreducible quartics through $9$ points, with three nodes 
one of them on a given point. The total number of such quartics is $105$ 
(which we can obtain by using our formula for $N(A_1^{\delta}, n_1, m_1, m_2)$). 
The reducible configuration is a cubic through any of the $8$ points and the given 
special point and a line through the remaining ninth point and the given special point. 
The total number of irreducible quartics through $9$ point, with three nodes, one of them 
lying on a point is 
\[ 105 - \binom{9}{8} = 96. \]
Hence, one final configuration that can occur is a line through two points and an irreducible 
quartic through nine points with a node lying on the given point of intersection of the two lines. 
That is 
\[ \binom{11}{2} \times 96. \] 
Similarly, the $3+2$ configurations are  as follows: put a conic through five points 
and then a nodal cubic through any one of the two points of intersection of the conic and the line. 
Or, put a nodal cubic through $7$ points, with a node on a given line and put a conic through the 
remaining five points. The sum total is 
\[ \binom{11}{7} \times 6+\binom{11}{5} \times 2.  \] 
Note that the number of cubic through seven points with a node on a given line is $6$.\\ 
\hf \hf Finally, the quintic could break as $3+1+1$. We could put a cubic through $9$ points 
and a line each through the remaining two connecting the the point of intersection of the 
cubic with the given line. The number of such configurations is 
\[ \binom{11}{9} \times 3 \times \binom{2}{1} \times \frac{1}{2}. \] 
Note that the lines are indistinguishable here. 
Finally, we can also put a line through two points, a line through one more of the 
remaining points and connecting the point of intersection of the previous line with the 
given line 
and a cubic through the remaining eight points and the nodal point. The number of such 
configurations is 
\[ \binom{11}{2} \times \binom{9}{1}.  \] 
Note that the lines are distinguishable here. 
Adding up, we get that the desried number is 
\begin{align*}
& \binom{11}{9} \times 12 + \binom{11}{10} \times 768 + \binom{11}{2} \times 96 + 
\binom{11}{7} \times 6+\binom{11}{5} \times 2  \\
& + \binom{11}{9} \times 3 \times \binom{2}{1} \times \frac{1}{2} + 
\binom{11}{2} \times \binom{9}{1} = 17952. 
\end{align*}
This is precisely the number predicted by our formula!  

\subsection{Verifying the number $N(A_1^4 D_4, [\textnormal{pt}])$} 
\hf \hf Let us denote $N(A_1^4 D_4, [\textnormal{pt}])$ to the be the number of degree $d$ curves in $\mathbb{P}^2$ 
through the appropriate number of generic points that have $4$ nodes and one triple point lying on a 
point. Let us verify this number directly for quintics. We will directly compute the 
number of quintics through $10$ points that have four nodes and one triple point lying on a 
point \\ 
\hf \hf The curve could break as $4+1$. It could break as quartic through $8$ points with a 
triple point on a point and a line through the remaining two points. Total number of such 
configurations is 
\[ \binom{10}{8} \times 1. \] 
Next, we could put an irreducible quartic through $9$ points  
with three nodes, one of them on a given point.The total number of such configurations is 
\[ \binom{10}{9} \times 96. \] 
Note that $96$ is the number of irreducible quartics through $9$ points with 
three nodes, one on a given point (as we calculated earlier). \\ 
\hf \hf Next, the curve could break as $3+2$. We could put a conic through four points and the 
given a point and then put a cubic through six points with node on the given point. The total 
number of such configurations are 
\[ \binom{10}{6}.  \] 
Finally, the curve could break as $3+1+1$. Put a cubic through any of the $8$ points and the 
given special point. Put a line through any one of the remaining point and the given special point. 
And put a line through the last point and the given special point. Total number of such configurations 
is 
\[ \binom{10}{8} \times \binom{2}{1} \times \binom{1}{1} \times \frac{1}{2}. \] 
Note that the lines are indistinguishable here, hence we divide by two. Adding up 
we get that the desired number is 
\begin{align*}
 \binom{10}{8} + \binom{10}{9} \times 96 + \binom{10}{4} + 
 \binom{10}{8}\times \binom{2}{1}\times \frac{1}{2} 
 & = 1260.
\end{align*}
That is precisely the number predicted by our formula!

\section{Low degree checks in $\mathbb{P}^1 \times \mathbb{P}^1$} 
\label{ldc_van}
\hf \hf In this section, we correct a few minor oversights made by Vainsencher in \cite{Van}, 
while he makes low degree checks for his formulas for curves in $\mathbb{P}^1 \times \mathbb{P}^1$. 
First, let us recall that that in 
\cite{KoMa}, Kontsevich and Manin obtain a formula for the number of rational 
degree $\beta$ curves in any del-Pezzo surface through $\langle c_1(TX), ~[\beta] \rangle -1$ 
generic points. 
%The fact that their formula is enumerative is proven in *. 
In particular, 
Kontsevich-Manin's formula can be applied to $\mathbb{P}^1 \times \mathbb{P}^1$. With this in mind, let 
us look at a couple of low degree checks made by Vainsencher in his paper \cite{Van}. \\ 
\hf \hf In \cite[Page 17]{Van}, 
%(Page $17$), 
Vainsencher, directly computes the number of irreducible 
curves in $\mathbb{P}^1 \times \mathbb{P}^1$ of bi-degree $(2,5)$ passing through $13$ points 
that have $4$-nodes. The answer he obtains is $3684$. However, the number of rational curves 
of bi-degree $(2,5)$ through $13$ points is predicted to be $3840$ by Kontsevich-Manin's formula. 
The reason for this mismatch 
is a simple combinatorial oversight made by the author. The total number of curves of bi-degree $(2,5)$ 
through $13$-points, with $4$-nodes is indeed $7038$. Of these curves, we have to subtract off all 
the reducible curves. One could have binodal curves of type $(2,4)$ through $12$ points and a curves 
of type $(0,1)$ through the remaining point. The number of bi-nodal curves of type $(2,4)$ through 
$12$ points is indeed $252$. However, not all of them are irreducible! There is the following 
configuration which occurs; a curve of type $(2,3)$ through $11$ points and a curve of type 
$(0,1)$ through one point. There are a total of $\binom{12}{1}$ such configurations. Hence, the 
total number of irreducible binodal curves of type $(2,4)$ is $252-12 = 240$. The total number 
of such configurations occurring for the $6$-nodal curves is 
\[ \binom{13}{1} \times 240. \] 
Finally, the curve can break as $(2,3) + (0,1) + (0,1)$. The two classes $(0,1)$ are 
\textit{indistinguishable}! Hence, the total number of such configurations is 
\[ \frac{1}{2}\times \binom{13}{11} \times \binom{2}{1} \times \binom{1}{1}.\] 
Hence, the final answer is 
\[ 7038- \binom{13}{1} \times 240 - \frac{1}{2}\binom{13}{11} \times \binom{2}{1} \times \binom{1}{1} = 3840.  \] 
This is precisely the value predicted by Kontsevich-Manin's formula! \\ 
\hf \hf Next, we look at another low degree check made by Vainsencher. In 
%page $18$ of 
\cite[Page 18]{Van}, 
he directly compute the number of irreducible curves of type $(3,4)$ through $13$ points 
that have $6$ nodes. The answer he obtains is $90508$. However, the number of rational curves 
of bi-degree $(3,4)$ passing through $13$ generic points is 
predicted by Kontsevich-Manins's formula to be $87544$. The reason for this 
mismatch is because the author omitted to consider one more configuration that occurs. \\ 
\hf \hf The number of $6$-nodal curves of class $(3,4)$ through $13$ points is indeed 
$122865$. To get the desired number, we simply have to subtract off all the reducible configurations. 
First, there will be nodal curves of class $(2,3)$ through $10$ points and a curve of class 
$(1,1)$ through $3$ points. The total number of such configurations is 
\[ \binom{13}{3} \times 20. \] 
Next, there can be curves of class $(2,2)$ through $8$ points and a curve of class $(1,2)$ through 
$5$ points. The total number of such configurations is 
\[ \binom{13}{8}. \] 
Next, the curve can break as $(3,2) + (0,1) + (0,1)$. The total number of such configurations is 
\[ \frac{1}{2} \times \binom{13}{11} \times \binom{2}{1}. \] 
Note that we divide by $2$ since the two $(0,1)$ curves are indistinguishable! 
Now we warn the reader, that we will not separately consider the configuration $(2,3) + (1,0) + (0,1)$; 
this will be part of a configuration we will soon consider. We don't consider this separately, since 
all the classes are different. \\ 
\hf \hf Next, we consider trinodal curves of class $(3,3)$ and one curve of class $(0,1)$. 
The number of trinodal curves of class $(3,3)$ through $12$ points is $1944$. This can break in 
two possible ways; either as $(3,2) + (0,1)$ or $(2,3) + (1,0)$. Let us only subtract off those 
configurations, where the curve breaks as $(3,2) + (0,1)$. The number of total configurations with 
$6$-nodal curves,  
excluding the ones of type $(3,2) + (0,1) + (0,1)$ is 
\[ \binom{13}{12} \times (1944- \binom{12}{11}). \] 
Note that the above number does include configurations of type $(2,3) + (1,0) + (0,1)$. This configuration 
doesn't require any special treatment, since $(1,0)$ and $(0,1)$ are two distinct classes. \\ 
\hf \hf Finally, we come to the configuration that was overlooked by Vainsencher. We can have a 
bi-nodal curve of type $(2,4)$ through $12$ points and a curve of class $(1,0)$ through one point. 
We have already calculated the number of irreducible bi-nodal curves of type $(2,4)$; it is 
$252-12 = 240$. Hence, the number of configuration of type $(2,4) + (1,0)$ where the $(2,4)$ is 
irreducible and binodal is 
\[ \binom{13}{12} \times 240. \] 
Hence, the final answer is 
\[ 122865 - \Big( \binom{13}{8} \times 20 + \binom{13}{8} + \frac{1}{2} \times \binom{13}{11} \times \binom{2}{1} + 
\binom{13}{12} \times (1944-12) + 240 \times \binom{13}{12}\Big) = 87544. \] 
This is precisely the value predicted by Kontsevich-Manin's formula!

\newpage  
\bibliographystyle{siam}
\bibliography{Myref_bib.bib}

\begin{thebibliography}{10}

\bibitem{BM_closure_of_seven_points}
{\sc S.~Basu and R.~Mukherjee}, {\em Collision of upto seven singular points}.
\newblock in preparation.

\bibitem{BM13}
\leavevmode\vrule height 2pt depth -1.6pt width 23pt, {\em Enumeration of
  curves with one singular point}.
\newblock available at \url{http://arxiv.org/abs/1308.2902}.

\bibitem{BM_Detail}
\leavevmode\vrule height 2pt depth -1.6pt width 23pt, {\em Enumeration of
  curves with singularities: Further details}.
\newblock available at \url{https://www.sites.google.com/site/ritwik371/home}.

\bibitem{BM_linear_system_codim_seven}
\leavevmode\vrule height 2pt depth -1.6pt width 23pt, {\em Enumeration of
  curves with upto two singular points in a general linear system}.
\newblock available at \url{http://arxiv.org/abs/1501.01557}.

\bibitem{BM13_2pt_published}
\leavevmode\vrule height 2pt depth -1.6pt width 23pt, {\em Enumeration of
  curves with two singular points}, Bull. Sci. Math,  (2014), pp.~667--735.

\bibitem{BM13_one_singular_point_published}
\leavevmode\vrule height 2pt depth -1.6pt width 23pt, {\em Enumeration of
  curves with one singular point}, Journ of Geom and Phys, 104 (2016),
  pp.~175--203.

\bibitem{BVSRM}
{\sc I.~Biswas, S.~D'Mello, R.~Mukherjee, and V.~Pingali}, {\em Rational
  cuspidal curves on del-pezzo surfaces}.
\newblock available at \url{http://arxiv.org/abs/1509.06300}.

\bibitem{BMT2}
{\sc I.~Biswas, R.~Mukherjee, and V.~Thakre}, {\em Genus two enumerative
  invariants in del-pezzo surfaces with a fixed complex structure}.
\newblock available at \url{http://arxiv.org/abs/1511.04900}.

\bibitem{BMT1}
\leavevmode\vrule height 2pt depth -1.6pt width 23pt, {\em Genus one
  enumerative invariants in del-pezzo surfaces with a fixed complex structure},
  Comptes Rendus Mathematique,  (2016).

\bibitem{BoTu}
{\sc R.~Bott and L.~W. Tu}, {\em Differential forms in algebraic topology},
  vol.~82 of Graduate Texts in Mathematics, Springer-Verlag, New York, 1982.

\bibitem{CH}
{\sc L.~Caporaso and J.~Harris}, {\em Counting plane curves of any genus},
  Invent. Math., 131 (1998), pp.~345--392.

\bibitem{F}
{\sc W.~Fulton}, {\em Intersection theory}, vol.~2 of Ergebnisse der Mathematik
  und ihrer Grenzgebiete. 3. Folge. A Series of Modern Surveys in Mathematics
  [Results in Mathematics and Related Areas. 3rd Series. A Series of Modern
  Surveys in Mathematics], Springer-Verlag, Berlin, second~ed., 1998.

\bibitem{Ionel_genus_one}
{\sc E.~Ionel}, {\em Genus one enumerative invariants in {$\Bbb P^2$}}, Duke
  Math Jour, 94 (1998), pp.~279--324.

\bibitem{IP_SSF}
{\sc E.~Ionel and T.~Parker}, {\em The {S}ymplectic {S}um {F}ormula for
  {G}romov-{W}itten invariants}, Annals of Mathematics, 159 (2004),
  pp.~935--1025.

\bibitem{Kaz}
{\sc M.~{\`E}. Kazarian}, {\em Multisingularities, cobordisms, and enumerative
  geometry}, Uspekhi Mat. Nauk, 58 (2003), pp.~29--88.

\bibitem{Kl_survey2}
{\sc S.~Kleiman}, {\em Counting curves \`{a} la {G}\"ottsche: A celebration of
  algebraic geometry}, Clay Mathematics Proceedings, Amer. Math. Soc.,
  Providence, RI, 2013, pp.~451--456.

\bibitem{KP}
{\sc S.~Kleiman and R.~Piene}, {\em Enumerating singular curves on surfaces},
  in Algebraic geometry: {H}irzebruch 70 ({W}arsaw, 1998), vol.~241 of Contemp.
  Math., Amer. Math. Soc., Providence, RI, 1999, pp.~209--238.

\bibitem{KlSh}
{\sc S.~Kleiman and V.~Shende}, {\em On the {G}\"ottsche threshold: A
  celebration of algebraic geometry}, Clay Mathematics Proceedings, Amer. Math.
  Soc., Providence, RI, 2013, pp.~429--449.

\bibitem{JKock}
{\sc J.~Kock}, {\em Characteristic number of rational curves with cusp or
  prescribed triple contact}, Math. Scand., 92 (2003), pp.~223--245.

\bibitem{KoMa}
{\sc M.~Kontsevich and Y.~Manin}, {\em Gromov-{W}itten classes, quantum
  cohomology, and enumerative geometry [ {MR}1291244 (95i:14049)]}, in Mirror
  symmetry, {II}, vol.~1 of AMS/IP Stud. Adv. Math., Amer. Math. Soc.,
  Providence, RI, 1997, pp.~607--653.

\bibitem{KST}
{\sc M.~Kool, V.~Shende, and R.~P. Thomas}, {\em A short proof of the
  {G}\"ottsche conjecture}, Geom. Topol., 15 (2011), pp.~397--406.

\bibitem{L2}
{\sc A.-K. Liu}, {\em The algebraic proof of the universality theorem}.
\newblock available at \url{http://arxiv.org/abs/math/0402045}.

\bibitem{L1}
\leavevmode\vrule height 2pt depth -1.6pt width 23pt, {\em Family blowup
  formula, admissible graphs and the enumeration of singular curves}, J.
  Differential Geom., 56 (2000), pp.~381--579.

\bibitem{GregMikh}
{\sc G.~Mikhalikin}, {\em Enumerative tropical algebraic geometry in}, J. Amer.
  Math. Soc., 18 (2005), pp.~313--377.

\bibitem{WMP3}
{\sc M.~Mikosz, P.~Pragacz, and A.~Weber}, {\em Positivity of {T}hom
  polynomials $ii$; the lagrange singularities}, Fund. Math., 202 (2009).

\bibitem{RM_thesis}
{\sc R.~Mukherjee}, {\em Enumerative geometry via topological computations}.
\newblock PhD Thesis, Stony Brook University, 2011 (available on request).

\bibitem{RM_Hypersurfaces_published}
\leavevmode\vrule height 2pt depth -1.6pt width 23pt, {\em Enumeration of
  singular hypersurfaces on arbitrary complex manifolds}, International Journal
  of Mathematics,  (2016).
\newblock In press.

\bibitem{Rahul_genus_one}
{\sc R.~Pandharipande}, {\em Counting elliptic plane curves with fixed
  j-invariant}, Proceedings of AMS, 125 (1997), pp.~3471--3479.

\bibitem{Rahul1}
{\sc R.~Pandharipande}, {\em Intersections of {$\Bbb Q$}-divisors on
  {K}ontsevich’s moduli space and {E}numerative {G}eometry}, Trans. Amer.
  Math. Soc., 351 (1999), pp.~1481--1505.

\bibitem{Rahul_Gottsche}
{\sc R.~Pandharipande and L.~{G}\"ottsche}, {\em The quantum cohomology of
  blowup of {$\Bbb P^2$} and enumerative geometry}, Jour of Diff Geometry, 48
  (1998), pp.~61--90.

\bibitem{WP1}
{\sc P.~Pragacz and A.~Weber}, {\em Positivity of {S}chur function expansions
  of {T}hom polynomials}, Fund. Math., 195 (2007), pp.~85--95.

\bibitem{WP2}
\leavevmode\vrule height 2pt depth -1.6pt width 23pt, {\em Thom polynomials of
  invariant cones, {S}chur functions and positivity}, Algebraic Cycles,
  Sheaves, Shtukas, and Moduli, Trends in Mathematics, Birkh¨auser,  (2007),
  pp.~117--129.

\bibitem{Ran2}
{\sc Z.~Ran}, {\em On nodal plane curves}, Invent. Math., 86 (1986),
  pp.~529--534.

\bibitem{Ran1}
{\sc Z.~Ran}, {\em Enumerative geometry of singular plane curves}, Invent.
  Math., 97 (1989), pp.~447--465.

\bibitem{Rennemo3}
{\sc J.~Rennemo}, {\em Universal polynomials for tautological integrals on
  {H}ilbert schemes}, Journal of Geometry and Topology,  (2016).
\newblock In press.

\bibitem{Rimanyi1}
{\sc R.~Rimanyi}, {\em Thom polynomials, symmetries and incidences of
  singularities}, Inventiones, 143 (2001), pp.~499--521.

\bibitem{Rimanyi2}
\leavevmode\vrule height 2pt depth -1.6pt width 23pt, {\em Multiple point
  formulas a new point of view}, Pacific Jour of Mathematics, 202 (2002),
  pp.~475--490.

\bibitem{Rimanyi_Szucs3}
{\sc R.~Rimanyi and A.~Szucs}, {\em Pontryagin {T}hom type construction for
  maps with singularities}, Topology, 37 (1998), pp.~1177--1191.

\bibitem{RT}
{\sc Y.~Ruan and G.~Tian}, {\em A mathematical theory of quantum cohomology},
  J. Differential Geom., 42 (1995), pp.~259--367.

\bibitem{MT_AZ}
{\sc M.~Tehrani and A.~Zinger}, {\em On {S}ymplectic {S}um {F}ormulas in
  {G}romov-{W}itten theory}.
\newblock available at \url{http://arxiv.org/abs/1404.1898}.

\bibitem{Tz}
{\sc Y.-J. Tzeng}, {\em A proof of the {G}\"ottsche-{Y}au-{Z}aslow formula}, J.
  Differential Geom., 90 (2012), pp.~439--472.

\bibitem{Tzeng_Li}
{\sc Y.-J. Tzeng and J.~Li}, {\em Universal polynoials for singular curves on
  surfaces}, Compos. Math., 150 (2014), pp.~1169--1182.

\bibitem{Van}
{\sc I.~Vainsencher}, {\em Enumeration of {$n$}-fold tangent hyperplanes to a
  surface}, J. Algebraic Geom., 4 (1995), pp.~503--526.

\bibitem{Vakil_CH}
{\sc R.~Vakil}, {\em Counting curves on rational surfaces}, Manuscripta
  Mathematica, 102 (2000), pp.~53--84.

\bibitem{A_Weber}
{\sc A.~Weber}, {\em Characteristic classes in singularity theory}.
\newblock Slides, available on
  \url{http://homepages.math.uic.edu/~jaca2009/notes/Weber.pdf}.

\bibitem{Zeuthen}
{\sc H.~Zeuthen}, {\em Almindelige egenskaber ved systemer af plane kurver},
  Kongelige Danske Videnskabernes Selskabs Skrifter, 10 (1873), pp.~285--393.

\bibitem{Z1}
{\sc A.~Zinger}, {\em Counting plane rational curves: old and new approaches}.
\newblock available at \url{http://arxiv.org/abs/math/0507105}.

\bibitem{Zinger_Thesis}
\leavevmode\vrule height 2pt depth -1.6pt width 23pt, {\em Enumerative
  algebraic geometry via techniques of symplectic topology and analysis of
  local obstructions}.
\newblock PhD Thesis, MIT, 2002 (available from the author on request).

\bibitem{g2p2and3}
{\sc A.~Zinger}, {\em Enumeration of genus-two curves with a fixed complex
  structure in {$\Bbb P^2$} and {$\Bbb P^3$}}, J. Differential Geom., 65
  (2003), pp.~341--467.

\bibitem{g1}
\leavevmode\vrule height 2pt depth -1.6pt width 23pt, {\em Enumeration of
  one-nodal rational curves in projective spaces}, Topology, 43 (2004),
  pp.~793--829.

\bibitem{g0pr}
\leavevmode\vrule height 2pt depth -1.6pt width 23pt, {\em Counting rational
  curves of arbitrary shape in projective spaces}, Geom. Topol., 9 (2005),
  pp.~571--697 (electronic).

\bibitem{g3}
\leavevmode\vrule height 2pt depth -1.6pt width 23pt, {\em Enumeration of
  genus-three plane curves with a fixed complex structure}, J. Algebraic Geom.,
  14 (2005), pp.~35--81.

\bibitem{pseudo_cycle}
\leavevmode\vrule height 2pt depth -1.6pt width 23pt, {\em Pseudocycles and
  integral homology}, Transactions of AMS, 360 (2008), pp.~2741--2765.

\end{thebibliography}

\vspace*{0.4cm}

\hf {\small D}{\scriptsize EPARTMENT OF }{\small M}{\scriptsize ATHEMATICS, }{\small IISER KOLKATA, }
{\small WB} {\footnotesize 741246, }{\small INDIA}\\
\hf{\it E-mail address} : \texttt{basu.somnath@gmail.com,  somnath.basu@iiserkol.ac.in}\\

\hf {\small S}{\scriptsize CHOOL OF }{\small M}{\scriptsize ATHEMATICS, }{\small NISER, BHUBANESWAR (HBNI),}
%{\small T}{\scriptsize ATA  }
%{\small I}{\scriptsize NSTITUTE }
%{\small OF}
%{\scriptsize FUNDAMENTAL}{\small R} {\scriptsize ESEARCH} 
{\small O}{\scriptsize DISHA }
{\footnotesize 752050, }{\small INDIA}\\
\hf{\it E-mail address} : \texttt{ritwik371@gmail.com, ritwikm@niser.ac.in}\\[0.2cm]

\end{document}